\documentclass{article}
\usepackage{graphicx} 
\usepackage[top=0.75in, bottom=0.75in, left=0.75in, right=0.75in]{geometry}

\usepackage{amsmath}
\usepackage{amsthm}
\usepackage{multirow}
\usepackage{booktabs}

\DeclareMathOperator*{\argmin}{arg\,min}
\DeclareMathOperator*{\sign}{sign}
\usepackage{amsfonts}
\usepackage[table,dvipsnames]{xcolor}
\usepackage{setspace}
\usepackage[titletoc]{appendix}
\usepackage{hyperref}
\newtheorem{proposition}{Proposition}[section]
\usepackage[round]{natbib}
\usepackage{subcaption}
\usepackage{algorithm}
\usepackage{algorithmic}
\usepackage{tikz}
\usetikzlibrary{arrows.meta}

\title{ParaDeAl: Parallel Descent Algorithms for Sparse Statistics}
\author{Adam Deng\\MIT ORC\\ \texttt{adamdeng@mit.edu} \and Rahul Mazumder\\MIT ORC\\ \texttt{rahulmaz@mit.edu}}
  
\begin{document}

\vspace{-10pt}
\maketitle

\begin{abstract}
    Best subset selection is a central problem in statistical regression, drawing increasing attention from the broader optimization and machine learning communities. Solving the hard L0 sparsity constraint directly has been computationally challenging; as a result, most existing work proposes formulations which do not enforce the constraint, but instead use a mix of L0/L1/L2 penalties to encourage sparsity. Modern methods for solving such penalized formulations are tractable for $p \approx 10^5$, but give up direct control over sparsity. In this paper, we propose a two-stage proximal gradient descent-based approach which both enforces the sparsity constraint and incorporates L1/L2 penalties, leveraging GPUs for computational acceleration to find approximate solutions. Our method, ParaDeAl, offers the user the ability to select the problem hyperparameters, solves all configurations in parallel, and obtains solutions up to $15 \times$ faster than the best state-of-the-art CPU competitors while achieving almost identical test accuracy for the same support size for large ($p > 10^5$) synthetic and real-world datasets.
    
\end{abstract}

\section{Introduction}

\subsection{Existing Work} \label{subsec:existing_work}

The proliferation of high-dimensional data has increased scientific interest in linear and sparse regression. We start with the standard linear regression model $y = X\beta + \epsilon$, where $X \in \mathbb{R}^{n \times p}$ is the data matrix of $n$ samples with $p$ features, $y \in \mathbb{R}^n$ is the response vector, $\beta \in \mathbb{R}^p$ is the regression coefficients vector, and $\epsilon \in \mathbb{R}^n$ is an errors vector. We assume that each column of $X$ has been standardized to mean 0 and L2-norm 1. For explainability and parsimony reasons, it is often desirable to sparsify $\beta$. One such formulation is the best subsets estimator (\cite{mil02}), which solves the problem \begin{equation} \min_\beta L(\beta) \stackrel{\text{def}}= ||y-X\beta||_2^2 \ \text{s.t.} \ ||\beta||_0 \leq k, \label{eq:bss} \end{equation} where the L0-norm $||\beta||_0$ counts the number of nonzeros in $\beta$. Exact solutions to Problem \ref{eq:bss} can be computed through a mixed integer optimization (MIO) approach in (\cite{ber15}), but doing so is compute-intensive. Near-optimal solutions can be found for $p \approx 10^3$ in several minutes, with significantly worse runtimes as $k$ increases. Alternatively, one can get approximate solutions to Problem \eqref{eq:bss} using iterative hard thresholding (IHT) (\cite{blu08}, \cite{ber15}), but IHT also suffers from computational inefficiency due to the large matrix multiplications in and convergence slowness of gradient descent (GD).

One related line of work focuses on convex proxies for Problem \ref{eq:bss} without the L0 constraint, including L1 aka Lasso (\cite{tibs96}) and  L1+L2 (Elastic Net) (\cite{zou05}) which is implemented in \verb|glmnet| (\cite{bre11}). Another line of work considers a penalized version for L0 with additional regularizers, such as L0+L1 and L0+L2 (\cite{haz20}), implemented in \verb|L0Learn| as approximate procedures that produce good solutions. Such procedures, relying on coordinate descent, are generally much faster than the IHT procedures discussed in \cite{ber15}. Additionally, \verb|L0BnB| (\cite{haz21}) proposes a tailored branch-and-bound procedure for solving the L0+L2 problem exactly. While these formulations are computationally appealing, they do not directly control the support size of $\beta$, which Problem \eqref{eq:bss} does.

In recent years, graphical processing units (GPUs) have gained popularity due to their ability to perform magnitudes more computations than CPUs in the same time (\cite{sou24}). GPUs are already deeply integrated into machine learning (ML) contexts, with tasks like training neural networks (\cite{cir11}) and deep learning for high-dimensional data (\cite{mit19}). Additionally, GPUs are used in optimization problems and methods, like linear programming (\cite{wir22}, \cite{lu24}, \cite{spa09}) and concurrent work on Branch and Bound (\cite{men26}). GPU-based methods have been proposed for the noiseless version of Problem \eqref{eq:bss}, compressed sensing (\cite{bla13}), which finds the minimum $||\beta||_0$ such that satisfies $y = X\beta$ exactly.

\subsection{Approach} \label{subsec:approach}

Motivated by earlier work (\cite{maz22}), we consider the estimator \begin{equation} \min_\beta f(\beta) \stackrel{\text{def}}= \left( L(\beta) + \lambda_1 ||\beta||_1 + \lambda_2 ||\beta||_2^2 \right) \ \text{s.t.} \ ||\beta||_0 \leq k. \label{eq:our_estimator} \end{equation} Here $\lambda_1 \geq 0$ and $\lambda_2 \geq 0$ are regularization parameters controlling the L1 and L2 norms of $\beta$ and can help in terms of out-of-sample prediction performance, while the constraint $||\beta||_0 \leq k$ controls the number of nonzero coefficients in $\beta$. We obtain feasible approximate solutions to Problem \eqref{eq:our_estimator} using a two-stage approach: support selection followed by coefficient value optimization, also known as value modification.  We use a proximal gradient descent-type procedure (\cite{ber15}) to compute a good support for Problem \eqref{eq:our_estimator}, which reduces to IHT when $\lambda_1, \lambda_2 = 0$. After finding the support, we modify the values of the nonzero coefficients of $\beta$ by solving Problem \eqref{eq:our_estimator} restricted to this support. In the special case $\lambda_1 = 0$, value fine-tuning is computed using preconditioned conjugate gradient (PCG) (\cite{she94}), which is generally faster than support-restricted GD. (We do not calculate the exact solution directly using inverses because it is difficult to parallelize over different supports and $\lambda_2$s, and computing the inverse is expensive.) When $p$ is large ($p \geq 10^4$), we also employ an active set-type procedure inspired by \cite{fri10} and \cite{haz20} to achieve computational speedups.

To accelerate computation we use GPUs, as our approach is well-suited for exploiting two types of parallelism offered by GPUs. The foremost type is matrix-matrix and matrix-vector products, forming the backbone of our approach. Moreover, GPUs can parallelize across multiple instances of $(k, \lambda_1, \lambda_2, \eta, \beta_o)$ values ($\beta_o$ means the starting $\beta$)---collectively, hyperparameters---solving Problem \eqref{eq:our_estimator} simultaneously for many hyperparameter configurations. In sparse regression, it is desirable to obtain solutions for multiple $(\lambda_1, \lambda_2, k)$---collectively, regularization hyperparameters---to assess the effects of different imposed conditions on the resulting $\beta$ and obtain more accurate $\beta$s for regression. Varying $\beta_o$ and learning rate $\eta$ is also important because for a given $(\lambda_1, \lambda_2, k)$, the quality of solution depends on $\beta_o$ and $\eta$, and selecting many $\beta_o$s and $\eta$s offers us flexibility to better optimize Problem \ref{eq:our_estimator}.

We term our collective methodology ParaDeAl for Parallel Descent Algorithms.

\subsection{Contributions}

Our contributions in this paper are as follows:

\begin{itemize}
    \item \textbf{Massively Parallel Sparse Regression:} We propose a GPU-compatible tensor-based framework, ParaDeAl, which obtains good solutions to Problem \eqref{eq:our_estimator} by parallelizing across the entire hyperparameter space $(\lambda_1, \lambda_2, k, \beta_o, \eta)$. Unlike existing packages, which are CPU-based and rely on procedures that can only solve problem instances sequentially or iteratively, we simultaneously solve up to thousands of problem candidates.
    \item \textbf{Scalability Beyond Existing Solvers:} Our method scales to solve 1000 simultaneous problem instances for $(n \approx 10^4, p \approx 10^5)$ or a single instance for $(n \approx 10^5, p \approx 10^6)$ on a single L40S GPU (48GB storage) within an hour. The latter problem is larger by one to two orders of magnitude than what existing solvers can solve. The former setup takes the best existing algorithms over 40x as long as ParaDeAl to solve.
\end{itemize}

\subsection{Paper Structure}

We present the full ParaDeAl architecture, with algorithmic pseudocode on specific parallelization procedures, in Section \ref{sec:architecture}. We then address a set of computational challenges and propose solutions in Section \ref{sec:considerations}. ParaDeAl's methodological edge over iterative counterparts is reported in Section \ref{sec:method_results}, while ParaDeAl's speed advantage and comparable accuracy vis-a-vis sparse regression solvers is demonstrated in Section \ref{sec:compare_results} through a series of tests involving synthetic and real-world datasets. We offer concluding remarks in \ref{sec:summary}. To provide completeness, and due to space limits, we explore certain implementation experiments and results in the appendices.

\subsection{Notation} \label{subsec:notation}

Let $\mathcal{S} = S(\beta)$ denote the support of $\beta$, i.e. the set of indices of $\beta$ for which the element value is not zero. 

For any set $\mathcal{T}$ which contains a subset of the $p$ features, let $P_\mathcal{T}(\beta)$ denote the projection of $\beta$ onto the elements of $\mathcal{T}$, so that $\beta_i = 0$ if $i \notin \mathcal{T}$. Let $R_\mathcal{T}(\beta)$ denote the compressed vector of length $|\mathcal{T}|$ consisting only of elements of $\beta$ whose features are within $\mathcal{T}$. Similarly, let $R_\mathcal{T}(X, \beta)$ denote a submatrix of $X$ consisting only of the columns corresponding to features of $\beta$ within $\mathcal{T}$.

We regard $\lambda_1, \lambda_2, k$ as regularization hyperparameters; and $\lambda_1, \lambda_2, k, \beta, \eta$ as hyperparameters collectively. We define one problem instance, with a particular $(\lambda_1, \lambda_2, k, \beta_o, \eta)$, as a configuration. Each configuration corresponds to its own $\beta$-trajectory; note that the number of configurations, which equals the number of $\beta$-trajectories, may exceed the number of different initial $\beta$ values chosen.

Let $\mathcal{C}$ represent the set of all hyperparameter configurations, where each configuration is a list of 5 hyperparameter values $(k, \lambda_1, \lambda_2, \eta, \beta)$. For a particular hyperparameter $a$, let $\mathcal{C}_a$ be the set of different values. Define $C_a$ as $|\mathcal{C}_a|$, i.e. the number of different values hyperparameter $a$ takes; similarly, define $\prod C$ as $|\mathcal{C}|$. An element of $\mathcal{C}$ is denoted as $c \in \mathcal{C}$; note that $c = (k, \lambda_1, \lambda_2, \eta, \beta)$ for some values of those hyperparameters.

\section{ParaDeAl Architecture} \label{sec:architecture}

In this section, we describe our algorithm for ParaDeAl, the motivation for parallelizing over problem configurations, and GPU-friendly ways to parallelize. In order to achieve parallelism, we may have to introduce redundancies in computation. We thus observe a GPU workflow tradeoff between uninterrupted vectorization and resizing data structures to minimize computation.

\subsection{Preliminaries} \label{subsec:preliminaries}

Consider our optimization formulation in Problem \eqref{eq:our_estimator}. Define \begin{equation} F_D(\beta) = L(\beta) + \lambda_2\|\beta\|_2^2 \label{eq:convex_part} \end{equation} as the convex differentiable portion of $F(\beta)$. We observe that $F_D$ is a $L_0$-smooth convex function, satisfying \begin{equation} \|\nabla F_D(\beta) - \nabla F_D(\alpha)\|_2 \leq L_0 \|\beta - \alpha\|_2 \ \ \forall \beta, \alpha \in \mathbb{R}^p \label{eq:lipschitz_inequality} \end{equation} where $L_0 = 2(\sigma_{\max}^2(X) + \lambda_2)$ and $\sigma_{\max}(X)$ refers to the largest singular value of $X$. Using the Descent Lemma (\cite{nes04}), for any $L' \geq L_0$, we obtain the following upper bound on $F_D(\beta)$: 
\begin{equation} F_D(\beta) \leq F_D(\alpha) + \left< \nabla F_D(\alpha), \beta-\alpha \right> + \frac{L'}{2} \|\beta-\alpha\|_2^2 := U_{L'}(\beta, \alpha) \ \ \forall \beta, \alpha \in \mathbb{R}^p. \label{eq:descent_lemma} \end{equation} Instead of directly minimizing Problem \ref{eq:our_estimator}, we minimize a local model of $F(\beta) = F_D(\beta) + \lambda_1||\beta||_1$ around candidate point $\alpha$. This leads to the following optimization problem: \begin{equation} \min_{\|\beta\|_0 \leq k} U_{L'}(\beta, \alpha) + \lambda_1 \|\beta\|_1 = \min_{\|\beta\|_0 \leq k} F_D(\alpha) + \frac{L'}{2} \left\| \beta - \left( \alpha - \frac{1}{L'} \nabla F_D(\alpha) \right) \right\| _2^2 - \frac{1}{2L'} \|\nabla F_D(\alpha)\|_2^2 + \lambda_1\|\beta\|_1 \label{eq:deriving_final_opt_problem} \end{equation} which is equivalent to solving \begin{equation} \min_{\|\beta\|_0 \leq k} \frac{L'}{2} \left\| \beta - \left( \alpha - \frac{1}{L'} \nabla F_D(\alpha) \right) \right\| _2^2 + \lambda_1\|\beta\|_1 \label{eq:final_opt_problem} \end{equation} because we only include optimization terms dependent on $\beta$. 

To solve~\eqref{eq:final_opt_problem} we need to consider the following thresholding operator:
\begin{equation} T(u, k, \lambda) := \argmin_{\beta \in \Re^p: \|\beta\|_0 \leq k} \frac{1}{2} \|\beta-u\|_2^2 + \lambda \|\beta\|_1\label{eq:t_thresholding} \end{equation} 
which can be computed in closed form as shown in the following proposition.

\begin{proposition} Let $T$ be the thresholding operator defined in~\eqref{eq:t_thresholding}. If $\hat{\beta}_{i}$ denotes the $i$th coordinate of $T$, then 
$$\hat{\beta}_i = \begin{cases}
0, & \text{if} \ i \notin \mathcal{R}_k(u) \\
\sign(u_i) \cdot \max \left( |u_i|-\lambda, \ 0\right), & \text{if} \ i \in \mathcal{R}_k(u)
\end{cases}$$ where $\mathcal{R}_k(u)$ denotes the set of $k$ indices corresponding to the $k$ largest absolute values of $u$. \label{prop:thresholding_form} \end{proposition} 

A proof is given in Appendix~\ref{app:proofs}. When there are ties among the magnitudes of $u$, Problem~\eqref{eq:t_thresholding} can have multiple minimizers, so $T$ is in general set-valued. Our algorithm is inspired by IHT \citep{blu08, ber15}.
Starting from an initial solution $\beta_{0} \in \mathbb{R}^p$, we apply an iterative thresholding procedure resulting in a sequence of solutions $\{ \beta_t \}$, given by \begin{equation} \beta_{t+1} \in T(\beta_t - \eta \nabla F_D(\beta_t), k, \eta \lambda_1), \label{eq:iterative_T} \end{equation} where the learning rate is $\eta = 1/L'$. 
The procedure is run until a stopping criterion is met; the criteria used in each stage are given in Section~\ref{subsec:two_stage_desc}.

By Proposition \ref{prop:thresholding_form}, we have the explicit form of~\eqref{eq:iterative_T}: 
\begin{equation} \label{eq:full_explicit_form} \left( \beta_{t+1} \right)_i = \begin{cases}
0, & \text{if} \ i \notin \mathcal{R}_k(\beta_t) \\
\sign(\beta_t - \eta \nabla F_D(\beta_t)) \cdot \max \left( |\beta_t - \eta \nabla F_D(\beta_t)|-\eta \lambda_1, \ 0\right), & \text{if} \ i \in \mathcal{R}_k(\beta_t)
\end{cases}. \end{equation} 

\begin{proposition} The updates given in \ref{eq:iterative_T} will converge to a particular $\beta$. \label{prop:converge} \end{proposition}

We prove Proposition \label{prop:converge} in Appendix \ref{app:proofs}.

We now describe how to parse the computation in Equation \eqref{eq:full_explicit_form} so that we can take advantage of GPU parallelism. We first perform the gradient update step \begin{equation} \tilde{\beta}_t \gets \beta_t - \eta \nabla F_D(\beta_t) \label{eq:1_grad_step} \end{equation} followed by the soft thresholding step (applied elementwise) \begin{equation} \tilde{\beta}_{t} \gets \text{soft}(\tilde{\beta}_t, \eta \lambda_1) := \text{sign}(\tilde{\beta}_t) \cdot \max(|\tilde{\beta}_t - \eta \lambda_1, 0) \label{eq:2_soft_threshold} \end{equation} followed by the hard thresholding operator \begin{equation} \beta_{t+1} = H_k(\tilde{\beta_{t}}) \label{eq:3_hard_threshold} \end{equation} where \begin{equation} H_k(\beta_{i}) = \begin{cases} 
0, & \text{if $i \notin \mathcal{R}_k(\beta_i)$} \\ \beta_{i}, & \text{if $i \in \mathcal{R}_k(\beta_i)$}, \end{cases} \label{explicit_beta_zeros} \end{equation} and $\mathcal{R}_k$ is defined in Proposition \ref{prop:thresholding_form}.

We refer to update~\eqref{eq:full_explicit_form} as a \textbf{hard-soft thresholding} update, as it consists of soft thresholding followed by hard thresholding. It updates the support and the coefficient values on the support at the same time, and every iterate $\beta_t$ satisfies the constraint $\|\beta_t\|_0 \leq k$. Moreover, when $\eta \leq 1/L_D$ the objective $F(\beta_t)$ is nonincreasing in $t$ (Proposition~\ref{prop:descent} in Appendix~\ref{app:proofs}).

Hard-soft thresholding updates the support and the values of the coefficients on the support. For a fixed support, Problem~\eqref{eq:our_estimator} is convex and can be solved with first-order methods (e.g., proximal gradient descent or PCG). This motivates the two-stage approach we describe here. To obtain good solutions to Problem~\eqref{eq:our_estimator}, we use two steps. (i) Support selection: we apply update~\eqref{eq:full_explicit_form} until a stopping criterion is met, and obtain a support. (ii) Value modification: we solve Problem~\eqref{eq:our_estimator} restricted to the support found in (i), which is a convex problem. In the special case $\lambda_1 = 0$ we use PCG, and when $\lambda_1 \neq 0$, we use a proximal gradient procedure (possibly with a different step-size than update~\eqref{eq:full_explicit_form}) as described in Section~\ref{subsubsec:value_mod}.

\subsection{ParaDeAl's Two-Stage Methodology for Problem~\eqref{eq:our_estimator}} \label{subsec:two_stage_desc}

In this section we explicitly write the details for getting solutions to Problem~\eqref{eq:our_estimator}, building on the preliminaries of Section~\ref{subsec:preliminaries}. 

We consider two separate cases: $\lambda_1 \neq 0$ and $\lambda_1 = 0$ as they have differences in implementation methodology. For each case we break down our discussion into support selection and value modification.

In support selection, we choose a \textbf{candidate set} $\mathcal{K}$ with $|\mathcal{K}| = k$ such that $\beta_{i} = 0 \ \text{if} \ i \notin \mathcal{K} \ \forall 1 \leq i \leq p$. We call $\mathcal{K}$ a candidate set rather than a support because it is possible that elements within $\mathcal{K}$ happen to be zero. The benefits of using a candidate set rather than defining $\mathcal{K}$ as the strict support with the possibility that $|\mathcal{K}| < k$ will be observed later, in value modification. We define our $\beta$ after support selection as $\beta_w$ (for ne\textbf{w} beta). After the candidate set is chosen, we perform value modification on it, obtaining

\begin{equation}
\text{Value Modification}~~~~~~~~\beta^* \in \argmin_\beta F(\beta) \ \ \text{s.t.} \ \ \beta_{i} = 0, i \notin \mathcal{K}. \label{eq:vm_equation} \end{equation} A diagram of the full ParaDeAl procedure, broken down into cases, is provided below in Figure \ref{fig:paradeal_arrow_diagram}. 

\begin{figure}[htbp]
    \centering
    \includegraphics[scale=0.35]{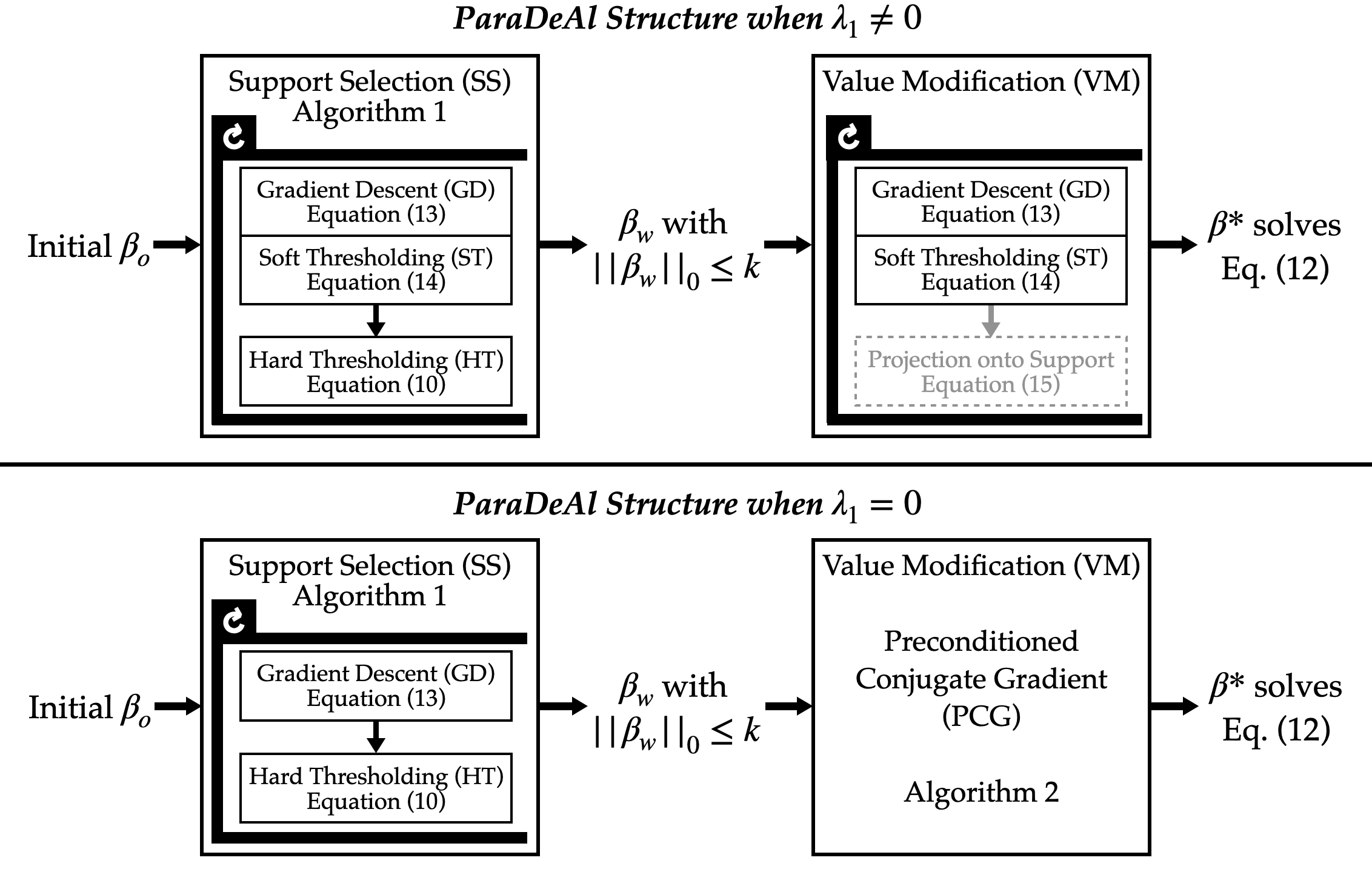}
    \caption{ParaDeAl structure. The gray-colored projection onto the candidate set $\mathcal{K}$, shown in value modification for $\lambda_1 \neq 0$, may or may not be performed depending on implementation choices described in Section \ref{subsec:parallelizing_vm}.}
    \label{fig:paradeal_arrow_diagram}
\end{figure}

\subsubsection{Support Selection} \label{subsubsec:support_selection}

The goal of support selection is to quickly find a high-quality support $\mathcal{K}$ that corresponds to an optimal or near-optimal objective for Problem \eqref{eq:our_estimator}. We split our method into two cases: If $\lambda_1 \neq 0$, we apply exactly the procedure described in Section~\ref{subsec:preliminaries}. If $\lambda_1 = 0$, we apply the same procedure except we remove the soft thresholding, as seen in Figure \ref{fig:paradeal_arrow_diagram}. 

We also adjust the stopping criterion. Support selection (Algorithm~\ref{alg:support_selection}) stops as soon as one update leaves the support unchanged, that is, once two consecutive iterates have the same support. The rule depends on the support only: the role of support selection is to find a support, and the coefficient values are set by value modification. More generally, one can require the support to remain unchanged for $r$ consecutive updates; the rule above is $r = 1$. Appendix~\ref{app:stability_r} examines the effect of $r$.

Once the support stabilizes, we output the final $\beta_w$ to be passed to the value modification function. Pseudocode for support selection is given below in Algorithm \ref{alg:support_selection}, with steps used only when $\lambda_1 \neq 0$ marked in \color{Cerulean}cerulean\color{black}. We also use an upper limit on the number of iterations (denoted by a pre-specified parameter $M_{SS}$). 

\paragraph{Returning the lowest-objective support.} \label{subsec:lowest_objective} With the large step sizes of Section~\ref{subsec:paradeal_finds_eta}, the objective does not always fall from one update to the next, so the last iterate of support selection is not always the best one seen. We therefore return the iterate with the smallest objective $F(\beta)$ seen during support selection, rather than the one at which support selection ends. The stopping rule is unchanged. Appendix~\ref{app:lowest_objective} shows that this choice helps most when $k$ is very small, and never makes the result meaningfully worse.

Keeping track of this iterate is done while the gradient steps run, and its cost is low. Like the rest of support selection, it is inherently parallel. Updating the running minimum is one elementwise comparison per configuration, followed by a masked overwrite of the stored iterate (\verb|torch.where|), so no configuration needs its own branch. The only extra work is computing $F(\beta)$, which costs at most one additional product with $X$ (or with $X_\mathcal{A}$ inside the active set of Section~\ref{subsec:active_set}) per iteration. For clarity, Algorithm~\ref{alg:support_selection} shows the version that returns the last iterate.

\renewcommand{\algorithmicensure}{\textbf{Return:}}

\begin{algorithm}[h]
\caption{Support Selection}
\label{alg:support_selection}
\begin{algorithmic}[1]
\REQUIRE $X, y$, initial $\beta_o$, $\lambda_1, \lambda_2, k, \eta$, max number of iterations $M_{SS}$
\ENSURE $\beta_w$ which satisfies $\|\beta\|_0 \leq k$
\STATE $\beta_w \gets \beta_o$ and $i \gets 1$
\WHILE{$i \leq M_{SS}$}
\STATE Update $\beta_w$ through gradient descent, Equation \ref{eq:1_grad_step}
\STATE \color{Cerulean} Update $\beta_w$ per soft thresholding operator, Equation \ref{eq:2_soft_threshold} \COMMENT{Only if $\lambda_1 \neq 0$} \color{black}
\STATE $\beta_w \gets H_k(\beta_w)$ per Equation \ref{explicit_beta_zeros} \COMMENT{$|S(\beta_w)| \leq k$ is guaranteed}

\IF{$S(\beta_w) = S(\beta_o)$}
\RETURN $\beta_w$
\ENDIF

\STATE $\beta_o \gets \beta_w$ and $i \gets i + 1$
\ENDWHILE
\RETURN $\beta_w$ \color{white}
\end{algorithmic}
\end{algorithm}

\subsubsection{Value Modification} \label{subsubsec:value_mod}

Once support selection creates a candidate set $\mathcal{K} \supseteq S(\beta_w)$, we solve Problem \ref{eq:vm_equation}. 

An important consideration, which we explain in Section \ref{subsec:parallelizing_vm}, emerges: Can we parallelize value modification using $\beta$s of length $k$ (as opposed to e.g. $p$), given that different $\beta$s likely have different supports? If not, we must project $\beta$s onto their supports. For the sake of completeness, we show the projection step in this section.

Once again, we consider the two cases $\lambda_1 \neq 0$ and $\lambda_1 = 0$ separately. 

\noindent {\bf Case $\lambda_1 \neq 0$}: If $\lambda_1 \neq 0$, we use the iterative algorithm given by Equation~\eqref{eq:iterative_T} with two minor changes. First, since support has been determined by algorithm \ref{alg:support_selection}, the hard thresholding step~\eqref{eq:3_hard_threshold} is replaced by a projection step $P_\mathcal{K}(\beta)$, where \begin{equation} P_\mathcal{K}(\beta_{i}) = \begin{cases} 0, & i \notin \mathcal{K} \\ \beta_{i}, & i \in \mathcal{K}. \end{cases} \label{eq:projection_split} \end{equation}

Second, instead of equal support as the termination criterion, we measure whether $||\beta_w - \beta_o||_2^2$ is sufficiently small. If $||\beta_w - \beta_o||_2^2$ is smaller than a tolerance threshold $T_{\text{VM}}$, we terminate the algorithm and return $\beta_w$ as our $\beta^*$. 

To summarize, updates for value modification are given by \begin{equation} \beta_{t+1} \in P_\mathcal{K}(\beta_t - \eta \nabla F_D(\beta_t), k, \eta \lambda_1), \label{eq:iterative_vm} \end{equation} with a convergence criterion of $\| \beta_{t+1} - \beta_t \|_2^2 < T_{\text{VM}}$.

(Compared with the pseudocode in Algorithm \ref{alg:support_selection}, value modification for $\lambda_1 \neq 0$ involves two changes: line 5 is either removed or replaced with $\beta_w \gets P_\mathcal{K}(\beta_w)$; and line 6 is replaced with the criterion $\| \beta_w - \beta_o\|_2^2 < T_{\text{VM}}$.)

\noindent {\bf Case $\lambda_1 = 0$}: If $\lambda_1=0$, to solve value modification~\eqref{eq:vm_equation}, we could use direct matrix inversion. However, when applied over large $\prod C$, finding a way to parallelize appears difficult in our experience and inversion becomes too computationally expensive. We could also use the proximal gradient-type algorithm described in Equation~\eqref{eq:iterative_vm}, but empirically we find it slower than PCG, which is what we use.

The optimality condition for Problem~\eqref{eq:vm_equation}, found by its derivative, is given by
\begin{equation} (X^TX + \lambda_2 I) \beta = X^Ty \ \text{and} \ \beta_i = 0, i \notin \mathcal{K}
\label{eq:pcg_derivative_equality} \end{equation} 
By restricting $\beta$ to the coefficients in $\mathcal{K}$, the above is equivalent to
\begin{equation} (X_\mathcal{K}^TX_\mathcal{K} + \lambda_2 I)\beta_\mathcal{K} = X_\mathcal{K}^Ty \label{eq:support_only_pcg} \end{equation} where $X_\mathcal{K} := R_S(X, \beta_w)$ and $\beta_{\mathcal{K}} := R_S(\beta_w)$. 

Problem~\eqref{eq:support_only_pcg} is of the form $A\beta = b$ with $A = X_\mathcal{K}^TX_\mathcal{K} + \lambda_2 I$, $b = X_\mathcal{K}^Ty$; and 
$A$ is symmetric positive-definite. Therefore, we can solve~\eqref{eq:support_only_pcg} with PCG (\cite{she94}). However, when we parallelize over multiple $\beta$, we face the same issue as with the $\lambda_1 \neq 0$ case noted above: different $\beta$s are likely to have different supports. As a result, from an implementation standpoint, we can adapt PCG to involve $\beta$s whose length exceeds the support size and include a projection operator where necessary to make sure that the coefficients of $\beta$ outside $\mathcal{K}$ are zero. The details are seen in Algorithm \ref{alg:pcg}. 

For the remainder of our PCG discussion, we include projections for the sake of completeness; if they are unnecessary, then projection steps can be disregarded, and $X$ can be replaced with $R_\mathcal{K}(X, \beta_w)$.

Our PCG procedure uses initializes the residual $r_o = P_\mathcal{K}(X^Ty - X^TX\beta_w - \lambda_2 \beta_w)$, then create a preconditioner $M = \text{diag}(X^TX) + \lambda_2 I$ (where 
$\text{diag}(X^TX)$ is a diagonal matrix) for the matrix $X^TX + \lambda_2I$. Next, we introduce preconditioned residual $z_o = P_\mathcal{K}(M^{-1}r_o)$, with \begin{equation} M^{-1} = \begin{bmatrix}
M_{1,1}^{-1} & 0 & \cdots & 0 \\
0 & M_{2,2}^{-1} & \cdots & 0 \\
\vdots & \vdots & \ddots & \vdots \\
0 & 0 & \cdots & M_{p,p}^{-1}
\end{bmatrix} \label{eq:diagional_matrix_inverse} \end{equation} due to the diagonal nature of $M$. Finally, we initialize search direction $\theta = z_o$. The iterative steps of PCG proceed as in Algorithm \ref{alg:pcg} and involve updating $\alpha$ (step size), $\beta_w$, $r$, $z$, and $\mu$ (conjugacy coefficient). 

\begin{algorithm}[h]
\caption{Preconditioned conjugate gradient (PCG) on candidate set $\mathcal{K}$}
\label{alg:pcg}
\begin{algorithmic}[1]
\REQUIRE{$\beta_w$, preconditioner $M$, inaccuracy tolerance $T_{\text{PCG}}$}, max number of PCG iterations $M_{PCG}$
\ENSURE{$\beta^*$ which minimizes $f_D(\beta)$}
\STATE Initialize $r_o, z_o, \theta$ as described above
\STATE $i \gets 1$
\WHILE{$i \leq M_{PCG}$}
\STATE $g \gets P_\mathcal{K}(X^TX\theta) + \lambda_2 \theta$ \COMMENT{$g$ for ease of calculation/avoid redundancy. $S(\theta)$ already equals $\mathcal{K}$}
\STATE $\alpha = \left( r_o^Tz_o \right) / \left( p^T g \right)$ \COMMENT{$\alpha$ is the step size}
\STATE $\beta_w \gets \beta_w + \alpha \theta$
\STATE $r_w = r_o - \alpha g$
\STATE $z_w = P_\mathcal{K}(M^{-1} r_w)$
\IF{$\|r_w\|_2 \leq T_{\text{PCG}} \cdot \| P_\mathcal{K}(X^Ty) \|_2$}
\RETURN $\beta_w$
\ENDIF
\STATE $\mu_o \gets r_o^Tz_o, \ \mu_w \gets r_w^Tz_w$ \COMMENT{$\mu$ is the conjugacy coefficient; search direction $p_w$ is $A$-orthogonal to previous $p$s}
\STATE $\theta \gets z_w + \frac{\mu_w}{\mu_o} \theta$

\STATE $z_o \gets z_w, \ r_o \gets r_w$ \COMMENT{$z$ and $r$ need old and new versions within an iteration of the loop}

\STATE $i \gets i+1$

\ENDWHILE

\RETURN $\beta_w$
\color{white}
\end{algorithmic}
\end{algorithm}

It is worth clarifying that when we refer to $\eta$, we might have to admit two separate learning rates, because the learning rate for support selection should not necessarily be the same as the learning rate for value modification. (It is only necessary to have two learning rates if $\lambda_1 \neq 0$.) Therefore, when $\lambda_q \neq 0$, we can think of $\eta$ as a cross product between $\eta_{SS}$ and $\eta_{VM}$. However, for simplicity, we use the term $\eta$ throughout the paper. The distinction between support selection and value modification learning rates is further explored in Section \ref{subsec:paradeal_finds_eta}.

\paragraph{The best support before value modification is not always the best support after it.} Support selection keeps the $k$ features with the largest coefficients during the hard-soft thresholding updates, but the quality of the final answer is the objective after value modification, and the two do not always agree. When two features are strongly correlated, their coefficients during support selection can be almost the same size, so which of them survives the hard thresholding step is close to a coin flip. After value modification, the feature that was left out may turn out to give the smaller objective. We checked this directly by refitting several nearby candidate supports with value modification and keeping the best (Appendix~\ref{app:multisupport}): a candidate with a smaller objective than the chosen support exists for most configurations on correlated data. Trying more supports in this way recovers only part of the gap, at a higher cost. The pair swap procedure of Section~\ref{subsec:pair_swap}, which exchanges features in and out of the support after value modification whenever doing so lowers the objective, solves this problem to a very large extent.

\subsection{The Parallelism Framework} \label{subsec:parallelism_paradigm}

We wish to parallelize Problem \eqref{eq:our_estimator} to speed up computation across different hyperparameter configurations. 

First, Problem \ref{eq:our_estimator} involves three regularization parameters $\lambda_1, \lambda_2, k$. These parameters collectively control model accuracy (prediction performance) and model compactness (number of nonzero coefficients), and practitioners need to obtain solutions to Problem \ref{eq:our_estimator} across a wide range of these parameters, because they can result in different solutions in terms of optimizing Problem~\eqref{eq:our_estimator}. This helps them select suitable tradeoffs in model accuracy and interpretability. Additionally, Algorithm~1 depends upon initialization $\beta_0$ and learning rate $\eta$; due to the nonconvexity of Problem~\ref{eq:our_estimator}, different choices of $\beta_0$ and $\eta$ can result in solutions of different quality. 

We intend to solve Problem~\ref{eq:our_estimator} on a collection of hyperparameter configurations denoted by $\mathcal{C}$.  We will use GPU parallelism to compute a batch of solutions to Problem \ref{eq:our_estimator} across the five hyperparameters $(k, \lambda_1, \lambda_2, \eta, \beta_0)$.  We allow $\mathcal{C}$ to be a grid over all hyperparameter choices (i.e. cross product) or a listed subset of pre-specified hyperparameter configurations. The former construct is called a grid; the latter, a tuple.

ParaDeAl takes advantage of GPU parallelism to get speedups from:

\begin{itemize}
    \item \textbf{Matrix-matrix or matrix-vector multiplication.}  A core building block of our ParaDeAl approach is matrix-matrix or matrix-vector multiplications, which can be performed much faster on GPUs than CPUs. This is because GPUs contain up to tens of thousands of cores, which can collectively process information much faster than CPUs, which only have a few dozen (\cite{men26}, \cite{app26}). Moreover, crucially, these multiplications are embarrassingly parallel: they are easy to subdivide or split across large numbers of cores, allowing the GPU to take full advantage and use most or all of its cores.
    
    \item \textbf{Multi-configuration computation.} The key advantage of ParaDeAl over existing algorithms is parallel processing across multiple problem setups indexed by $\mathcal C$. By calling PyTorch functions like \verb|vmap|, ParaDeAl can achieve large speedups over CPUs --- for comparisons, see the experiments in Section 4.1. In popular statistics/machine learning software, such as scikit-learn (\cite{ped18}), hyperparameter configurations are often done sequentially and/or independently of one another --- we refer the reader to Section 4.3 for empirical results and the speedups we obtain due to parallelism. 
\end{itemize}

We use two complementary ways to implement parallelism:

\begin{itemize}
    \item \textbf{\texttt{vmap} Method.}
Here, we write a function for a single problem instance, then call \texttt{vmap} to execute this function across a collection of instances. 
    For embarrassingly parallel workflows in which there is no communication across different instances and the same function is applied to all instances, we batch all instances for that workflow together using \verb|vmap|. An example is the gradient descent update in Equation~\eqref{eq:1_grad_step}: given shared (static) matrices $X$ and $y$, and different instances of $\beta_o$ (initialization), learning rate $\eta$, tuning parameter $\lambda_2$, we can write a function that computes the gradient descent update for one instance of $\beta_o, \eta, \lambda_2$, then call \verb|vmap| over all instances, specifying $X$ and $y$ as static. (We note that this operation is computationally  equivalent to creating large tensors including all possibilities for each of $\beta_o, \eta, \lambda_2$ and then performing the update in~\eqref{eq:1_grad_step} once with the large tensors instead of a single $\beta_o, \eta, \lambda_2$. However, where \verb|vmap| is available, it is preferred due to its simplicity.)

    The syntax for \verb|vmap| is \begin{equation} \verb|vmap(func_name, in_dims=(args_status))(args)| \label{eq:vmap_syntax} \end{equation} where \verb|func_name| is the name of the function to be called, \verb|args_status| represents whether each argument, in order, is to be treated as a static argument across all instances (\verb|None|) or parallelized over (\verb|0| for our purposes), and \text{args} is the list of arguments to be applied to the function. For example, if we are applying Equation \ref{eq:1_grad_step} using the function $$\verb|GD(X, y, beta, eta, lambda2)|$$ over the same $X, y, \eta$ but different $\lambda_2$ and $\beta$, then the batched procedure would look like $$\verb|vmap(GD, in_dims=(None, None, 0, None, 0))(X, y, betas, eta, lambda2s)|$$ where \verb|betas| is the tensor of all $\beta$s and \verb|lambda2s| is the tensor of all $\lambda_2$s organized in the same order as \verb|betas|.

    When tuple search is used, we apply Equation \ref{eq:vmap_syntax}. When grid search is used, we may have to extend the list of hyperparameter values. Using the same example above, supposing we have 1 value each of $k$, $\eta$, and $\lambda_1$, but 15 $\beta_o$s and 12 $\lambda_2$s, then the lists \verb|betas| and \verb|lambda2s| must be of length 180 to match the fact that we have 180 problem instances. To extend the list of $\lambda_2$s from length 12 to 180, we do not have to allocate additional memory; using operations like \verb|view|, \verb|expand|, and \verb|reshape| elongates this list without incurring extra storage cost. In this particular case, adopting the $B$ storage convention in Section \ref{subsec:beta_B_tensor}, the correct expression is $$\verb|l2s.view(-1, 1).expand(-1,|C_\eta C_\beta\verb|).reshape(-1).view(-1, 1).expand(-1,|C_kC_{\lambda_1}\verb|).T.reshape(-1)|$$ where the $C_{(\cdot)}$ are defined in Section \ref{subsec:beta_B_tensor} and \verb|l2s| is the original list of 12 $\lambda_2$s. Another example of extension is demonstrated in Equation \ref{eq:k_extension}.
    
    \item \textbf{Overarching Method.} This corresponds to the most general (overarching) case which cannot be captured in the \verb|vmap| method. We modify the single-instance function to fit all instances. Not all ParaDeAl operations can be batched, especially those which reshape matrices, batch over indexing, or require control-dependent branching (different iteration count for different instances, likely due to \verb|if| statements). As a result, it is not advisable to directly \verb|vmap| over all instances, but some level of parallelism can still be achieved if the \verb|vmap|able individual steps are still \verb|vmap|ped, while control of logic (e.g. \verb|if| statements) is not parallelized. In effect, the single-instance function changes into an overarching function, where large tensors are explicitly created to handle all instances of non-static arguments at once. An example of this is Algorithm \ref{alg:topk_parallel}.
\end{itemize}

We note that if a function is \verb|vmap|able, then we should use the \verb|vmap| Method, otherwise we should use the more general Overarching Method. Oftentimes, PyTorch code can be carefully rewritten to move from the latter category to the former. For example, if implementing the hard thresholding operation (Equation~\eqref{eq:3_hard_threshold}) over multiple $\beta$ for the same $k$ in Equation \ref{explicit_beta_zeros}, using \verb|beta[~support] = 0| to set non-support-values of $\beta$ to zero is unacceptable because we batch over a data-dependent binary overlay \verb|support| (i.e. it changes depending on $\beta$); on the contrary, writing this same operation as \verb|beta = torch.where(support != 0, beta, 0)| avoids this issue because it calls a PyTorch function \verb|where| amenable to batching.

\subsection{A multi-dimensional tensor representation of $\beta$s} \label{subsec:beta_B_tensor}

Every configuration $c = (\beta_o, \lambda_1, \lambda_2, \eta, k)$ affects the trajectory of the sequence of solutions $\beta_t$ obtained by ParaDeAl.

We create a multidimensional tensor $B$ to store $\beta$s corresponding to every $c$ as they are being updated through iterations. (We refer to $B$ as a tensor following the PyTorch conventions.)

We consider two main formats of $\mathcal{C}$. It can be (i) a grid, a cross product over each of the 5 dimensions of $\mathcal{C}_a$; or (ii) a tuple, a list of values of $c$. This creates two general search regions of interest: grid search and tuple search.

In grid search, we initialize a tensor $B$ consisting of $\prod C= C_kC_{\lambda_1}C_{\lambda_2}C_\eta C_\beta$ different coefficient vector $\beta$s, with each $\beta$ corresponding to a configuration $c \in {\mathcal C}$. We represent $B$ as a multidimensional matrix with $\prod C$ elements.

This multidimensional representation facilitates the indexing of distinct values of $\beta$ by their corresponding hyperparameters. (Note that certain dimensions of $B$ may have length one if we are considering only one instance of that hyperparameter.
For example, a user may wish to perform ParaDeAl with only L2 penalties; then, $C_{\lambda_1} = 1$.)

\begin{figure}[htbp]
    \centering
    \includegraphics[scale=0.35]{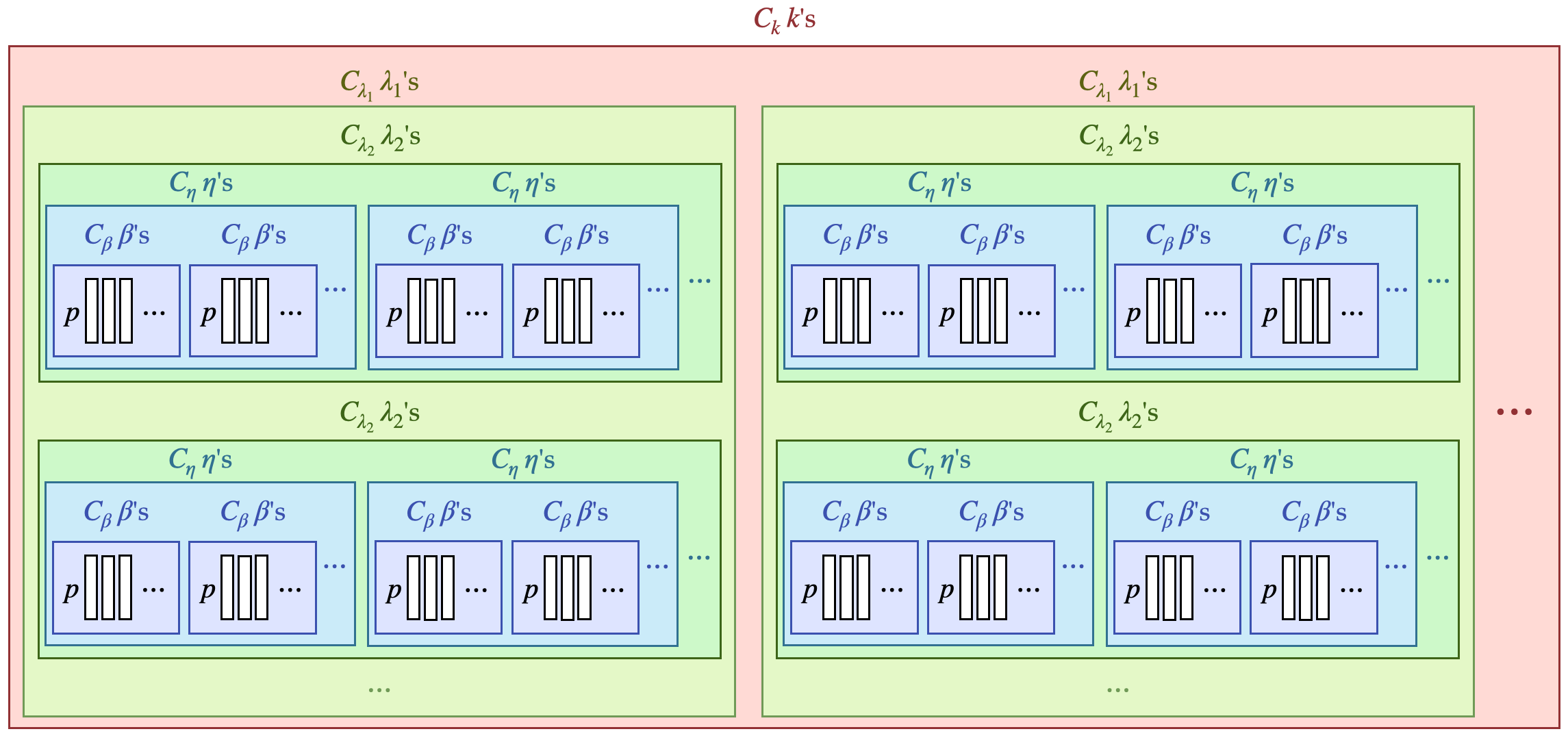}
    \caption{A concept of a 6-dimensional $B$. Each individual $\beta$ is a white length-$p$ vector. The nested structure is a byproduct of our representation of $B$ and can be described thus: Each label and its color represent the smallest unit that encompasses all values of that label from the dimension it is given. For example, since $B$ is $(C_k \times C_{\lambda_1} \times C_{\lambda_2} \times C_\eta \times C_\beta \times p)$, a single blue $C_\eta$ block is a 3-D tensor consisting of all $C_\eta$ variations of the lower dimensions, in which each purple block (labeled as ``$C_\beta$ $\beta$s") represents one $C_\eta$ over all the lower dimensions (in this case, just $C_\beta$).}
    \label{fig:7d_case}
\end{figure}

One example of a 6-dimensional tensor $B$ is given in Figure \ref{fig:7d_case}. We choose the ordering of the hyperparameter coordinates to facilitate computation: for instance, in the sequence of operations (11)--(13), the last operation is hard thresholding, so $k$ is the outermost dimension. 

The chief benefit of this representation is that all $\beta$s corresponding to any of the five hyperparameter coordinates can be indexed by specifying the indices for that variable. For example, if we wish to reference all $\beta$s using the third value of $\eta$ we specified, this can be done with \verb|B[:, :, :, 2, :, :]|. 

When hyperparameters are stored as a list in the tuple form, the tensor $B$ is simply a list of $\beta$s, with dimension $(\prod C \times p \times m)$.

\subsection{Parallelizing Support Selection} \label{subsec:parallelizing_ss}

We use \texttt{vmap} to parallelize the gradient descent and soft thresholding steps in Equations \ref{eq:1_grad_step} and \ref{eq:2_soft_threshold}. To make implementation PyTorch-compatible, if $k$ is the same across all configurations, we can use \texttt{vmap} on hard thresholding. However, if there are multiple values of $k$, we cannot use \texttt{vmap}. In Algorithm \ref{alg:topk}, we demonstrate the \texttt{vmap}-enabled case, which forms the basis for the more general approach in Algorithm \ref{alg:topk_parallel}.

\begin{algorithm}[h]
\caption{Single-$\beta$ Hard Thresholding (can be used with vmap when $k$ is same across instances)}
\label{alg:topk}
\begin{algorithmic}[1]
\REQUIRE single $\beta$, support size $k$
\ENSURE new $\beta$ with at most $k$ nonzeros

\STATE Obtain the top-$k$ indices from $|\beta|$
\STATE Generate binary overlay $O$ of 1s and 0s based on obtained top-$k$ indices \COMMENT{1 if the element is in the top-$k$}
\RETURN $O \odot \beta$ \COMMENT{$\odot$ denotes Hadamard, i.e. element-wise, product} \color{white}
\end{algorithmic}
\end{algorithm}

When different configurations have different values of the support size $k$, one approach is to run Algorithm \ref{alg:topk} in a \texttt{for} loop across different $k$ values. We can use binary indicators (1/0) which we refer to as a binary overlay to indicate which elements of $\beta$ are among the largest $k$ magnitudes. This is not GPU-friendly, so we use an alternative approach described in Algorithm \ref{alg:topk_parallel}. We can operate on the $B$ matrix of size $(\prod C \times p)$ which contains all $\beta$s as introduced in Section~\ref{subsec:beta_B_tensor}. This procedure, like Algorithm \ref{alg:topk}, uses a binary overlay, i.e. a tensor of 1s and 0s which either allows or denies each element, but this time it is a matrix of $(\prod C \times p)$ indicating whether each element of the flattened $B$ should be included or not in the final output.

\begin{algorithm}[h]
\caption{Correct Overlay Approach for Parallel Hard Thresholding}
\label{alg:topk_parallel}
\begin{algorithmic}[1]
\REQUIRE tensor of $\beta$'s $B$, support size list per $\beta$ $K$
\ENSURE new $B$ with at most $k \in K$ nonzeroes per $\beta \in B$

\STATE $B \gets \verb|torch.abs|(B)$
\STATE For each row of $B$, obtain the sort indices sort\_indices using \verb|torch.sort| \COMMENT{The smallest-magnitude element of a $\beta$ is indexed 0, and so on until $p-1$}
\STATE $R \gets [[0, 1, \dots, p-1]]$ \COMMENT{Ranged vector of size $(1 \times p)$ for overlay production}
\STATE $K$ becomes $(\prod C \times 1)$ \COMMENT{Previously it was a vector of length $\prod C$}
\STATE binary\_elts $\gets R < K$ \COMMENT{PyTorch automatically scales $R$ and $K$ to $(\prod C \times p)$; binary\_elts has the correct number of 1s and 0s but not in the right indices}
\STATE Initialize overlay $O$ of shape $(\prod C \times p)$
\STATE $O \gets $ scattered binary\_elts across $O$ using sort\_indices with row dominance \COMMENT{$O$ is the final 0-1 overlay of whether to keep the element or not}
\RETURN $O \odot B$ \color{white}

\end{algorithmic}
\end{algorithm}

Algorithm \ref{alg:topk_parallel} applies to both the cases of grid search and tuple search. In tuple search, $K$ is the list of $k$s. In grid search, as there are $C_k$ values of $k$ yet $\prod C$ configurations (thus the input $K$ needs to be of length $\prod C$), each $k$ will be repeated $R_k = (\prod C)/C_k$ times. We extend the original list of $k$s, here denoted $K_0$, by a factor of $R_k$ without using additional memory through a sequence of \verb|torch| operations, namely 
\begin{equation} K \gets K_0.\verb|view|(-1, 1).\verb|expand|(-1, R_k).\verb|reshape|(-1). \label{eq:k_extension} \end{equation} 
Note that in the above none of \verb|view|, \verb|expand|, or \verb|reshape| require additional memory allocation because they use pointer and stride manipulation. (Another possibility is to use \verb|torch.repeat_interleave| which is shorter codewise, but we do not use it because it requires more memory as each element in the new tensor gets stored in its own memory block.) Notice that in~\eqref{eq:k_extension} we selected $C_k$ as the outermost dimension of $B$: the order of $\beta$s in $B$ in Algorithm \ref{alg:topk_parallel} exactly matches the order of the corresponding $k$s in $K$. (If this were not the case, a more cumbersome extension would have to be developed.)

When the algorithm is run across multiple configurations, different instances may require different numbers of iterations to satisfy the termination criterion. Consequently, the control flow depends on the instance. Such instance-dependent control flow cannot be parallelized using standard PyTorch operations. We refer to this challenge as the \textbf{control-dependent parallelization problem (CDPP)} and propose several approaches for parallelizing support selection when the required number of iterations varies across instances.

\begin{enumerate}
    \item \textbf{Max Number of Iterations}: We run every instance of support selection for exactly $M_{SS}$ iterations; $M_{SS}$ is chosen a priori. This method allows the entire support update to be \verb|vmap|ped across all configurations, providing the most parallelism possible. However, if $M_{SS}$ is not well-chosen, the method will involve redundant computations if $M_{SS}$ is too large, 
    or cause early termination if $M_{SS}$ is too small.
    \item \textbf{Flag Approach}: We run every instance of support selection until all of them have stabilized support. At every iteration in the overarching function, we flag whether all supports have stabilized using a computationally cheap vectorized boolean operation; once all flags are true, we terminate. This method preserves \verb|vmap| structure across all iterations while adapting termination to the worst-case configuration.
    \item \textbf{Only Unstabilized}: At every iteration, we run only the instances that have yet to stabilize support. We also use the flag from the Flag Approach. After each iteration, if any instances have converged, we shrink the tensor holding the heretofore unconverged $\beta$s to remove those that just converged. While this method minimizes the number of calculations by entirely removing superfluous ones, it incurs extra cost due to shrinking $B$ a large number of times. 

    \item \textbf{Dynamic Shrinking (see Section \ref{subsec:use_dynamic_shrinking}):} This is a compromise solution. The Flag Approach requires minimal modification to $B$, but it may perform substantial redundant computation on instances that have already converged. On the other hand, Only Unstabilized eliminates this redundancy, but incurs overhead by repeatedly modifying $B$. Dynamic shrinking is a compromise between these two approaches. Given a shrinkage factor $s\in(0,1)$, we continue processing the full tensor $B$ until the number of unconverged $\beta$s falls below $s$ times the current tensor size. At that point, we compact the tensor by retaining only the unconverged instances, and proceed. We repeat this procedure until all instances have converged. In practice, we set a threshold (default 100) under which the cost of shrinking outweighs its benefits, and we proceed as in the Flag Approach.
\end{enumerate}

Since the Max Number of Iterations method requires guessing $M_{SS}$, we disregard it. Empirically, we find out that the overhead of maintaining a flag vector of size $\prod C$ is minimal, and the $B$-shrinking operation is much cheaper than an iteration of support selection or value modification, as parallel matrix-vector multiplication is more expensive. We compare the viability of the Flag Approach, Only Unstabilized, and Dynamic Shrinking in Section \ref{subsec:resolve_cdpp}.

\subsection{Parallelizing Value Modification} \label{subsec:parallelizing_vm}

The general analysis of how to parallelize value modification is the same for the $\lambda_1 \neq 0$ (proximal gradient-like procedure, Section \ref{subsubsec:value_mod}) and $\lambda_1 = 0$ (PCG, Algorithm \ref{alg:pcg}) scenarios. Because both methods rely on \verb|if| statements for termination, the CDPP applies to value modification as well.

However, there are additional challenges associated with value modification, specifically in vectorization. Once candidate sets $\mathcal{K}_1, \mathcal{K}_2, \dots$ have been chosen for the $\beta_w$s (see Section \ref{subsubsec:support_selection}), we can exploit the fact that they may be similar to condense the length of $\beta$. As we noted in Section \ref{subsubsec:value_mod}, while it is true that for a single $\beta$ one can programmatically shrink $\beta$ to just the $k$ features of its candidate set to reduce the amount of computation, for $\prod C$ different $\beta$s, the candidate sets may not consist of the same feature indices; therefore, it may be more difficult to compress the length of the $\beta$s to reduce computation. Nevertheless, we can still systematically analyze the different possible cases. We call this the \textbf{Computation-Memory Tradeoff Problem (CMTP)}, as it involves tradeoffs between amount of computation and memory footprint. (All solutions can be fully parallelized except for control-dependent branching.) We propose three ways to address the computation-memory tradeoff, all of which use the Overlay Method due to the \verb|if| criteria:

\begin{enumerate}

    \item \textbf{Resizing to Candidate Set:} We parallelize $\beta_w$s over only their candidate sets; that is, over $R_\mathcal{K}(\beta_w)$ and $R_\mathcal{K}(X, \beta_w)$ (defined in \ref{subsec:notation}), so that the support constraint is enforced by the shape of the operands rather than by an explicit projection at every iteration. Upon receiving the $\beta_w$s of length $p$ from support selection, we restrict each one to its $k$ coordinates of largest magnitude, $\beta_\mathcal{K} \gets R_\mathcal{K}(\beta_w)$, and restrict $X$ to the matching columns, producing a separate $(n \times k)$ matrix $X_\mathcal{K} \gets R_\mathcal{K}(X, \beta_w)$ for each configuration. Stacking these gives a $(\prod C \times k)$ tensor and a $(\prod C \times n \times k)$ tensor, which we supply as the dynamic arguments of the \verb|vmap|-based value modification procedure, with $y$ and the hyperparameters static.

    \item \textbf{Maintaining Length $p$:} We perform value modification over length-$p$ $\beta_w$s using projections to maintain the sparsity pattern ($\beta_i = 0, \ i \notin \mathcal{K}$; see Section \ref{subsubsec:value_mod}). To avoid creating the new $X_\mathcal{K}$ matrices, we can use the projection-based value modification procedure (see Section \ref{subsubsec:value_mod}) directly. To compute $P_\mathcal{K}(\beta)$, we use a binary overlay method: we create a length-$p$ binary overlay vector $O_S$ for each configuration, which is 1 for elements of $\beta_w$ in the support and 0 otherwise. 
    In other words, for a vector $v$, the operation $P_\mathcal{K}(v)$ is equivalent to performing $O_S \odot v$ --- this can be parallelized using the \verb|vmap| procedure.

    \item \textbf{Union Method:} We form a union of each $\beta_w$'s supports, i.e. $\mathcal{U} = \cup_{\prod C} S(\beta_w)$, resize all $\beta_w$s and $X$ to $\mathcal{U}$ once, and parallelize over length-$|\mathcal{U}|$ $\beta_w$s with projection onto $\mathcal{U}$. By removing features that are not found in any of the $\prod C$ supports, we may be able to eliminate a large portion of the $p$ features; in particular, the upper bound for the number of remaining features is  $\max(\prod C \cdot k, p)$. Then, we proceed as in the above item {\it Maintaining Length $p$}.
\end{enumerate}

Notice that in {\it Resizing to Candidate Set}, we do not restrict updates to $\beta_w$ to $S(\beta_w)$ itself, but rather $\mathcal{K}$, which is of size $k$ exactly. This is beneficial because restricting to supports, which may be of size smaller than $k$, might yield operands of different widths which can be hard to batch using \verb|vmap|. On the contrary, candidate sets are guaranteed to be of size $k$. Moreover, allowing certain zero coordinates to potentially change value may actually be beneficial, because it enlarges the feasible region of Problem \ref{eq:vm_equation} and still obeys the cardinality constraint of Problem \ref{eq:our_estimator}.

All three approaches to the computation-memory tradeoff problem must accommodate different values of $k$. For \textit{Maintaining Length $p$} and the \textit{Union Method}, this can be achieved by applying \verb|vmap| across $O_S$. For \textit{Resizing to Candidate Set}, however, direct use of \verb|vmap| is more difficult because each $\beta_w$ must be resized to its corresponding $k$, implicitly producing a ragged tensor that is difficult to parallelize. Thus, we use an alternative approach inspired by Algorithm~\ref{alg:topk_parallel}. There, multiple values of $k$ were handled using a single overlay matrix $O$, with one row per configuration. The present setting introduces an additional complication: unlike in hard thresholding, where all $\beta_w$ vectors had the same length, the vectors here have different lengths. To preserve parallelism, we pad each $\beta_w$ to length $k_M$, where $k_M$ denotes the maximum value of $k$, and apply the projections through overlay matrices using \verb|vmap|. Thus, all three approaches implement the required operations as projections.

The methods outlined above have their advantages and shortcomings. We compare these items in the table below. 
\renewcommand{\arraystretch}{1.2}
\begin{table}[htbp]
    \centering
    \begin{tabular}{|c|c|c|c|}
    \hline
        Attribute & Resizing to Cand. Set & Maintaining Length $p$ & Union Method \\ \hline
        Required preprocessing & \cellcolor{red!15}Most & \cellcolor{green!20}Least & \cellcolor{yellow!25}Medium \\ \hline
        Superfluous calculations & \cellcolor{green!20}Least & \cellcolor{red!15}Most & \cellcolor{yellow!25}Some \\ \hline
        Computational cost & \cellcolor{green!20}$\propto k_M$ & \cellcolor{red!15}$\propto p$ & \cellcolor{yellow!25}$\propto |\mathcal{U}|$ \\ \hline
        Auxiliary data structure storage & $O(\prod C \cdot nk_M)$ & $O(\prod C \cdot p)$ & $O(|\mathcal{U}|(\prod C + n))$ \\ \hline
        
    \end{tabular}
    \caption{Comparison of three methods to resolve Computation-Memory Tradeoff Problem. Green is best; red, worst.}
    \label{tab:cfrmtp_table}
\end{table}

The values of $n$, $k_M$, $p$, and $\prod C$ determine the memory requirements of the three approaches. In general, large $n$ and/or $k_M$ favor \textit{Maintaining Length $p$}, whereas small unified support $|\mathcal{U}|$ favors the \textit{Union Method} over \textit{Resizing to Candidate Set}. Runtime comparisons are less straightforward because preprocessing complexity and parallel computation complexity are inversely related across the approaches. Preprocessing is a one-time cost but incurs cost on GPUs because it requires indexing, whereas the parallelized computations are iterative but efficient on GPUs. When $k_M \ll p$, \textit{Resizing to Candidate Set} and the \textit{Union Method} are generally preferable.

\subsection{Active Set Procedure} \label{subsec:active_set}

If $p$ is large (e.g $p \geq 10^4$), support selection becomes computationally intensive. It must operate across all $p$ features across all iterations. Anticipating a sparse solution because $k \ll p$, we employ an active set procedure which is commonly used in sparse learning (\cite{fri10}, \cite{haz20}). We develop a procedure for active set support selection, or active set, which can be operated on GPUs. Our active set procedure works as follows. First, it considers a small subset $\mathcal{A}$ among all $p$ features with size $|\mathcal{A}| < p$ on which to perform support selection. After support selection on this active set $\mathcal{A}$ is done, we see if performing one more iteration of support selection on the \textit{entire} $p$-dimensional $\beta$ changes the support. If so, we update $\mathcal{A}$ to be the $|\mathcal{A}|$ elements of largest magnitude, and continue active set. If not, then support has stabilized and we exit active set. Notice that active set requires two learning rates: one corresponding to the specific $\mathcal{A}$ that accounts for its smaller size; and the ``general" learning rate for all $p$ features.

Pseudocode for active set is provided in Algorithm \ref{alg:active_set}. We keep track of the old ($A_o$) and new ($A_n$) active sets, both of length $p$, indicating which element indices are in the active set. In line 6, we obtain computational savings by performing support selection on $X_A$ of size $(n \times |\mathcal{A}_w|)$ and $\beta_A$ of length $|\mathcal{A}_w|$. In line 7 we define $E$ as an elongation operator: After line 6, $\beta_o$ is of length $|\mathcal{A}_w|$. We then restore $\beta_o$ to its original length $p$ by scattering the $|\mathcal{A}_w|$ elements back into their original indices which were precisely the 1s in $A_w$ (the other $p - |\mathcal{A}_w|$ indices correspond to zeroes). We then perform one more iteration of support selection on the full $\beta_o$ to see if the support changes; lines 9-12 cover the two cases described above.

\begin{algorithm}[h]
\caption{Active Set}
\label{alg:active_set}
\begin{algorithmic}[1]
\REQUIRE $X, y, \lambda_1, \lambda_2, k, \eta, \beta_o$, max number of iterations $M_{AS}$
\ENSURE new $\beta_w$ with $||\beta_w||_0 \leq k$

\STATE Initialize $\beta_w$ with Equation \ref{eq:2_soft_threshold} \COMMENT{equivalent to Eq. \eqref{eq:1_grad_step} if $\lambda_1 = 0$}
\STATE Initialize $A_o \gets 0_p$ and ${A}_w \gets \mathcal{O}(H_k(\beta_w))$ and $i \gets 1$ \COMMENT{$0_p$ is a vector of zeros with length $p$}
\WHILE{$A_w \neq A_o$ and $i \leq M_{AS}$}
\STATE $\beta_o \gets \beta_w$
\STATE Create $X_A \gets R(X, \mathcal{A}_w)$, $\beta_A \gets R(\beta_o, \mathcal{A}_w)$ \COMMENT{$\mathcal{A}_w = \{i | (A_w)_i = 1 \}$}
\STATE $\beta_o \gets$ Algorithm \ref{alg:support_selection} applied on $X_A$ instead of $X$ and $\beta_A$ instead of $\beta_w$ \COMMENT{support selection}
\STATE $\beta_o \gets E(\beta_o)$ \COMMENT{$E$ defined in this section as an elongation operator}
\STATE $\beta_w \gets H_k($Equation \eqref{eq:2_soft_threshold} applied on $\beta_o)$
\IF{$S(\beta_w) = S(\beta_o)$}
\RETURN $\beta_w$
\ENDIF
\STATE $A_o \gets A_w$ and $A_w \gets \mathcal{O}(\beta_w)$ and $i \gets i+1$
\ENDWHILE \color{white}

\end{algorithmic}
\end{algorithm}

The active set procedure is a loop over an iterative method, support selection which operates on a reduced set of features as updating $p$ features is computationally expensive. While we may need to perform support selection multiple times due to changes in $\mathcal{A}$, we expect this cost to be outweighed by the runtime reduction of performing support selection with smaller tensors. There is a tradeoff between the size of the active set $|\mathcal{A}|$ and the number of active set iterations: if $|\mathcal{A}|$ is large, we can expect a few iterations of high runtime; if $|\mathcal{A}|$ is small, many iterations of small runtime. Numerically, we find that $|\mathcal{A}|$ is a few multiples of $k$ and is much smaller than $p$.

Parallelizing over active set uses the same techniques as previously seen. We use the Overarching Method (see Section \ref{subsec:parallelism_paradigm}) and apply the Flag Approach or Dynamic Shrinking (see Section \ref{subsec:parallelizing_ss}) across configurations to deal with the \verb|if| condition in line 9 of Algorithm \ref{alg:active_set}. To parallelize line 6, which is support selection, we use the overarching function explained in Section \ref{subsec:parallelizing_ss} to update our $\prod C$ configurations. 

No additional machinery is required for the nested control flow. Each branch decision---including the outer \verb|if| on line 9 and the stopping criterion within support selection---is evaluated separately for each configuration and can therefore be represented by a Boolean overlay across the $\prod C$ configurations. Nested conditions are then implemented as elementwise products of the corresponding overlays.

\section{Computational Variants \& Additions to ParaDeAl} \label{sec:considerations}

Different variants of ParaDeAl emerge 
depending upon the problem dimensions and configurations.
Here we discuss two key choices: whether to use dynamic shrinkage for large $\prod C$ (Section \ref{subsec:use_dynamic_shrinking}), and whether to calculate $(X^TX)\beta$ or $X^T(X\beta)$ in gradient update (Section \ref{subsec:xtxbeta}). We also devise an algorithm for finding good learning rates $\eta$ in Section \ref{subsec:paradeal_finds_eta}. Finally, we describe pair swaps, an addition that improves the support after value modification (Section~\ref{subsec:pair_swap}).

\subsection{Dynamic Shrinking: Addressing the Control-Dependent Parallelization Problem} \label{subsec:use_dynamic_shrinking}

Support selection and value modification are iterative procedures applied simultaneously to all $\prod C$ configurations. As established in Section~\ref{subsec:parallelizing_ss}, these procedures generally require different numbers of iterations to converge across configurations. The {\it Flag Approach} mentioned earlier continues updating all trajectories until the slowest trajectory has converged. Although this preserves full vectorization, it performs redundant computations for configurations that have already converged. Alternatively, trajectories may be removed as soon as they converge, thereby eliminating redundant computations. However, the indexing operations required to identify and isolate the remaining trajectories may incur its own overhead.

Dynamic shrinking provides a compromise between these approaches, analogous to dynamic-array resizing. Rather than removing converged trajectories individually or retaining all trajectories until termination, it removes them in geometrically sized batches. A working tensor $B_R$ stores the active trajectories and their associated hyperparameters. Trajectories satisfying the convergence tolerance are flagged and frozen but remain in $B_R$ until the active count reaches a prescribed threshold. At that point, the converged $\beta$ values are written to their final positions, and $B_R$ is compacted to contain only the remaining trajectories. The shrink factor $s$ specifies the maximum fraction of active trajectories that are allowed to remain until compaction. For example, if $\prod C=1{,}000$ and $s=0.7$, the first resize occurs the moment fewer than or equal to $700$ trajectories remain active. If $686$ is the first number of active trajectories that dips below $700$, the next occurs below $480 = \lfloor 686 \cdot 0.7 \rfloor$. Thus, $B_R$ is approximately reduced to 70\% of its size at each resize. Since each resize removes a fixed fraction of the active trajectories, the number of resizes is logarithmic in $\prod C$. For $10^6$ configurations, halving and quartering require at most $14$ and $7$ resizes, respectively, before the active count falls below $N_{\min}=100$. Below this threshold, we avoid further compaction because of its additional bookkeeping overhead. 

Dynamic shrinking directly reduces computation and storage: after the $j$th resize, each gradient step operates on a tensor of size $(1-s)^j$ relative to the original. It may also reduce wall-clock time, as demonstrated in Section~\ref{subsec:resolve_cdpp}. Pruning at every iteration requires $O(M_{SS})$ gather operations, each interrupting vectorized execution, whereas dynamic shrinking requires only $O(\log \prod C)$ gathers, independently of the iteration count. Thus, the number of costly gather operations is independent of the computation length, while their benefits persist over subsequent iterations. Dynamic shrinking is therefore most effective when $\prod C$ is large and iteration counts vary substantially across configurations.

We note that after the first resize, the rows of $B_R$ no longer form a Cartesian product of the hyperparameter values; hence, the grid representation of Section~\ref{subsec:beta_B_tensor} cannot be retained. Subsequent updates instead use a tuple representation with an explicit configuration list. Algorithm~\ref{alg:dynamic_shrinking} presents pseudocode for a generic parallel update $\mathcal{P}$, applicable without modification to support selection, value modification, and the active-set procedure. Section~\ref{subsec:resolve_cdpp} evaluates the validity of dynamic shrinking.

\begin{algorithm}[h!]
\caption{Dynamic Shrinking}
\label{alg:dynamic_shrinking}

\begin{algorithmic}[1]
\REQUIRE $B$, the $\prod C$ initial $\beta$s; hyperparameter lists $H$; parallel update $\mathcal{P}$; convergence test $\mathcal{T}$; shrink factor $s$; floor $N_{\min}$; iteration cap $M$
\ENSURE $B_F$, the final $\beta$ for each of the $\prod C$ configurations
\STATE $B_R \gets B$ \COMMENT{$B_S$ is the working tensor, holding only trajectories still running}
\STATE Initialize $B_F$ of size $(\prod C \times p)$ \COMMENT{filled in as trajectories retire}
\STATE $\iota \gets (1, 2, \dots, \prod C)$ \COMMENT{$\iota_j$ is the row of $B_F$ to which row $j$ of $B_R$ belongs}
\STATE $\alpha \gets 1_{\prod C}$ \COMMENT{$\alpha_j = 1$ while trajectory $j$ is still running; $1_q$ is a vector of ones of length $q$}
\STATE $a \gets \prod C$ and $a_R \gets \prod C$ \COMMENT{number active now, and number active at the most recent resize}
\STATE $i \gets 1$ and grid2tuple $\gets$ False
\WHILE{$i \leq M$}
\STATE $B_N \gets \mathcal{P}(B_R, H)$ \COMMENT{one iteration, taken in parallel over every row of $B_R$}
\STATE $\delta \gets \mathcal{T}(B_R, B_N)$ \COMMENT{$\delta_j = 1$ if trajectory $j$ has just met its tolerance}
\STATE $B_R \gets \alpha \odot B_N + (1 - \alpha) \odot B_R$ \COMMENT{already-retired rows keep their value; \texttt{torch.where}, so no indexing}
\STATE $\alpha \gets \alpha \odot (1 - \delta)$ and $a \gets \sum_j \alpha_j$

\IF{$a = 0$}
\STATE $B_F[\iota] \gets B_R$
\RETURN $B_F$
\ENDIF

\IF{$a \leq (1-s) \cdot a_c$ AND $a \geq N_{\min}$}
\STATE $B_F[\{\iota_j : \alpha_j = 0\}] \gets \{(B_R)_j : \alpha_j = 0\}$ \COMMENT{retire every converged trajectory at once}
\STATE $B_R \gets \{(B_R)_j : \alpha_j = 1\}$, $H \gets \{H_j : \alpha_j = 1\}$, $\iota \gets \{\iota_j : \alpha_j = 1\}$ \COMMENT{the resize: one gather, incurred $O(\log \prod C)$ times in total}
\STATE $\alpha \gets 1_a$ and $a_R \gets a$
\STATE grid2tuple $\gets$ True \COMMENT{rows of $B_S$ are no longer a cross product, so $\mathcal{P}$ takes its tuple path}
\ENDIF

\STATE $i \gets i + 1$
\ENDWHILE

\STATE $B_F[\iota] \gets B_R$
\RETURN $B_F$ \color{white}
\end{algorithmic}
\end{algorithm}

\subsection{Computing $(X^TX)\beta$ Versus $X^T(X\beta)$} \label{subsec:xtxbeta}

Each gradient step in Equation~\eqref{eq:1_grad_step} requires $X^TX\beta$, where $X\in\mathbb{R}^{n\times p}$ and $\beta\in\mathbb{R}^{p}$. We know that there are two ways in which this operation can be performed, 
$(X^TX)\beta$ or 
$X^T(X\beta)$, with different computational costs. The first form permits $X^TX\in\mathbb{R}^{p\times p}$ to be precomputed and reused, whereas the right-associated form requires both multiplications at every iteration because $\beta$ changes.
The preferred form depends primarily on the dimensions of $X$. When $p\gg n$, forming and storing the Gram matrix is costly, favoring $X^T(X\beta)$. When $n\gg p$, precomputing $X^TX$ is comparatively inexpensive and its cost is amortized across iterations. In intermediate regimes, the choice also depends on the iteration count and $\prod C$. We evaluate this trade-off empirically, with further analysis provided in Appendix~\ref{app:xtx_or_xbeta}.

\subsection{Learning Rate $\eta$ Selection} \label{subsec:paradeal_finds_eta}

The analysis of Section \ref{subsec:preliminaries} provides a conservative estimate of the step size $\eta$ in Update \ref{eq:iterative_T} that is guaranteed to converge. That estimate is likely too small because $\beta$ is sparse ($k \ll p$). As a result, we would like a more aggressive choice of $\eta$, which is difficult to determine through theoretical analysis. If no $\eta$ set has been provided by the ParaDeAl user, we can introduce our own learning rate selection process which works as follows: For the operation requiring a learning rate (e.g. support selection), we test a series of learning rates for a small number $T_{LR}$ of iterations, then select the largest one that has not diverged.

To select a learning rate for support selection, we take a collection of $\eta$s geometrically spaced around $\eta_0$ and run Update \ref{eq:iterative_T} for $T_{LR} = 15$ iterations. Afterwards, we select the largest $\eta$, $\eta^*$, that has not diverged, giving us our learning rate for support selection. (We then continue running support selection with this single $\eta^*$.)

For value modification when $\lambda_1 \neq 0$, the $\eta$ selection procedure is identical: select a set of candidate learning rates; run value modification for $T_{LR} = 15$ iterations; and select the largest $\eta$ that does not diverge.

Our $\eta$ selection procedure can also be performed within active set; this will be discussed in Section \ref{subsec:active_set_results}. When active set is present, there are two separate operations that require learning rates: value modification within an active set iteration (corresponding to $\eta_{VM}$) and value modification on the iteration checking whether we change the active set (corresponding to $\eta_{SS}$). As we explained in Section \ref{subsec:two_stage_desc}, there can be several learning rates of different magnitudes depending on the particular algorithm or operation being applied.

Performing learning rate selection is computationally attractive, as it avoids performing ParaDeAl in its entirety across several learning rates. It is worth noting that $\eta^*$ may depend on two hyperparameters: $\lambda_2$, which affects the learning rate bound in Section \ref{subsec:preliminaries}; and $k$, which affects the practical amount by which the true $\eta^*$ will be larger than that bound.

\subsection{Pair Swaps: Improving the Support After Value Modification} \label{subsec:pair_swap}

Support selection followed by value modification returns a point that the hard-soft thresholding update cannot improve. A better support can still be close by, especially when features are correlated (see the note at the end of Section~\ref{subsec:two_stage_desc}). Coordinate descent methods such as L0Learn face the same issue, and \cite{haz20} address it with a local search that tries to exchange features in and out of the support. We add a step of the same kind to ParaDeAl, run after value modification and designed so that every configuration is handled at once.

\paragraph{Moves and their exact cost.} Consider one configuration, with solution $\beta$ after value modification, support $\mathcal{S}$ and residual $r = y - X\beta$. A \emph{move} removes a set $\mathcal{D} \subseteq \mathcal{S}$ of features from the support and adds a set $\mathcal{J}$ of features from outside it. We allow at most two features on each side, which is why we call the procedure \emph{pair swap}: add one feature, exchange one for one, add two, drop one and add two, and exchange two for two. Moves that add features are allowed only if the new support still has at most $k$ features.

The appeal of these moves is that their effect on the objective can be computed exactly and cheaply. Recall that every column of $X$ has unit norm, and let $c = X^Tr$ be the correlations of the features with the current residual. If the coefficients outside $\mathcal{D} \cup \mathcal{J}$ are held fixed, adding a single feature $j$ with the best possible coefficient $u$ lowers the objective by
\begin{equation}
g(c_j) = -\min_u \left( -2uc_j + (1+\lambda_2)u^2 + \lambda_1|u| \right) = \frac{\left(2|c_j| - \lambda_1\right)_+^2}{4(1+\lambda_2)}, \label{eq:swap_gain}
\end{equation}
where $(a)_+ = \max(a, 0)$, and removing a single feature $i$ raises the objective by
\begin{equation}
d_i = 2\beta_ic_i + \beta_i^2 - \lambda_1|\beta_i| - \lambda_2\beta_i^2. \label{eq:swap_cost}
\end{equation}
Once feature $i$ is removed, the residual becomes $r + x_i\beta_i$ and the correlation of feature $j$ becomes $c_j + (x_j^Tx_i)\beta_i$. The exact change in the objective of exchanging $i$ for $j$ is therefore
\begin{equation}
\Delta F(i \to j) = d_i - g\big(c_j + (x_j^Tx_i)\,\beta_i\big). \label{eq:swap_exchange}
\end{equation}
Equation~\eqref{eq:swap_exchange} shows why exchanges help most when features are correlated. If $x_j$ is nearly a copy of $x_i$, then $x_j^Tx_i \approx 1$, so feature $j$ inherits almost all of the correlation that feature $i$ carried and can take over its role, possibly at a lower objective. Hard thresholding cannot see this, because it compares coefficient sizes rather than the objective after a refit. Moves involving two added features need the best pair of coefficients instead of a single one, which is a $2 \times 2$ problem with a closed form when $\lambda_1 = 0$ and a few exact coordinate steps otherwise.

\paragraph{One round.} To keep the number of moves small, we only consider adding the $M = 16$ outside features with the largest single-feature gain $g(c_j)$. One-for-one exchanges pair each of these with every feature in the support, and the moves that drop two features use the $D = 8$ support features that are cheapest to remove. Each configuration then applies its best move, provided that the move lowers the objective by more than a relative $10^{-7}$ (smaller changes are rounding noise in \verb|fp32|). Value modification then refits all coefficients on the new support, starting from the current coefficients, and the refit is kept if it lowers the objective further. Because this refit follows every accepted move, each configuration ends pair swaps at the value modification solution on its final support, and no further value modification is needed. A configuration stops when no move lowers its objective, or after at most $M_{PS} = 200$ rounds.

\paragraph{What pair swaps guarantee.} Every accepted move strictly lowers the objective, so pair swaps can never make a solution worse, and since there are finitely many supports, the procedure stops. When it stops, no move it considers lowers the objective: the solution is a local minimum with respect to exchanging up to two features, which is a stronger property than being a fixed point of IHT. This is the same kind of guarantee as the swap-based local minima of \cite{haz20}, restricted to the candidate features above.

\begin{figure}[htbp]
\centering
\begin{tikzpicture}[
  >=Stealth,
  feat/.style={circle, draw, minimum size=5.4mm, inner sep=0pt, font=\footnotesize},
  insup/.style={feat, fill=gray!35},
  flowbox/.style={draw, rounded corners, align=center, font=\footnotesize, text width=5.0cm, minimum height=9mm, execute at begin node={\hyphenpenalty=10000}}
]
\node[font=\small\bfseries, anchor=west] at (-1.6, 1.0) {(a) One exchange};
\node[font=\footnotesize, anchor=east] at (-0.5, 0) {before};
\node[font=\footnotesize, anchor=east] at (-0.5, -1.9) {after};
\foreach \i in {1,...,8} {
  \pgfmathsetmacro{\x}{(\i-1)*0.7}
  \ifnum\i<5 \node[insup] (b\i) at (\x, 0) {$\i$}; \else \node[feat] (b\i) at (\x, 0) {$\i$}; \fi
}
\foreach \i in {1,2,3,6} { \pgfmathsetmacro{\x}{(\i-1)*0.7} \node[insup] (a\i) at (\x, -1.9) {$\i$}; }
\foreach \i in {4,5,7,8} { \pgfmathsetmacro{\x}{(\i-1)*0.7} \node[feat] (a\i) at (\x, -1.9) {$\i$}; }
\draw[->, thick] (b4.south) -- (a4.north);
\draw[->, thick] (b6.south) -- (a6.north);
\node[font=\footnotesize, align=center] at (2.45, -3.0) {drop feature 4, add feature 6\\ (shaded: in the support, $k = 4$)};
\node[font=\small\bfseries, anchor=west] at (6.6, 1.0) {(b) One round, all configurations at once};
\node[flowbox] (s1) at (9.6, 0) {correlations $c = X^T(y - X\beta)$\\ (one product for all configurations)};
\node[flowbox] (s2) at (9.6, -1.35) {score every move exactly\\ (Equations~\eqref{eq:swap_gain}--\eqref{eq:swap_exchange})};
\node[flowbox] (s3) at (9.6, -2.7) {apply each configuration's best move;\\ refit with value modification};
\node[draw, rounded corners, align=center, font=\footnotesize, text width=2.3cm, execute at begin node={\hyphenpenalty=10000}] (stop) at (14.4, -1.35) {no move lowers $F$:\\ this configuration stops};
\draw[->, thick] (s1) -- (s2);
\draw[->, thick] (s2) -- (s3);
\draw[->, thick] (s2) -- (stop);
\draw[->, thick] (s3.south) -- ++(0, -0.4) -| (6.75, 0) -- (s1.west);
\node[font=\scriptsize, anchor=north] at (9.6, -3.65) {next round};
\end{tikzpicture}
\caption{Pair swap. (a) An exchange removes one feature from the support and adds one from outside it; moves with two features on either side work the same way. (b) Each round computes the correlations of all configurations with one matrix product, scores every candidate move exactly, and applies each configuration's best move followed by a refit. A configuration leaves the batch once no move lowers its objective.}
\label{fig:pair_swap}
\end{figure}
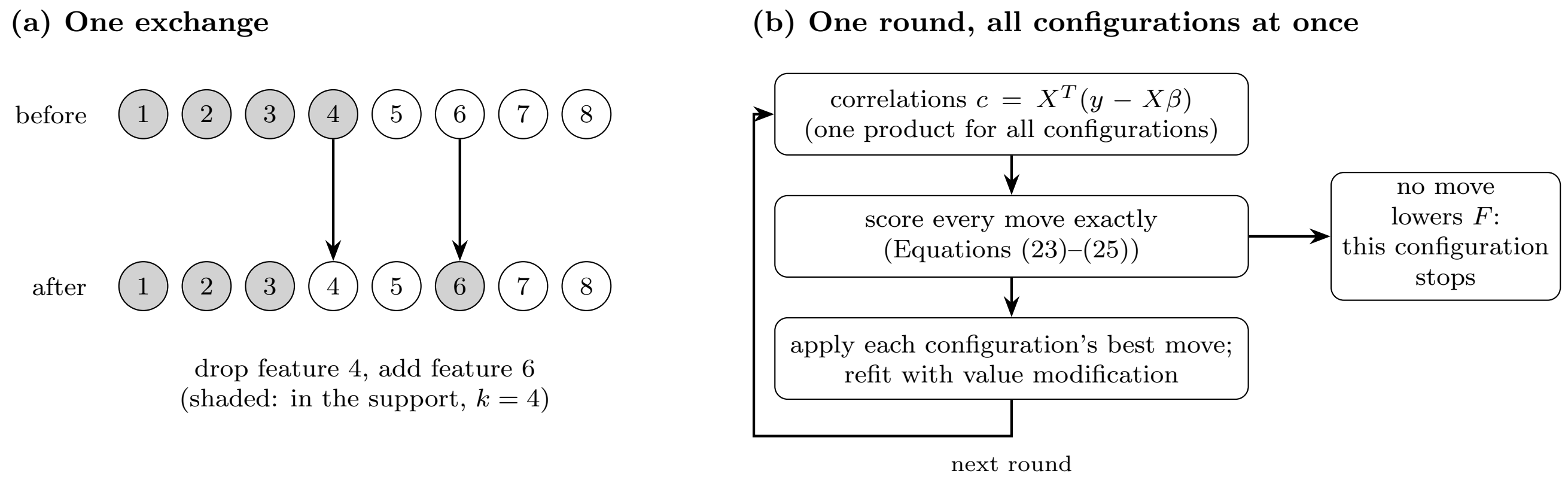

\paragraph{Pair swaps are inherently parallel.} Every step of a round is the same computation for every configuration, so a round runs for all configurations at once, as summarized in Figure~\ref{fig:pair_swap} and Algorithm~\ref{alg:pair_swap}. The correlations of all configurations come from one matrix product $X^T R$, where $R$ holds all residuals, which is the same product as in a support selection update. The scores in Equations~\eqref{eq:swap_gain}--\eqref{eq:swap_exchange} are elementwise operations on tensors of shape (configurations $\times$ candidates $\times$ support), using small blocks of $X^TX$ gathered for each configuration. The best move of each configuration is a row-wise minimum, and the refit is the parallel value modification of Section~\ref{subsec:parallelizing_vm}. Configurations with no improving move leave the batch, as in Dynamic Shrinking (Section~\ref{subsec:use_dynamic_shrinking}). A round therefore costs about as much as one full-width support selection update, plus a refit on supports of size $k$.

\begin{algorithm}[h]
\caption{Parallel Pair Swap}
\label{alg:pair_swap}
\begin{algorithmic}[1]
\REQUIRE $X, y$; solutions $B$ of all configurations after value modification; per-configuration $k, \lambda_1, \lambda_2$; candidate counts $M, D$; maximum number of rounds $M_{PS}$
\ENSURE $B$, with every row's objective no larger than before
\STATE $\alpha \gets 1_{\prod C}$ \COMMENT{$\alpha_j = 1$ while configuration $j$ is still swapping}
\FOR{$t = 1, \dots, M_{PS}$}
\STATE $R \gets y - XB$ and $G \gets X^TR$ for the rows with $\alpha_j = 1$ \COMMENT{one full-width product for all of them}
\STATE In parallel over rows: candidate additions $\mathcal{J}$ (the $M$ largest $g(c_j)$ outside the support) and drops (the $D$ smallest $d_i$)
\STATE Gather the blocks $X_\mathcal{J}^TX_\mathcal{S}$ and $X_\mathcal{J}^TX_\mathcal{J}$ for every row \COMMENT{batched products of size $M \times k$ and $M \times M$}
\STATE Score every move with Equations~\eqref{eq:swap_gain}--\eqref{eq:swap_exchange}; $\Delta F^*_j \gets$ row-wise minimum
\STATE $\alpha \gets \alpha \odot \mathbf{1}\{\Delta F^*_j < -10^{-7} F_j\}$
\IF{$\sum_j \alpha_j = 0$}
\RETURN $B$
\ENDIF
\STATE Apply the best move to every row with $\alpha_j = 1$; refit those rows with value modification and keep the refit where it lowers $F$
\STATE Remove rows with $\alpha_j = 0$ from the working batch \COMMENT{as in Dynamic Shrinking}
\ENDFOR
\RETURN $B$ \color{white}
\end{algorithmic}
\end{algorithm}

\paragraph{Effect.} Pair swaps matter most when features are strongly correlated. On synthetic data with independent features, ParaDeAl without pair swaps already reaches the objective of L0Learn and pair swaps change almost nothing. On synthetic data in which nearby features are strongly correlated, with correlation $0.9^{|i-j|}$ or $0.99^{|i-j|}$ between features $i$ and $j$, ParaDeAl without pair swaps reaches on average only 62 to 79\% of the objective quality of L0Learn (L0Learn's objective divided by ParaDeAl's), and the shortfall comes from supports of moderate size. Pair swaps raise this ratio to between 0.99 and 1.03 over three random datasets, so that ParaDeAl is often slightly better than L0Learn, and they recover about as many of the true features as L0Learn does. They add 0.1 to 0.7 seconds per grid of 900 configurations, and ParaDeAl with pair swaps still takes 0.5 to 1.6 seconds, against 5.7 to 7.4 seconds for L0Learn. Section~\ref{subsec:real_results} reports the effect on real data, and Appendix~\ref{app:pair_swap_results} gives the full synthetic results.

\section{Experimental Analysis of ParaDeAl's Components} \label{sec:method_results}

This section evaluates the components of ParaDeAl on synthetic data. We study three questions: the benefits of parallelism (on CPU and GPU) and the utility of the GPU (Section~\ref{subsec:benefits}); the relative performance of the solutions to the Control-Dependent Parallelization Problem and the Computation-Memory Tradeoff Problem (Sections~\ref{subsec:resolve_cdpp} and~\ref{subsec:resolve_cmtp}); and the gains from the active set procedure (Section~\ref{subsec:active_set_results}). Section~\ref{sec:compare_results} then compares ParaDeAl with state-of-the-art CPU methods.

For all of the following experiments, our times reported are the median of 3 trials. For the experiments in this section in particular, we have $X$ be an i.i.d. Gaussian with mean 0 and standard deviation $1$. We standardize all features of $X$ so that all features have mean $0$ and standard deviation $1$. Unless otherwise mentioned, all computation is in \verb|fp32| on a single NVIDIA L40S.

\subsection{Benefits of Parallelism (CPU, GPU) and GPU} \label{subsec:benefits}

ParaDeAl's central benefit is that it can get solutions to multiple problems in parallel, whereas competing methods can only do so iteratively. This advantage is compounded with ParaDeAl's ability to be used on both GPUs and CPUs, whereas competing methods are CPU-only.

To measure the individual and joint benefits of these two design choices, we use a $2 \times 2$ design crossing hardware \{GPU, CPU\} with execution strategy \{parallel aka ParaDeAl, iterative\}.

We consider a grid of 3,125 configurations, with each of the 5 hyperparameters receiving 5 different values. Specifically, $k \in \{5, 10, 25, 50, 100\}$, $\lambda_1, \lambda_2 \in \{0.001, 0.01, 0.1, 1, 10\}$, $\eta \in \{5, 2, 1, 0.5, 0.2\} \times 10^{-4}$, and $\beta_o$ is one of the following five choices: all 0s; all 1s; all -1s; standard Gaussian; uniform on $[-1, 1]$. There is a single $(X, y)$ with $n=500$ and $p=1000$. $y$ is generated from a 25-sparse signal at a signal-to-noise ratio $\|X\beta\|_2^2/\|\epsilon\|_2^2$ of 25 (\cite{maz22}, see Appendix \ref{app:generate_snr}). We force all iterative loops to run for exactly 1,000 iterations with no early termination (tolerance 0), so that each setup performs the same arithmetic and any difference in runtime is due to the implementation alone; that is, we use the Max Number of Iterations approach to the Control-Dependent Parallelization Problem. For the Computation-Memory Tradeoff Problem, we use {\it Maintaining Length $p$}. Each wall-clock runtime result is the median of three repetitions following one untimed warm-up, and every GPU measurement is preceded by a device synchronization, without which the recorded time reflects kernel launch rather than kernel execution. All computation is in \verb|fp32| on a single NVIDIA L40S GPU and 16 CPU threads of the same node, and the Gram matrix $X^TX$ is formed once and is shared across the four setups.

Of the 3,125 configurations, 3,124 obtain the same support under the iterative and parallel strategies, with nearly identical coefficients (absolute differences below $10^{-2}$). Upon examination, the exception happens because two different supports are selected due to floating point tolerance discrepancy on the order of $10^{-7}$. Moreover, it is GPU Iterative which differs from the other three methods (GPU ParaDeAl, CPU ParaDeAl, CPU Iterative) by 0.63\% on the objective $\|y-X\beta\|_2^2 + \lambda_1 \| \beta \|_1 + \lambda_2 \| \beta \|_2^2$.


\renewcommand{\arraystretch}{1.15}

\begin{table}[htbp]
\centering
\begin{tabular}{|c|c|c|c|}
\hline
\textbf{Hardware} & \textbf{ParaDeAl} & \textbf{Iterative} & \textbf{Speedup} \\
 & (all configurations at once) & (one configuration at a time) & (ParaDeAl / iterative) \\ \hline
GPU (L40S) & \cellcolor{green!25} 0.488s \ (0.156ms) & \cellcolor{red!25} 297.723s \ (95.271ms) & $609.6\times$ \\ \hline
CPU (16 threads) & 17.104s \ (5.473ms) & 207.035s \ (66.251ms) & $12.1\times$ \\ \hline
\textbf{Speedup} (GPU / CPU) & $35.0\times$ & $0.70\times$ & $423.9\times$ \\ \hline
\end{tabular}
\caption{Wall-clock seconds to solve all 3,125 configurations, with per-configuration milliseconds in parentheses. Numbers are the median of three trials. The marginal entries give the speedup along each factor; the lower-right entry compares ParaDeAl GPU versus iterative CPU. Green marks the fastest cell and red the slowest.}
\label{tab:four_way}
\end{table}

The row and column labeled ``Speedup" in Table~\ref{tab:four_way} marginalize/separate the two effects. Moving from CPU to GPU while using ParaDeAl yields $35.0\times$ speedup; adopting ParaDeAl while staying on the GPU yields a speedup of $609.6\times$. Thus, parallelization across configurations, rather than the hardware alone, accounts for most of the acceleration. We emphasize that ParaDeAl is not just a GPU technique: it speeds up CPU-based computation by $12.1\times$, since a single matrix product over all configurations replaces 3,125 separate products and amortizes memory traffic that the iterative loop pays repeatedly.

The slowest cell is iterative GPU, not iterative CPU; the CPU is $43\%$ faster. Iterative GPU combines the disadvantages of both settings: the loop over configurations seems to remove the simultaneity that the GPU is supposed to, while each individual iteration, a single $p \times p$ by $p \times 1$ product, is too small to occupy the device, so iterative GPU becomes bound by kernel launch overhead rather than arithmetic. This deficit depends on the problem size: at $n=1000$, $p=2000$ with 2{,}500 configurations, the same comparison inverts to $1.28\times$ in the GPU's favor, because the larger products begin to fill the device. As the dimensions (especially $p$) of $X$ grow, the GPU appears more favorable overall.

\paragraph{Does a faster GPU help?} Table~\ref{tab:gpu_types} repeats two workloads on three NVIDIA GPUs, from the L40S to the newer and faster H100 and H200, and reports the runtime and peak memory of each: the ParaDeAl GPU cell of Table~\ref{tab:four_way}, and one full grid from Section~\ref{sec:compare_results} (the Crime dataset with $30 \times 30$ configurations, solved for both penalty families exactly as in that section). The small problem of Table~\ref{tab:four_way} takes the same time on every card: it is limited by the overhead of launching kernels rather than by arithmetic, so a faster card has nothing to speed up. The larger grids run about $2\times$ faster on the H100 and H200 than on the L40S. Peak memory is the same on all three cards and far below the capacity of any of them, so the larger memory of the newer cards brings no benefit here, and the three GPUs return the same supports on all but at most 2 of the 3{,}125 configurations, with objectives that agree to about $10^{-9}$ in relative terms.

\begin{table}[htbp]
\centering
\small
\begin{tabular}{lccc}
\toprule
 & \multicolumn{3}{c}{Runtime in seconds} \\ \cmidrule(lr){2-4}
\textbf{GPU (memory)} & Table~\ref{tab:four_way} workload & Crime, L0+L2 & Crime, L1+L2 \\
\midrule
L40S (48 GB) & 0.484 & 4.04 & 8.90 \\
H100 (80 GB) & 0.488 & 2.23 & 4.49 \\
H200 (141 GB) & 0.468 & 2.12 & 4.22 \\
\midrule
Peak memory used & 0.12--0.14 GB & 9.6 GB & 2.9--3.0 GB \\
\bottomrule
\end{tabular}
\caption{The same workloads on three GPUs. Each time is the median of three trials after an untimed warm-up, in \texttt{fp32}. The L0+L2 time includes the support search with the lowest-objective rule (Section~\ref{subsec:lowest_objective}) and the pair swaps (Section~\ref{subsec:pair_swap}), which take about 0.12 seconds on every card; the L1+L2 time is ParaDeAl with FISTA, as in Section~\ref{sec:compare_results}.}
\label{tab:gpu_types}
\end{table}

\subsection{Comparing Control-Dependent Parallelization Schemes}
\label{subsec:resolve_cdpp}

Section~\ref{subsec:parallelizing_ss} introduced the Control-Dependent Parallelization Problem (CDPP): different configurations take different numbers of iterations to converge, both in support selection and in value modification. Section~\ref{subsec:use_dynamic_shrinking} described Dynamic Shrinking, which reduces redundant computation by periodically removing converged configurations from the working tensor. Here we compare the {\it Flag Approach} with {\it Dynamic Shrinking}. Note that {\it Only Unstabilized} corresponds to the limiting case of compacting after every iteration (the row $s = 1$ (guarded) below).

We consider 7 different values each for $(k, \lambda_1, \lambda_2, \eta, \beta_0)$, yielding $\prod C = 7^5 = 16,807$. We set $n=10^3$, $p=10^4$, a relative tolerance of $10^{-4}$, and an iteration cap of $10^4$ for each iterative algorithm. For \textit{Dynamic Shrinking}, we sweep across values of the shrink factor $s$, and we stop once the number of configurations dips below 100. Moreover, the grids are specified as scale-free fractions, so the same specification is meaningful at any problem size: $\eta$ as a fraction of the baseline value $1/\sigma_{\max}^2(X) = 5.796\times 10^{-5}$, $\lambda_1$ as a fraction of $\lambda_{\max} = 2\|X^Ty\|_\infty = 6132$, $\lambda_2$ as a fraction of $\sigma_{\max}^2(X) = 1.725\times 10^{4}$, and $k$ as a fraction of $p$, giving $k \in \{10, 25, 50, 100, 200, 350, 500\}$ against a true support of $k^\star = 50$. All configurations converge for these choices of $\eta$.
$\sigma_{\max}^2$ is obtained by power iteration, and the gradient is always computed as $X^T(X\beta)$ because $p > n$.

The setup is chosen so that the iteration counts vary widely across configurations. In particular, the fastest configuration reaches convergence tolerance in 13 iterations, the median in 236, the 90th percentile in 1547, and the slowest in 5,125, a factor of 394 between the extremes. All 16,807 configurations converge within the cap, ensuring that the cap does not interfere with a fair comparison. The flag approach must carry every configuration for all 5,125 iterations, performing 86.1M row-iterations where 9.5M would suffice, an excess of 811\%. Dynamic shrinking removes this excess.

Table~\ref{tab:ds_axis7} reports results across different shrink factors. Every shrink factor $s$ improves on the flag approach, and the improvement is large: the best setting, $s = 0.99$, completes in 36.122~s against 332.740~s, a speedup of 9.21$\times$, while the least effective setting tested, $s = 0.25$, still gives 4.56$\times$. The observed gather counts track the prediction $\log(\prod C / N_{\min}) / \log(1/s)$ closely, confirming that the number of resizes is logarithmic in $\prod C$ and independent of the number of iterations run. The excess column exhibits the mechanism directly: as $s$ increases, the working tensor is compacted more often and the redundancy in computation decreases. Meanwhile, the percentage of the runtime consumed by gather operations remains low.

\begin{table}[htbp]
    \centering
    \small
    \begin{tabular}{|c|c|c|c|r|r|r|r|c|c|}
    \hline
        \multicolumn{10}{|c|}{\textbf{Dynamic Shrinking}} \\ \hline
        $\mathbf{s}$ & \textbf{wall (s)} & \textbf{spread} & \textbf{speedup} & \textbf{gathers} & \textbf{pred.} & \textbf{row-iters} & \textbf{excess} & \textbf{gather (s)} & \textbf{\% run} \\ \hline
        \textbf{flag} & \cellcolor{red!15}332.740 & 0.1\% & \cellcolor{red!15}1.00$\times$ & 0 & -- & 86.1M & 811\% & 0.000 & 0.0\% \\ \hline
        $0.25$ & 73.038 & 0.0\% & 4.56$\times$ & 3 & 4 & 19.2M & 103\% & 0.006 & 0.0\% \\ \hline
        $0.4$ & 57.436 & 0.1\% & 5.79$\times$ & 5 & 6 & 15.1M & 59\% & 0.007 & 0.0\% \\ \hline
        $0.5$ & 50.905 & 0.1\% & 6.54$\times$ & 7 & 7 & 13.3M & 41\% & 0.008 & 0.0\% \\ \hline
        $0.6$ & 45.919 & 0.2\% & 7.25$\times$ & 10 & 10 & 12.1M & 28\% & 0.010 & 0.0\% \\ \hline
        $0.75$ & 41.194 & 0.1\% & 8.08$\times$ & 17 & 18 & 10.9M & 15\% & 0.015 & 0.0\% \\ \hline
        $0.9$ & 37.583 & 0.0\% & 8.85$\times$ & 47 & 49 & 10.0M & 5\% & 0.035 & 0.1\% \\ \hline
        $0.95$ & 36.674 & 0.0\% & 9.07$\times$ & 96 & 100 & 9.7M & 3\% & 0.068 & 0.2\% \\ \hline
        $0.99$ & \cellcolor{green!20}36.122 & 0.0\% & \cellcolor{green!20}9.21$\times$ & 434 & 510 & 9.5M & 0\% & 0.300 & 0.8\% \\ \hline
        $1$ (guarded) & 37.437 & 0.1\% & 8.89$\times$ & 2436 & -- & 9.5M & 0\% & 1.758 & 4.7\% \\ \hline
    \end{tabular}
    \caption{Dynamic shrinking against the flag approach at $\prod C = 16{,}807$, $n = 1{,}000$, $p = 10{,}000$ on one NVIDIA L40S. ``pred.'' is the predicted gather count $\log(\prod C/N_{\min})/\log(1/s)$; ``excess'' is row-iterations above the zero-redundancy floor; ``\% run'' is the share of GPU time spent in gathers. Median of 3 repetitions.}
    \label{tab:ds_axis7}
\end{table}

To check that the picture does not depend on the particular shape of the instance, we repeated the sweep with 5 values per axis ($\prod C = 3{,}125$) at twice the number of features, $p = 20{,}000$. The qualitative picture remains similar: every shrink factor beats the flag approach, and the best setting is again $s = 0.99$, with a speedup of 8.36$\times$.

We observe that dynamic shrinking is most attractive when $\prod C$ is large and the hardest configurations require far more iterations than the easiest.

The resize is the only operation that dynamic shrinking adds, and it is much cheaper than the size of the problem instance would suggest. At $\prod C = 16{,}807$ and $p = 10{,}000$ the working tensor occupies 0.63~GiB, and the first resize must read it, write the surviving portion, and scatter the retired rows into their final positions. This first resize takes 3.097~ms on average across shrink factors, against 64.9~ms for a single full-width iteration: a resize costs 0.048 of one iteration. Summed over an entire run, gathers account for at most 4.7\% of GPU time over all tested shrink factors. The reason is that a boolean compaction is not a serial operation. Selecting the live rows is a parallel primitive, a prefix sum over the mask followed by a single scatter, and it runs at streaming memory bandwidth. A resize therefore costs essentially one pass over the tensor, which is negligible next to the $O(np\prod C)$ arithmetic of a gradient step. It appears that reindexing does not hurt computational speed. A resize breaks the Cartesian-product representation of the configuration grid, not vectorization.

The repeated sweep shows what that ratio depends on. With $p = 20{,}000$ and $\prod C = 3{,}125$, a resize costs 0.055 of one full-width iteration, against 0.048 above, even though $p$ doubles and $\prod C$ shrinks by a factor of 5.4. This agrees with a simple cost count. A gather moves $O(p \prod C)$ bytes once, whereas an iteration performs $O(n p \prod C)$ floating-point operations, so their ratio scales as $1/n$ and is insensitive to both $p$ and $\prod C$. 

In practice, the relative cost of a resize is controlled by the number of observations alone: it does not grow with the size of the configuration grid or the number of features, which are exactly the quantities that are large when dynamic shrinking is worth using.

The shrink factor trades the number of gathers against the redundancy left in the tensor between them. Small $s$ waits until most configurations have retired and pays for few, large resizes; large $s$ compacts as soon as a small fraction has finished and pays for many small ones. Because a gather is so cheap, the optimum sits far toward the aggressive end: $s = 0.99$ here, with 434 gathers costing 0.300~s, or 0.8\% of that arm's runtime. Beyond the optimum, the curve turns over: at $s = 1$ the tensor is compacted 2436 times and comes within 0\% of the floor, but the gathers now consume 4.7\% of the run and the arm is slower, 37.437~s against 36.122~s. This is the only regime in which the cost of reindexing is visible at all. The practical guidance is that $s$ should be chosen large. 
The penalty for over-shrinking is mild and bounded by the gather bill, whereas the penalty for under-shrinking is the full redundancy of the flag approach.

We note two limitations. First, kernel launch costs form an overhead. Once the working $B$ becomes small, an iteration is dominated by launch overhead rather than arithmetic, which is why $N_{\min}$ exists. For the values of $p$ chosen, the overhead is not large enough relative to arithmetic. However, on smaller instances, dynamic shrinking will look less attractive. Second, dynamic shrinking affects the order in which floating-point operations are reassociated, resulting in bitwise different results to those done by the flag approach, and consequently different supports and objective values. This issue is further explored in Appendix \ref{app:bitwise_differences}.

\subsection{Analyzing Computation-Memory Tradeoffs} \label{subsec:resolve_cmtp}

Section~\ref{subsec:parallelizing_vm} introduced the Computation-Memory Tradeoff Problem: after support selection, different configurations have different supports, and there are several ways to operate on them in parallel. The best method depends on the problem dimensions. We compare the three methods in two setups. In both setups we use five different starting points $\beta_o$, so that the candidate sets coming out of support selection actually differ across configurations, which is what the Computation-Memory Tradeoff Problem must centrally address:
\begin{itemize}
    \item $n = 1,000, \ p = 10,000, \ k = 100, \ C_{\lambda_1} = C_{\lambda_2} = 100, \ C_{\beta} = 5, \ \prod C = 50,000$. We have a sizable number of features with a large number of different L1 and L2 penalties, leading to tens of thousands of configurations. This setup tests whether memory becomes prohibitive for Resizing to Candidate Set, which requires a total of $50{,}000 \cdot 100 = 5{,}000{,}000$ gathered columns.
    \item $n = 1,000, \ p = 10,000, \ k = \{ 10, 20, 30, 50, 100 \}, \ C_{\lambda_1} = C_{\lambda_2} = 5, \ C_{\beta} = 5, \ \prod C = 625$. With a much smaller $\prod C$, this setup measures the benefit of working with shorter vectors, as in Resizing to Candidate Set and the Union Method, rather than computing over all $p$ features, as in Maintaining Length $p$.
\end{itemize}

In both setups, we generate $y$ through a 50-sparse signal with an SNR of 25 (see Appendix \ref{app:generate_snr}). We choose our penalty parameters ($\lambda_1, \lambda_2$) identical to Section \ref{subsec:resolve_cdpp}. For all configurations, $\lambda_1 \neq 0$, so value modification uses the proximal gradient update in Equation~\eqref{eq:iterative_vm}. The 5 starting points are the same as in Section \ref{subsec:benefits}. Each configuration uses learning rate $\eta = 0.9 \times 2/L_D$.

We run support selection (Algorithm \ref{alg:support_selection}) in parallel over all configurations, with a cap of $M_{SS} = 1,000$ iterations. All configurations stabilized well under $M_{SS}$ iterations. The candidate set of each configuration is computed once from this output and given to all three methods, so they all solve exactly the same restricted problem. To equalize the cost of all three methods, we set the value modification tolerance to zero and run exactly 500 iterations. Each single-configuration update is written once and batched with \verb|vmap| and gradient is computed as $X^T(X\beta)$. The reported time includes the one-time preprocessing of each method (building the overlays, and for \textit{Resizing to Candidate Set} and \textit{Union Method}, gathering the columns of $X$) and the final step that writes the $\beta$s back to length $p$. Peak memory is measured above what is already stored before value modification begins, namely $X$, $y$, the output of support selection, and the hyperparameter lists. The results are in Table \ref{tab:cmtp}.

\begin{table}[htbp]
    \centering
    \begin{tabular}{|c|c|c|c|c|c|c|}
    \hline
        \textbf{Setup} & \textbf{Method} & \textbf{Width} & \textbf{Preproc. (ms)} & \textbf{Iters (s)} & \textbf{Total (s)} & \textbf{PM (GiB)} \\ \hline
        \multirow{3}{*}{\shortstack{$\prod C = 50{,}000$ \\ $k = 100$}} & Resizing to Cand. Set & 100 & 474.9 & 29.752 & 30.230 & \cellcolor{red!15}20.527 \\ \cline{2-7}
         & Maintaining Length $p$ & 10,000 & 17.1 & 69.208 & \cellcolor{red!15}69.225 & 13.234 \\ \cline{2-7}
         & Union Method & 499 & 7.4 & 3.612 & \cellcolor{green!20}3.629 & \cellcolor{green!20}2.057 \\ \hline
        \multirow{3}{*}{\shortstack{$\prod C = 625$ \\ $k \in \{10, \dots, 100\}$}} & Resizing to Cand. Set & 100 & 1.8 & 0.324 & 0.326 & 0.257 \\ \cline{2-7}
         & Maintaining Length $p$ & 10,000 & 0.6 & 0.567 & 0.568 & 0.166 \\ \cline{2-7}
         & Union Method & 429 & 0.3 & 0.063 & 0.064 & 0.028 \\ \hline
    \end{tabular}
    \caption{The three Computation-Memory Tradeoff Problem methods on 500 iterations of value modification, $n = 1{,}000$, $p = 10{,}000$, five starting points $\beta_o$, one NVIDIA L40S. ``Width'' is the length of each $\beta$ during the iterations ($k_M$, $p$, and $|\mathcal{U}|$, respectively). ``Preproc.'' is the one-time setup before the iterations (building overlays and gathering columns of $X$), ``Iters'' is the time for the 500 value modification iterations, and ``Total'' is their sum plus writing the $\beta$s back to length $p$. ``PM'' is the peak memory above what is stored before value modification begins. Medians of three trials.}
    \label{tab:cmtp}
\end{table}

Empirically, we find that the final objective values $\|y-X\beta\|_2^2 + \lambda_1 \| \beta\| _1 + \lambda_2 \| \beta \|_2^2$ and $\beta$s are near-identical across the three methods, with a maximum relative difference in $\|y-X\beta\|_2^2$ of $4.7 \times 10^{-7}$ and a maximum absolute difference in any coefficient of $1.1 \times 10^{-6}$ (the largest coefficient is about 3 in absolute value). 

The first setup shows the memory cost of {\it Resizing to Candidate Set}. Gathering $k = 100$ columns of $X$ for each of the 50,000 configurations creates $5 \times 10^6$ columns, a $(50{,}000 \times 1{,}000 \times 100)$ tensor of 18.6 GiB, which is about 500 times the size of $X$ itself. Together with the length-$p$ output, this brings the peak to 20.5 GiB above the baseline. This still fits on the 48 GB L40S, but since the tensor grows as $O(\prod C \cdot n k_M)$, it would run out of memory at roughly $\prod C = 120{,}000$ for the same $n$ and $k$, less than three times the size of this setup. {\it Maintaining Length $p$} needs 13.2 GiB, which comes from holding several $(\prod C \times p)$ tensors ($\beta$s, the overlay $O_S$, and the intermediate results of each step). The \textit{Union Method} needs 2.1 GiB, almost all of which is the length-$p$ output. 

The explanation for the runtimes is less obvious. Table \ref{tab:cfrmtp_table} ranks \textit{Resizing to Candidate Set} as doing the fewest superfluous calculations, and it does do the least arithmetic of the three. Yet it is 8.3 times slower than the \textit{Union Method} in the first setup, and only 2.3 times faster than \textit{Maintaining Length $p$}, even though it works with vectors 100 times shorter. We find that almost all (98\%) of the runtime is spent in the two batched products $X_\mathcal{K}\beta_\mathcal{K}$ and $X_\mathcal{K}^Tr$, which take 30.0 ms and 28.7 ms. Each of them reads the whole 18.6 GiB tensor of gathered columns at 670 to 700 GB/s, which is about 80\% of the L40S's peak bandwidth of 864 GB/s. This means the GPU is mostly waiting on memory, not computing. The reason is that each configuration has its own copy of its columns, and a column that has been read is used for only one configuration. If all configurations share one matrix, the same arithmetic turns into a single matrix-matrix product that reads each column of $X$ once and uses it for every configuration. On the same GPU this takes 1.4 ms instead of 58.6 ms. \textit{Maintaining Length $p$} and the \textit{Union Method} both work this way, but the \textit{Union Method} has a much shorter width, and it ends up 19 times faster than \textit{Maintaining Length $p$}. \\

The width of the union depends on how much the candidate sets overlap. In the first setup, the 50,000 configurations have 375 distinct candidate sets, and their union $\mathcal{U}$ has 499 features, or 5\% of $p$. This figure is $\beta_o$-dependent. The 10,000 configurations that start from $\beta_o = 0$ give only 3 distinct candidate sets with a union of 101 features, since the first update ranks the features by $|X^Ty|$ for all of them. Each of the other four starting points gives between 67 and 131 candidate sets, with a union of 211 to 254 features. We use several starting points in this experiment for exactly this reason. With $\beta_o = 0$ alone the candidate sets would be nearly the same, and the Computation-Memory Tradeoff Problem would hardly come up. If the candidate sets have less in common, for example with many more random starting points or a weak signal, $|\mathcal{U}|$ gets close to $\min(\prod C \cdot k, p)$ and the \textit{Union Method} is no better than \textit{Maintaining Length $p$}. \textit{Resizing to Candidate Set} does not have this problem, but it needs much more memory.
  
The second setup is small. All three methods finish in under 0.6 s and use at most 0.26 GiB, but the order is the same as before. The \textit{Union Method} is 5.1 times faster than \textit{Resizing to Candidate Set} and 8.9 times faster than \textit{Maintaining Length $p$}. With 625 configurations and $|\mathcal{U}| = 429$, one iteration of the \textit{Union Method} takes 0.13 ms, which is about as fast as a GPU iteration gets because of kernel launch overhead. The other two methods are big enough that their actual work sets the time. Padding to $k_M$ with the overlay, which we use to handle several values of $k$, gives the same $\beta$s as the other two methods. Since we are talking about fractions of a second here, the choice of method only really matters when $\prod C$ is large.
  
Based on these results, we use the \textit{Union Method} by default. It is the fastest in both setups, it used the least memory, and its preprocessing took under 8 ms. When $|\mathcal{U}|$ is close to $p$, we switch to \textit{Maintaining Length $p$}, which uses less memory than \textit{Resizing to Candidate Set} in both setups. We only use \textit{Resizing to Candidate Set} when $|\mathcal{U}|$ is large, $k_M \ll p$, and the $O(\prod C \cdot n k_M)$ tensor fits easily in GPU memory.

\subsection{Benefits of the Active Set Procedure} \label{subsec:active_set_results}

We call an update that operates on all $p$ features a \emph{full-width} update. Without an active set, every support selection update is a full-width update, which costs $O(np)$. The active set procedure of Section~\ref{subsec:active_set} runs support selection on a subset $\mathcal{A}$ of the features instead, followed by a single full-width update that checks whether any feature outside $\mathcal{A}$ should enter the support. This section shows that, for large $p$, this gives substantial speedups without making the solution worse.

\paragraph{Setup.} For each test case we run ParaDeAl twice independently, once with the active set and once without. In each run, ParaDeAl parallelizes across the five starting points $\beta_o$ of Section~\ref{subsec:benefits} (all zeros, all ones, all negative ones, standard Gaussian, and uniform on $[-1, 1]$), with one value each of $k$, $\lambda_1$ and $\lambda_2$. We take $n = 1{,}000$ and $p \in \{10^5, \ 3 \times 10^5, \ 10^6\}$, and generate the response from a 50-sparse coefficient vector with an SNR of 25. We consider $k \in \{10, 50, 100\}$ at every $p$, and $k \in \{250, 500\}$ at $p \geq 3 \times 10^5$, one value of $k$ at a time. The penalties are $\lambda_1 = 10^{-3}\lambda_{\max}$ and $\lambda_2 = 10^{-3}\sigma_{\max}^2(X)$. The gradient is computed as $X^T(X\beta)$, since forming the whole gram matrix $X^TX$ would take 4 TB in \verb|fp32| at $p = 10^6$.

The arm without an active set needs one learning rate, for the full-width updates. The arm with the active set needs a second one, for support selection on $X_\mathcal{A}$. Both are chosen by the learning rate selection of Section~\ref{subsec:paradeal_finds_eta}, with candidates $\eta_{\text{pred}} \times \{1/16, 1/8, 1/4, 1/2, 1, 2\}$, each run for 10 iterations. Here $\eta_{\text{pred}} = 0.9/(\sigma_k^2 + \lambda_2)$, where $\sigma_k^2 = n(1 + \sqrt{k/n})^2$ approximates $\sigma_{\max}^2$ of a $k$-column submatrix of $X$. Both loops use the stopping rule of Section~\ref{subsubsec:support_selection}: support selection on $X_\mathcal{A}$ stops when one update leaves the support unchanged, and the outer loop stops when the full-width check leaves the support unchanged. To break cycles of the form $A, B, A$, each loop also stops when the support (in the outer loop, the active set) returns to the one from two steps back. The iteration number limits are $M_{SS} = 1{,}000$ and $M_{AS} = 100$, which are never exceeded in any test.

To compare how good the supports provided by non-active set and active set are, we run value modification on the support each arm returns and compute the objective $F(\beta)$ of Problem~\eqref{eq:our_estimator}, including both penalty terms. Value modification is not timed.

\begin{table}[htbp]
    \centering
    \small
    \setlength{\tabcolsep}{4.5pt}
    \begin{tabular}{|c|c|c|c|c|c|c|c|c|c|c|c|}
    \hline
        \multicolumn{12}{|c|}{\textbf{Active Set}} \\ \hline
        \multirow{2}{*}{$p$} & \multirow{2}{*}{$k$} & \multirow{2}{*}{$|\mathcal{A}|$} & \multicolumn{3}{c|}{\textbf{No active set}} & \multicolumn{3}{c|}{\textbf{Active set}} & \multicolumn{2}{c|}{\textbf{Speedup}} & \textbf{obj.} \\ \cline{4-11}
         & & & \textbf{total} & \textbf{solve} & \textbf{upd.} & \textbf{total} & \textbf{solve} & \textbf{upd.} & \textbf{total} & \textbf{solve} & \textbf{diff. ($\downarrow$)} \\ \hline
        \multirow{3}{*}{$10^5$} & 10 & $4k$ & 23.5 & 9.7 & 6 & 21.1 & 6.0 & 3 & 1.11$\times$ & 1.62$\times$ & $+0.2\%$ \\ \cline{2-12}
         & 50 & $5k$ & 42.6 & 28.7 & 21 & 28.8 & 13.7 & 5 & 1.48$\times$ & 2.09$\times$ & $-1.0\%$ \\ \cline{2-12}
         & 100 & $2k$ & 18.2 & 4.3 & 3 & 26.2 & 11.2 & 4 & \cellcolor{red!15}0.69$\times$ & \cellcolor{red!15}0.39$\times$ & $-69.5\%$ \\ \hline
        \multirow{4}{*}{$3 \times 10^5$} & 10 & $4k$ & 70.8 & 28.5 & 7 & 57.7 & 14.1 & 3 & 1.23$\times$ & 2.02$\times$ & $+0.1\%$ \\ \cline{2-12}
         & 50 & $5k$ & 123.5 & 81.3 & 20 & 70.5 & 26.9 & 5 & 1.75$\times$ & 3.03$\times$ & $-0.3\%$ \\ \cline{2-12}
         & 100 & $2k$ & 54.5 & 12.3 & 3 & 65.3 & 21.8 & 4 & \cellcolor{red!15}0.84$\times$ & \cellcolor{red!15}0.56$\times$ & $-58.6\%$ \\ \cline{2-12}
         & 500 & $5k$ & 160.8 & 118.4 & 29 & 68.0 & 24.4 & 3 & 2.37$\times$ & 4.85$\times$ & $-30.7\%$ \\ \hline
        \multirow{5}{*}{$10^6$} & 10 & $3k$ & 226.5 & 87.6 & 6 & 177.5 & 40.6 & 3 & 1.28$\times$ & 2.16$\times$ & $-0.1\%$ \\ \cline{2-12}
         & 50 & $5k$ & 317.0 & 181.4 & 14 & 206.0 & 69.1 & 5 & 1.54$\times$ & 2.62$\times$ & $-0.4\%$ \\ \cline{2-12}
         & 100 & $2k$ & 433.0 & 297.3 & 23 & 179.5 & 42.6 & 3 & 2.41$\times$ & 6.98$\times$ & $-34.9\%$ \\ \cline{2-12}
         & 250 & $4k$ & 419.5 & 284.3 & 21 & 243.4 & 106.7 & 7 & 1.72$\times$ & 2.66$\times$ & $+1.7\%$ \\ \cline{2-12}
         & 500 & $2k$ & 1170.4 & 1034.9 & 80 & 202.5 & 65.8 & 4 & \cellcolor{green!20}5.78$\times$ & \cellcolor{green!20}15.72$\times$ & $-1.7\%$ \\ \hline
    \end{tabular}
    \caption{Active set against full-width support selection, five starting points in parallel, $n = 1{,}000$. Times are in milliseconds. ``total'' includes learning rate selection, ``solve'' excludes it, and ``upd.'' counts full-width updates. ``obj.\ diff.'' is the relative change in the objective $F(\beta)$ after value modification, averaged over the five starting points; negative values mean the active set found the better support. For each $(p, k)$, $|\mathcal{A}|$ is the size from the sweep $\{2k, 3k, 4k, 5k, 8k, 10k, 15k\}$ of Appendix~\ref{app:as_size} that gives the best speed without a worse objective. Red cells mark the two rows where the arm without an active set stops early on a much worse support.}
    \label{tab:as_main}
\end{table}

Table~\ref{tab:as_main} shows the results. In ten of the twelve rows the active set is faster, by 1.1 to 5.8 times end to end and by 1.6 to 15.7 times on the solve alone (i.e. excluding time for the learning rate search), and its objective is within 2\% of the arm without an active set or better. The speedup grows with $p$ and $k$. The reason is the number of full-width updates. Without an active set, every update is full-width, and more of them are needed as $k$ grows: 80 at $p = 10^6$, $k = 500$. With an active set, there is one full-width update per round, and the active set needs only 3 to 7 of them in these rows. The rest of the work runs on at most $|\mathcal{A}| = 2{,}500$ features instead of $10^6$. The gap between ``total'' and ``solve'' is learning rate selection, a fixed cost that both arms pay once per batch. As $\prod C$ grows this cost is spread over more configurations, and the end-to-end speedup approaches the solve-only speedup (Appendix~\ref{app:as_amortize}).

In the two red rows, at $k = 100$ with $p = 10^5$ and $3 \times 10^5$, the arm without an active set stops after 3 updates. It is faster only because it stops early on a much worse support: the active set lowers the objective by 59 to 70\%. In these rows the active set trades a slower run for a much better answer. Appendix~\ref{app:as_supports} shows why the two arms can end on different supports. In short, from a dense starting point the full-width step is too small to move away from the support given by $\beta_o$. The inner loop of the active set works on $X_\mathcal{A}$, where it can take larger steps.

\begin{figure}[htbp]
    \centering
    \includegraphics[width=\textwidth]{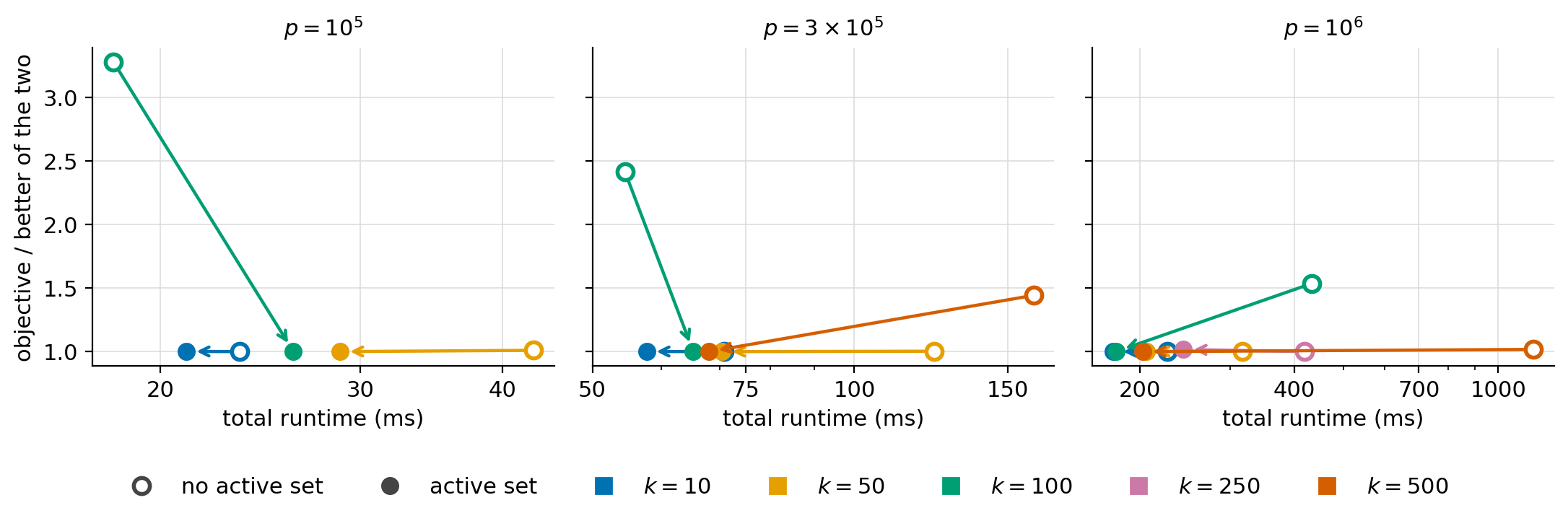}
    \caption{Runtime against objective for the rows of Table~\ref{tab:as_main}. The objective is divided by the better of the two runs, so 1 marks the better run. Open circles: no active set; filled circles: active set; each arrow points from the first to the second. Colors mark $k$; the horizontal axis is logarithmic. Accompanying numbers in Table~\ref{tab:as_main} present some additional information.}
    \label{fig:as_runtime_objective}
\end{figure}

Figure~\ref{fig:as_runtime_objective} shows the same rows as runtime against objective. Each arrow goes from the run without the active set to the run with it. Most arrows point to the left along the line at 1, so the active set gives the same objective in less time. The two arrows that point to the right and down are the red rows of Table~\ref{tab:as_main}, where the active set is slower but finds a much better support. In two further rows ($p = 3\times10^5$, $k = 500$ and $p = 10^6$, $k = 100$), the active set is both faster and better.

The choice of $|\mathcal{A}|$ matters more for accuracy than for speed. Every row of Table~\ref{tab:as_main} uses an $|\mathcal{A}|$ between $2k$ and $5k$, which is the range we recommend. Larger active sets tend to stop the outer loop too early. At $p = 10^6$ and $k \geq 100$, $|\mathcal{A}| = 10k$ and $15k$ give objectives 8 to 157\% worse, and at $k = 100$ and $500$ they stop after one or two rounds. Appendix~\ref{app:as_size} gives the full sweep over $|\mathcal{A}|$, the variation between runs, and the limitations of the procedure.

Peak memory is similar for the two arms, since it is dominated by the length-$p$ tensors of the five configurations: at $p = 10^6$, 149 MiB without the active set and 175 to 195 MiB with it. Throughout this section each batch has a single value of $k$. Batching several values of $k$ together is possible if $|\mathcal{A}|$ is shared across the batch (Appendix~\ref{app:as_across_k}).

\section{ParaDeAl Versus State-of-the-Art} \label{sec:compare_results}

This section compares ParaDeAl with established CPU-based algorithms on three real datasets, which we expand to between $1.9\times10^4$ and $2.8\times10^5$ features. We then repeat the comparison on three further datasets whose features are less strongly correlated, to see how much the results depend on correlation. In every comparison the competitor chooses the hyperparameter configurations and ParaDeAl solves the same configurations, so speed and accuracy are compared configuration by configuration.

\subsection{The Competitors}
\label{sec:competitors}

Table~\ref{tab:competitors} lists the methods. L0Learn \citep{haz20,haz23} is a state of the art algorithm for $L_0$-regularized regression and is the benchmark for the L0+L1 and L0+L2 modes. For the elastic net (L1+L2), glmnet \citep{fri10} a popular solver, scikit-learn \citep{ped18} is also a widely used Python implementation, and skglm \citep{ber22} is a recent Python library designed for speed. L0Learn and glmnet compute a warm-started regularization path, in which each solution is the starting point for the next. scikit-learn and skglm fit one configuration at a time, as a user sweeping a grid would. None of the four parallelizes across hyperparameter configurations, which is the axis ParaDeAl is built around, and none solves the L0+L1+L2 mode.

Two differences between ParaDeAl and the competitors shape the comparison. First, L0Learn penalizes the support size instead of constraining it, so we let L0Learn compute its path and give ParaDeAl the support size and penalty of every solution that L0Learn returns. We run L0Learn without its swap-based local search method which is more compute-intensive. Second, the four competitors run on a single CPU core, since their core algorithms are sequential, while ParaDeAl runs on one NVIDIA L40S GPU. The speedups below therefore compare one GPU with one CPU core. We include no mixed integer solver: at these values of $p$, an exact solver such as Gurobi does not finish a single configuration in the time the other methods take for a whole grid. Appendix~\ref{app:competitors} describes the algorithms, the software versions and settings, and how each method's penalty parameters convert to those of ParaDeAl.

\paragraph{Warm starts.} The methods differ in what they reuse from one configuration to the next. L0Learn and glmnet warm-start: along each path, every fit starts from the previous solution, which is a large part of what makes a path cheap. scikit-learn and skglm start every configuration from zero, as a user fitting one configuration at a time would. ParaDeAl starts every configuration from $\beta = 0$ as well. The only things it shares across configurations are the step-size search, which is run once per distinct $(k, \lambda_2)$, and, within pair swaps, each refit starts from the current coefficients of the same configuration. ParaDeAl therefore gets no help from neighboring solutions, so on runtime the comparison with the two path methods is, if anything, tilted in their favor.

\begin{table}[htbp]
\centering\footnotesize\setlength{\tabcolsep}{4pt}
\begin{tabular}{l|c|c|c|l}
\toprule
Method & Problem & Algorithm & Across configurations & Language \\
\midrule
L0Learn & $L_0$ penalty (+ $L_1$ or $L_2$) & coordinate descent & warm-started path & C++ (R) \\
glmnet & elastic net & coordinate descent & warm-started path & Fortran (R) \\
scikit-learn & elastic net & coordinate descent & one at a time & Cython (Python) \\
skglm & elastic net & working sets + accelerated CD & one at a time & Numba (Python) \\
ParaDeAl & $L_0$ constraint and/or $L_1$, $L_2$ & proximal gradient, IHT, CG & all at once (GPU) & PyTorch \\
\bottomrule
\end{tabular}
\caption{The methods compared.}
\label{tab:competitors}
\end{table}

\subsection{Datasets}
\label{subsec:datasets}

We use one dataset from each of three application areas (Table~\ref{tab:datasets}). In healthcare, the Diabetes 130-US Hospitals data \citep{uci_diabetes,str14} are used to predict the time a patient spends in hospital. In socioeconomics, the Communities and Crime data \citep{uci_crime,red02} are used to predict the rate of violent crime in a community. In meteorology, the Air Quality data \citep{uci_airquality,dev08} are used to predict the hourly temperature.

For each of these datasets, we expand the number of features by using pairwise interactions (the original columns together with every product $x_ix_j$, $i<j$) and nearest-neighbor features (each row is augmented with the features of its $h$ nearest rows). Every column of the training design is then scaled to mean 0 and unit $L_2$ norm, and the test design receives the same transformation. Appendix~\ref{app:data_construction} gives the full construction. The resulting datasets are in Table~\ref{tab:datasets}.

\paragraph{Datasets with softer correlations.} The neighbor features above are close copies of one another, so these three datasets contain many pairs of nearly identical columns: in Air Quality every feature has another feature with correlation above 0.9, and in Diabetes 41\% of the features have a near-exact duplicate. Strong correlations are exactly the setting in which the choice of support is hardest (Section~\ref{subsec:pair_swap}). To see how the comparison changes when correlations are weaker, we build three further datasets from polynomial features alone, without neighbors. Each base feature is first replaced by its normal score (its rank, mapped to the corresponding quantile of a standard normal distribution), and the features are then all products of one, two or three distinct base features. We apply this to Communities and Crime again (giving the overlap with the first group), to a healthcare dataset from the US National Health and Nutrition Examination Survey (NHANES) \citep{nhanes}, with systolic blood pressure as the response, and to a European weather dataset \citep{huber_weather}, with the next day's temperature in Basel as the response. In NHANES and Weather no two columns have correlation above 0.99, and only 15 to 19\% of the features have a partner with correlation above 0.9. Crime keeps strong correlations even without neighbors: 86\% of its features have a partner above 0.9 and 20\% a partner above 0.99, so it is closer to the first group than to the other two. The datasets are in the lower half of Table~\ref{tab:datasets}, and Appendix~\ref{app:data_construction} gives the construction.

\begin{table}[htbp]
\centering\footnotesize\setlength{\tabcolsep}{3pt}
\begin{tabular}{llllrrr}
\toprule
Dataset & Area & Response & Expansion & $n_\text{train}$ & $n_\text{test}$ & $p$ \\
\midrule
Diabetes & Healthcare & time in hospital & 47 features, interactions, 250 neighbors & 2{,}500 & 1{,}000 & 283{,}128 \\
Communities and Crime & Socioeconomics & violent crimes per capita & 122 features, interactions, 10 neighbors & 1{,}500 & 250 & 82{,}533 \\
Air Quality & Meteorology & temperature & 12 features, 15 neighbors, interactions & 1{,}000 & 1{,}000 & 18{,}528 \\
\midrule
Communities and Crime & Socioeconomics & violent crimes per capita & 122 features, products of up to 3 & 1{,}500 & 250 & 302{,}743 \\
NHANES & Healthcare & systolic blood pressure & 52 features, products of up to 3 & 2{,}500 & 1{,}000 & 23{,}478 \\
Weather & Meteorology & Basel next-day temp. & 128 features, products of up to 3 & 2{,}500 & 1{,}123 & 100{,}000 \\
\bottomrule
\end{tabular}
\caption{The datasets after feature expansion. Top: the three datasets with neighbor features. Bottom: the three datasets with polynomial features only, whose correlations are weaker.}
\label{tab:datasets}
\end{table}

\subsection{Setup}
\label{subsec:setup}

\paragraph{Modes and configurations.} With $L(\beta) = \|y - X\beta\|_2^2$, the L1+L2 mode minimizes $L(\beta) + \lambda_1\|\beta\|_1 + \lambda_2\|\beta\|_2^2$, the L0+L1+L2 mode adds the constraint $\|\beta\|_0 \le k$, and the L0+L1 and L0+L2 modes keep the constraint and one of the two penalties. In the L0 modes, L0Learn computes paths over a grid of $G$ values of its $L_1$ or $L_2$ penalty and $G$ values of $\lambda_0$, and each solution on the paths defines one ParaDeAl configuration $(k, \lambda_1)$ or $(k, \lambda_2)$. In the L1+L2 mode, glmnet computes a path of $G$ values of its overall penalty strength $\lambda$
for each of $G$ values of its mixing parameter, and ParaDeAl, scikit-learn and skglm solve the same $(\lambda_1, \lambda_2)$ pairs. The L0+L1+L2 mode uses these pairs with $k = 50$. We write $10\times10$, $30\times30$ and $100\times100$ for the grids with $G = 10, 30, 100$.

\paragraph{ParaDeAl.} All configurations of a run are solved together on one L40S GPU in fp32, in chunks sized to a 30\,GB memory budget, and every loop uses Dynamic Shrinking (Section~\ref{subsec:use_dynamic_shrinking}). In the L0 modes, the step size of each configuration comes from the parallel search of Section~\ref{subsec:paradeal_finds_eta}, and support selection runs inside the active set procedure (Section~\ref{subsec:active_set}) with $|A| = \max(250, 5k)$. Support selection stops by the support-unchanged rule and returns the lowest-objective iterate (Section~\ref{subsec:lowest_objective}), and the value modification solvers stop once the relative change in $\beta$ falls below $5\cdot10^{-4}$ in the L0 modes and $10^{-3}$ in the L1+L2 mode. In the L0 modes, value modification is followed by pair swaps (Section~\ref{subsec:pair_swap}) with at most 200 rounds, and the reported runtime includes them. For L1+L2 we use FISTA \citep{beck2009fast}, the accelerated version of proximal gradient, with at most 5{,}000 iterations (note that for these problems unaccelerated proximal gradient with the same budget ends 2--18\% above the competitors' optimization objective). Appendix~\ref{app:setup_details} gives the full settings, with additional ablation studies.

\paragraph{Measurements.} We compare different procedures based on the optimization objective (after standardizing for the different scales used by different implementations, Appendix~\ref{app:competitors}). Additionally we consider the test error $\|y_\text{test} - X_\text{test}\beta\|_2^2$, relative to ParaDeAl.
For each configuration we divide the competitor's optimization objective by ParaDeAl's, and likewise for the test error; we report the mean of each ratio over configurations, so a ratio below 1 means that the competitor is more accurate. On Diabetes, the feature expansion creates columns that are nearly constant on the training rows, so a few test entries become very large after scaling. We therefore compute the test error with each test feature clipped to the range of its training column, which changes only the scoring, not the fitted coefficients (Appendix~\ref{app:data_construction}); on the other datasets, clipping has almost no effect. Runtime is the wall-clock time of the whole solve. For ParaDeAl, the runtime includes the step-size search and both stages. scikit-learn and skglm are too slow to fit the grids in full, so their runtimes are estimated from a random subset of at most 50 configurations of each grid.

When the features are independent, the comparison in the L0 modes is uneventful: on synthetic data with independent features, ParaDeAl matches the objective of L0Learn to within 0.5\% without pair swaps and exactly with them (Appendix~\ref{app:pair_swap_results}). The differences reported below therefore come from correlation between features.

\subsection{Results}
\label{subsec:real_results}

Table~\ref{tab:res_hard} gives the results on the three datasets with strong correlations, and Table~\ref{tab:res_soft} in Appendix~\ref{app:poly_results} gives the same results on the three polynomial datasets. Each block of rows is one dataset and each group of rows within it is one mode, with ParaDeAl in the first row of the group and its competitors below. The three runtime columns give the wall-clock seconds on each grid; for a competitor, the number in parentheses is its runtime divided by ParaDeAl's, so a number above 1 means that ParaDeAl is faster. The last two columns give the competitor's optimization objective and test error divided by ParaDeAl's, averaged over all configurations of the three grids, so a number below 1 means that the competitor is more accurate. In the L0 modes we also show ParaDeAl without pair swaps; for this row the ratios divide L0Learn's values by its own. Panel (b) of each table splits the accuracy ratios of the L0 modes by support size, into $k$ from 1 to 5, 6 to 20, 21 to 50, 51 to 100, and above 100. The L0 modes are meant to select a small number of features that can be interpreted, so the bins up to $k = 100$ are the ones that matter most in practice. The configurations are those of L0Learn, and supports larger than 100 occur only on the $100\times100$ grids, where they make up 23 to 46\% of the configurations.

\begin{table}[p]\centering\scriptsize\setlength{\tabcolsep}{3.5pt}
{\footnotesize (a) Runtime, and accuracy over all configurations}\\[3pt]
\begin{tabular}{llrrrcc}\toprule
 & & \multicolumn{3}{c}{Runtime (s) (speedup)} & \multicolumn{2}{c}{Competitor $\div$ ParaDeAl} \\ \cmidrule(lr){3-5}\cmidrule(lr){6-7}
Mode & Method & 10$\times$10 & 30$\times$30 & 100$\times$100 & Objective ($\uparrow$) & Test error ($\downarrow$) \\
\midrule
\multicolumn{5}{l}{\textbf{Diabetes} ($p = 283{,}128$)} & $n = 2{,}500$ & $n = 1{,}000$ \\
L0+L1 & ParaDeAl (pair swaps) & 0.690 & 6.01 & 794 & 1 & 1 \\
 & ParaDeAl (no swaps) & 0.664 & 3.80 & 112 & 0.637 & 1.30 \\
 & L0Learn (CPU) & 27.5 (39.8$\times$) & 332 (55.3$\times$) & 4{,}398 (5.54$\times$) & 0.964 & 1.01 \\
L0+L2 & ParaDeAl (pair swaps) & 1.81 & 33.3 & 420 & 1 & 1 \\
 & ParaDeAl (no swaps) & 1.75 & 31.8 & 258 & 0.685 & 1.22 \\
 & L0Learn (CPU) & 29.6 (16.3$\times$) & 385 (11.6$\times$) & 4{,}320 (10.3$\times$) & 0.975 & 1.00 \\
L1+L2 & ParaDeAl & 11.3 & 43.7 & 452 & 1 & 1 \\
 & glmnet (CPU) & 179 (15.9$\times$) & 570 (13.1$\times$) & 4{,}094 (9.05$\times$) & 1.000 & 1.00 \\
 & scikit-learn (CPU) & 22{,}101 (1{,}964$\times$) & 216{,}225 (4{,}953$\times$) & 2{,}288{,}215 (5{,}060$\times$) & 1.000 & 1.00 \\
 & skglm (CPU) & 461 (41.0$\times$) & 4{,}128 (94.6$\times$) & 56{,}136 (124$\times$) & 1.000 & 1.00 \\
L0+L1+L2 & ParaDeAl & -- & -- & 588 & -- & -- \\
\midrule
\multicolumn{5}{l}{\textbf{Communities and Crime} ($p = 82{,}533$)} & $n = 1{,}500$ & $n = 250$ \\
L0+L1 & ParaDeAl (pair swaps) & 0.147 & 0.748 & 294 & 1 & 1 \\
 & ParaDeAl (no swaps) & 0.132 & 0.522 & 35.6 & 0.684 & 1.27 \\
 & L0Learn (CPU) & 3.69 (25.1$\times$) & 33.3 (44.5$\times$) & 586 (2.00$\times$) & 0.943 & 1.06 \\
L0+L2 & ParaDeAl (pair swaps) & 0.375 & 4.02 & 109 & 1 & 1 \\
 & ParaDeAl (no swaps) & 0.332 & 3.89 & 63.4 & 0.723 & 1.12 \\
 & L0Learn (CPU) & 6.40 (17.1$\times$) & 45.5 (11.3$\times$) & 563 (5.15$\times$) & 0.945 & 1.01 \\
L1+L2 & ParaDeAl & 1.77 & 8.88 & 93.1 & 1 & 1 \\
 & glmnet (CPU) & 31.5 (17.8$\times$) & 82.3 (9.27$\times$) & 438 (4.70$\times$) & 1.000 & 1.00 \\
 & scikit-learn (CPU) & 1{,}574 (887$\times$) & 13{,}626 (1{,}535$\times$) & 103{,}379 (1{,}110$\times$) & 1.000 & 1.00 \\
 & skglm (CPU) & 66.5 (37.5$\times$) & 483 (54.4$\times$) & 5{,}768 (62.0$\times$) & 1.011 & 1.00 \\
L0+L1+L2 & ParaDeAl & -- & -- & 135 & -- & -- \\
\midrule
\multicolumn{5}{l}{\textbf{Air Quality} ($p = 18{,}528$)} & $n = 1{,}000$ & $n = 1{,}000$ \\
L0+L1 & ParaDeAl (pair swaps) & 0.0929 & 0.243 & 62.1 & 1 & 1 \\
 & ParaDeAl (no swaps) & 0.0544 & 0.0951 & 6.23 & 0.236 & 0.15 \\
 & L0Learn (CPU) & 0.663 (7.14$\times$) & 5.54 (22.8$\times$) & 93.5 (1.51$\times$) & 0.910 & 0.94 \\
L0+L2 & ParaDeAl (pair swaps) & 0.142 & 0.772 & 35.0 & 1 & 1 \\
 & ParaDeAl (no swaps) & 0.0671 & 0.688 & 15.3 & 0.485 & 0.46 \\
 & L0Learn (CPU) & 0.838 (5.92$\times$) & 5.34 (6.92$\times$) & 156 (4.45$\times$) & 0.956 & 1.01 \\
L1+L2 & ParaDeAl & 0.407 & 1.70 & 20.9 & 1 & 1 \\
 & glmnet (CPU) & 3.64 (8.94$\times$) & 11.4 (6.70$\times$) & 82.8 (3.96$\times$) & 0.998 & 1.00 \\
 & scikit-learn (CPU) & 223 (546$\times$) & 2{,}773 (1{,}632$\times$) & 22{,}036 (1{,}055$\times$) & 0.998 & 1.00 \\
 & skglm (CPU) & 12.6 (31.0$\times$) & 134 (78.9$\times$) & 1{,}263 (60.4$\times$) & 0.998 & 1.00 \\
L0+L1+L2 & ParaDeAl & -- & -- & 27.0 & -- & -- \\
\bottomrule\end{tabular}

\vspace{4pt}
{\footnotesize (b) Accuracy by support size $k$ in the L0 modes, pooled over the three grids}\\[2pt]
\includegraphics[width=0.8\textwidth]{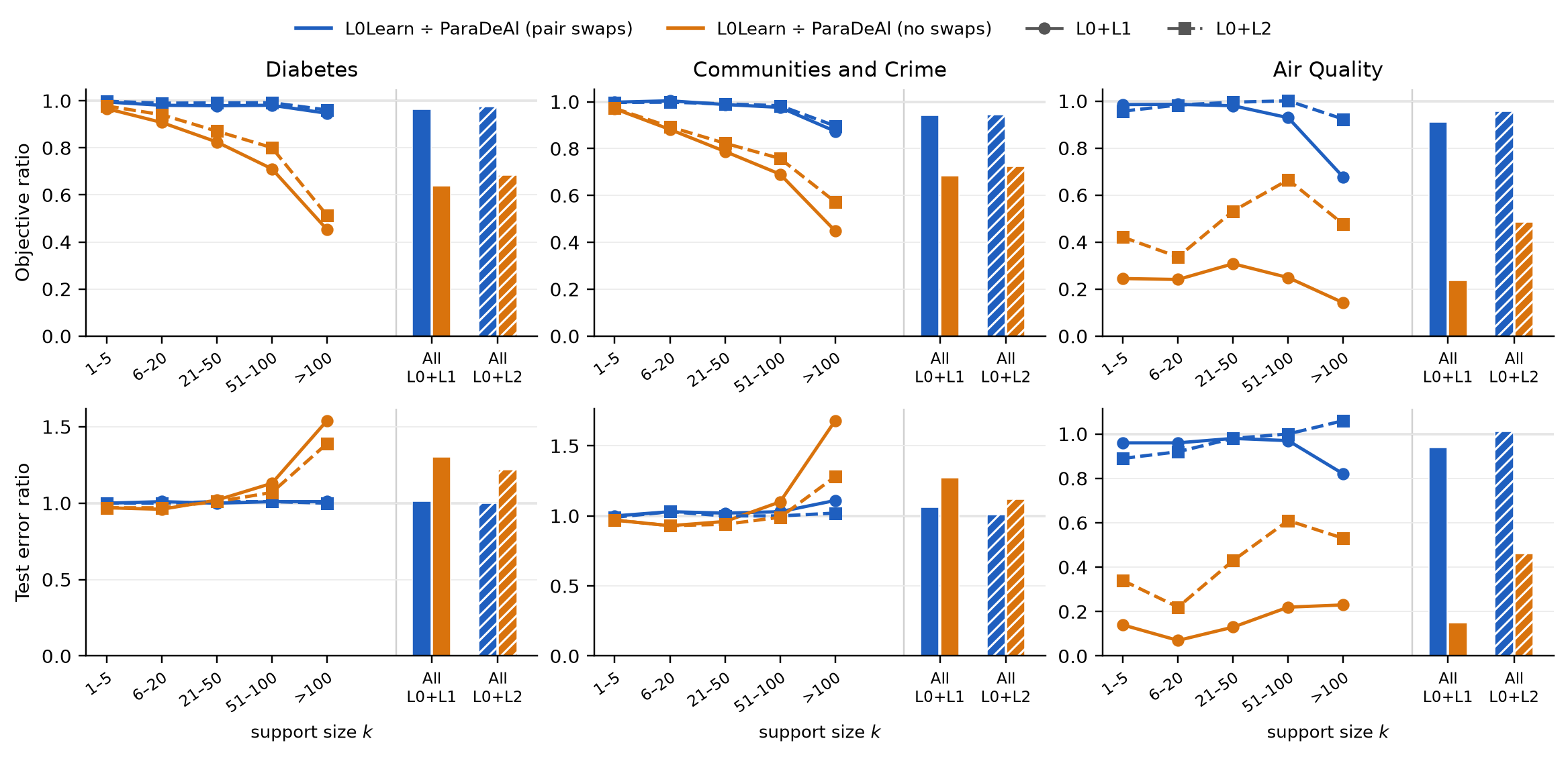}

\caption{Results on the three datasets with strong correlations. ParaDeAl uses the lowest-objective support (Section~\ref{subsec:lowest_objective}) followed by pair swaps (Section~\ref{subsec:pair_swap}) in the L0 modes, and FISTA with tolerance $10^{-3}$ in the L1+L2 mode; its runtimes are the median of three runs. Ratios are averaged over all configurations of the three grids; test errors are computed with test features clipped to their training range. The runtimes of scikit-learn and skglm are estimated from at most 50 configurations per grid. For ParaDeAl (no swaps), the ratios divide L0Learn's values by those of ParaDeAl without pair swaps. Panel (b): L0Learn's objective (top) and test error (bottom) divided by those of ParaDeAl with pair swaps (blue) or without them (orange), so values above 1 favor ParaDeAl. Lines follow the ranges of $k$ (solid: L0+L1, dashed: L0+L2); bars pool all configurations.}
\label{tab:res_hard}
\end{table}

ParaDeAl is faster than every competitor on every grid of every dataset, with one exception: on NHANES in the L0+L1 mode (Appendix~\ref{app:poly_results}), the pair swaps of the $100\times100$ grid take long enough that ParaDeAl and L0Learn finish in the same time. In the L0 modes, panel (b) shows that with pair swaps ParaDeAl matches L0Learn's optimization objective on sparse supports. Pooled over all configurations with $k \le 100$, the ratio is between 0.976 and 1.006, and within each bin up to $k = 100$ it is between 0.98 and 1.01, except on Air Quality and in a few cases with $k$ between 51 and 100 (down to 0.93). The lower pooled ratios of panel (a), 0.91 to 0.99, come from supports larger than 100. For $k \le 100$, L0Learn's test error divided by ParaDeAl's is between 0.89 and 1.07. In the L1+L2 mode, ParaDeAl with FISTA matches the optimization objective of all three competitors. The L0+L1+L2 mode has no competitor; ParaDeAl solves its $100\times100$ grid in 588\,s on Diabetes, 135\,s on Crime and 27\,s on Air Quality, and in 502\,s, 52\,s and 210\,s on the polynomial versions of Crime, NHANES and Weather.

\subsection{Summary of findings}

\paragraph{L0+L1 and L0+L2 against L0Learn: much faster, and as accurate on sparse supports once pair swaps are used.} On the $10\times10$ and $30\times30$ grids, ParaDeAl is 5.9--268$\times$ faster than L0Learn. On the $100\times100$ grids the advantage shrinks to 1.5--10$\times$ (and to a tie on NHANES in the L0+L1 mode), because these grids also contain supports of up to 500 features, on which pair swaps need many rounds; without pair swaps, ParaDeAl is 8.9--39$\times$ faster there.

The breakdown by support size in panel (b) of Table~\ref{tab:res_hard}, and of Table~\ref{tab:res_soft} in Appendix~\ref{app:poly_results}, shows where accuracy is gained and lost. Without pair swaps, ParaDeAl falls further behind L0Learn as $k$ grows. On Crime in the L0+L1 mode, for example, the objective ratio is 0.970 for $k$ from 1 to 5, 0.882 from 6 to 20, 0.788 from 21 to 50, 0.691 from 51 to 100 and 0.449 above 100, and every other dataset shows the same pattern. Pair swaps remove almost all of this gap on the sparse supports that the L0 modes are meant for. For $k$ up to 50, the ratio is between 0.98 and 1.01 on every dataset, except Air Quality in the L0+L2 mode with $k \le 5$ (0.957). For $k$ from 51 to 100 it is between 0.96 and 1.00, except Air Quality in the L0+L1 mode (0.929). A clear gap remains only above $k = 100$ (0.68 to 0.97), where the limit of 200 rounds stops pair swaps before they converge. Since these large supports make up a quarter to a half of the $100\times100$ grids, they also explain the lower pooled ratios of panel (a) (0.91 to 0.99). What remains of the gap is largest on Air Quality, in which every feature has a near-copy, and smallest on NHANES and Weather, the two datasets with the weakest correlations. This agrees with the synthetic experiments of Section~\ref{subsec:pair_swap}, where the gap appears only when features are correlated.

On held-out data the two methods are about equally accurate on sparse supports. With pair swaps and $k \le 100$, the L0Learn/ParaDeAl test error ratio is between 0.89 and 1.07, and it is below 0.95 only on Air Quality in the L0+L2 mode with $k \le 20$. Without pair swaps, ParaDeAl's test error is lower than L0Learn's on Diabetes and on both Crime datasets (pooled ratios 1.12--1.33) despite its worse objective, and this advantage grows with $k$, so part of L0Learn's lower objective on large supports comes with overfitting.

The two methods usually reach these similar objectives with different features. With pair swaps, ParaDeAl returns exactly L0Learn's support for up to 83\% of the configurations with $k$ from 1 to 5 (on NHANES), but for at most 26\% of the configurations with larger $k$, and for almost none on Diabetes (Appendix~\ref{app:support_agreement}). When features are strongly correlated, many different supports fit almost equally well, so this is to be expected.

\paragraph{L1+L2: FISTA matches the competitors and is faster than all of them.} glmnet, scikit-learn and skglm agree with each other, and ParaDeAl with FISTA matches their training objective to within 0.3\% on every dataset (pooled ratios 0.998--1.000; skglm on Crime is the one exception, at 1.010 in ParaDeAl's favor). With FISTA, ParaDeAl is faster than glmnet on every grid, by 4.0--17.8$\times$. The advantage is largest on the small grids and on Diabetes, and smallest on the large grids of Crime and Air Quality, where glmnet's warm-started path does very little work per additional $\lambda$.  Against the methods that solve one configuration at a time, the advantage is much larger: 31--124$\times$ against skglm and 546--5{,}060$\times$ against scikit-learn, at the same accuracy. On the polynomial datasets (Table~\ref{tab:res_soft} in Appendix~\ref{app:poly_results}) the picture is the same: FISTA matches all three competitors (pooled ratios 0.999--1.016), and it is 3.2--52$\times$ faster than glmnet, 23--320$\times$ faster than skglm and 744--1{,}829$\times$ faster than scikit-learn.

\paragraph{Storage.} On the small grids ParaDeAl uses less memory than L0Learn and glmnet. On the large grids ParaDeAl uses more memory than glmnet and, in most cases, than L0Learn; its GPU peak there is bounded by the 30 GB working budget it is given for chunking configurations, and a smaller budget would trade memory for time. L0Learn peaks at 0.2--30 GB of RAM, glmnet at 0.4--18 GB and skglm below 0.1 GB.

\section{Summary and Further Work} \label{sec:summary}

We present ParaDeAl, a GPU-based framework that solves large collections of hyperparameter configurations of the L0-constrained regression problem~\eqref{eq:our_estimator} simultaneously, typically within minutes. ParaDeAl supports L1 and L2 penalties, a hard L0 constraint, and both grid and tuple search. Potential applications include:

\begin{itemize}
    \item \textbf{Sparsity--accuracy tradeoff in ML modeling.} The level of compression that can be achieved in weights or activations depends on the acceptable loss of accuracy. An informed decision requires the objective as a function of $k$; computing it for many values of $k$ is expensive in both time and memory, even with large amounts of parallel hardware. A parallel algorithm that solves hundreds of values of $k$ at once makes it quick to locate the `elbow' of this curve, which can then be refined further, so choosing the sparsity level becomes a tractable problem.
    \item \textbf{Portfolio optimization.} Many financial models have the same structure as Problem~\eqref{eq:our_estimator}: a loss function with regularization, where the goal is to select a small number of assets from a large universe. Problems of this form can be solved with ParaDeAl. Variants include latent variable models, most commonly probit or logit models, in which additional variables represent underlying traits of the assets. With modification, ParaDeAl may also be applied to hyperparameter tuning for gradient-boosted machines (GBMs), whose accuracy depends strongly on their hyperparameters and for which little theory on hyperparameter selection exists.
    \item \textbf{Exploring nonconvex landscapes.} In high-dimensional nonconvex ML, a central challenge is to ensure that the local minimum returned by GD is the global minimum, or provably close to it. Many ML objective functions have complicated landscapes. Running many starting points and step sizes at once gives a quick map of the local minima that descent methods reach, which can help characterize these landscapes empirically.
\end{itemize}

Extensions of ParaDeAl include:

\begin{itemize}
    \item \textbf{Better procedures for the hard thresholding step.} From an implementation perspective, \verb|torch.vmap| does not support a varying $k$, and the existing alternatives are a binary overlay or a \texttt{for} loop, so a more efficient parallel algorithm would be beneficial. From an algorithmic perspective, top-$k$ selection by magnitude may not be the best sparsification rule. In ongoing work, we find that sparsifying a fixed proportion per row or per column improves accuracy over plain magnitude pruning for a set of reconstruction objectives. A data-aware rule in the spirit of Wanda \citep{sun24}, which scores coefficients using $X$ as well as $\beta$, is also parallelizable and may improve accuracy at some cost in runtime.
    \item \textbf{Incorporating quantization into IHT.} The default implementation of ParaDeAl uses single precision (fp32). For very large $X$ and $y$, it may be necessary to store them in reduced precision to fit them in GPU memory.
    \item \textbf{Sharding and distributed matrix computations.} Another option for very large $X$ is to split the matrix operations into smaller chunks, which can be distributed across several GPUs.
    \item \textbf{Adapting the learning rate $\eta$ during the iterations,} for example by a per-configuration backtracking rule, which would also address the step-size failures described in Section~\ref{sec:compare_results}.
\end{itemize}

\bibliography{sources}

\vspace{20pt}

\begin{appendices}

\noindent {\Large\bfseries Appendices}

\section{Using $(X^TX)\beta$ Versus $X^T(X\beta)$} \label{app:xtx_or_xbeta}

Each gradient step~\eqref{eq:1_grad_step} requires the product $X^TX\beta$ for every configuration. It can be computed either as $(X^TX)\beta$, where the Gram matrix $G = X^TX$ is formed once before the iterations and reused, or as $X^T(X\beta)$, where both products are recomputed at every iteration. We refer to these as the Gram form and the two-product form. This appendix compares their costs and reports timings on an NVIDIA L40S.

\subsection*{Cost analysis}

Let $B$ be the $(\prod C \times p)$ matrix whose rows are the $\beta$s of all configurations, and let $I$ be the number of iterations. Multiplying an $(a \times b)$ matrix by a $(b \times c)$ matrix takes $abc$ multiply--add operations. The Gram form computes $G$ once, at a cost of $np^2$ operations, and then computes $BG$ at every iteration, at a cost of $\prod C \cdot p^2$. The two-product form computes $BX^T$ and then $(BX^T)X$ at every iteration, at a cost of $2\prod C \cdot np$. The additional memory, beyond $X$ and $B$, is $p^2$ entries for $G$ plus $\prod C \cdot p$ for the output in the Gram form, and $\prod C \cdot n$ for the intermediate $BX^T$ plus $\prod C \cdot p$ for the output in the two-product form. Table~\ref{tab:xtx_vs_xbeta_runtime_storage} summarizes these costs.

\begin{table}[htbp]
\centering
\begin{tabular}{|c|c|c|}
    \hline
    Method & Operations over $I$ iterations & Additional memory \\ \hline
    $(X^TX)\beta$ & $np^2 + I \prod C \cdot p^2$ & $p^2 + \prod C \cdot p$ \\ \hline
    $X^T(X\beta)$ & $2 I \prod C \cdot np$ & $\prod C \cdot (n + p)$ \\ \hline
\end{tabular}
\caption{Arithmetic and memory costs of the two ways of computing $X^TX\beta$ for $\prod C$ configurations over $I$ iterations. The first term for $(X^TX)\beta$ is the one-time cost of forming $X^TX$.}
\label{tab:xtx_vs_xbeta_runtime_storage}
\end{table}

Comparing the operation counts, the Gram form is cheaper if and only if $np^2 + I \prod C \cdot p^2 < 2 I \prod C \cdot np$. Dividing by $p$ and rearranging gives
\begin{equation}
p < p^\ast(n) := \frac{2n \cdot I \prod C}{n + I \prod C}. \label{eq:xtx_threshold}
\end{equation}
When $I \prod C \gg n$, the cost of forming $G$ is amortized and $p^\ast \approx 2n$: the Gram form is preferable when $p < 2n$, and the two-product form when $p > 2n$. When $I \prod C$ is small relative to $n$, the cost of forming $G$ is not recovered and $p^\ast \approx 2 I \prod C$. The memory comparison is independent of $I$: the Gram form needs $p^2$ additional entries, which is prohibitive for large $p$ (about 40 GB in fp32 at $p = 10^5$), while the two-product form needs $\prod C \cdot n$, which grows with the number of configurations.

The operation count describes the cost on a GPU only when the products are large enough to keep the device busy. For small products, the runtime is dominated by the overhead of launching each kernel and by memory traffic, and the two-product form launches twice as many kernels per iteration. The count also ignores caching: a Gram matrix small enough to stay in the GPU's cache is cheaper to read than its size suggests. We therefore measure the two forms directly.

\subsection*{Experiments}

We time the product alone, without the other gradient terms, for $I = 1{,}000$ iterations, with 

\begin{figure}[htbp]
    \centering
    \begin{subfigure}{0.48\textwidth}
        \centering
        \includegraphics[width=\textwidth]{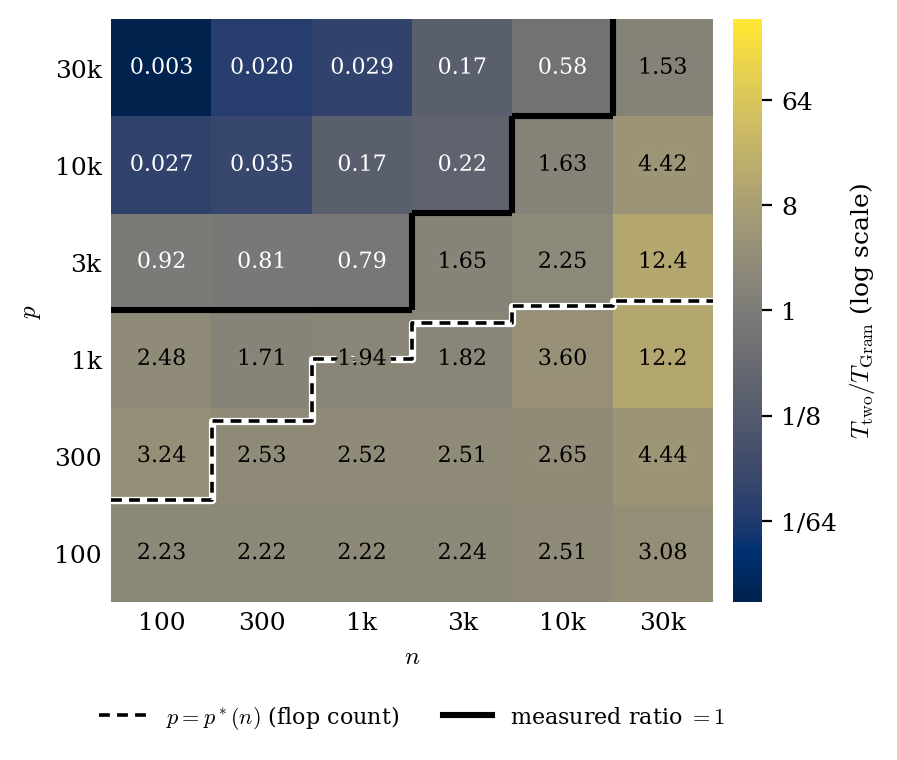}
        \caption{$\prod C = 1$.}
    \end{subfigure}
    \hfill
    \begin{subfigure}{0.48\textwidth}
        \centering
        \includegraphics[width=\textwidth]{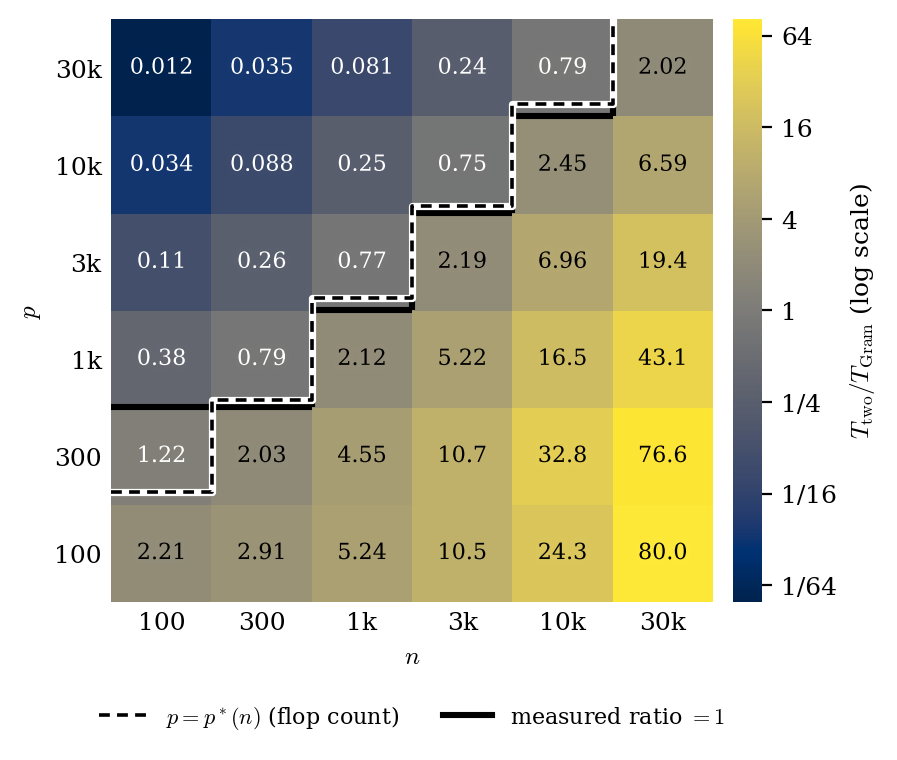}
        \caption{$\prod C = 1{,}000$.}
    \end{subfigure}
    \caption{Runtime of $X^T(X\beta)$ divided by runtime of $(X^TX)\beta$ over 1{,}000 iterations, including the time to form $X^TX$, on one NVIDIA L40S. Values above 1 favor $(X^TX)\beta$. The solid line marks where the measured ratio crosses 1, and the dashed line marks the threshold $p = p^\ast(n)$ of~\eqref{eq:xtx_threshold}.}
    \label{fig:xtx_vs_xbeta}
\end{figure}

$n, p \in \{100, 300, 1{,}000, 3{,}000, 10{,}000, 30{,}000\}$, and with $\prod C = 1$ and $\prod C = 1{,}000$ configurations. The entries of $X$ are standard Gaussian and its columns are standardized to mean 0 and standard deviation 1. All computation is in fp32 (TF32 disabled) on one NVIDIA L40S (48 GB) with PyTorch 2.7.1 and CUDA 12.6. Times are medians of 3 repetitions, with a device synchronization before each reading. The time of the Gram form includes forming $G$. The two forms agree to a relative difference of at most $1.5 \times 10^{-5}$. Figure~\ref{fig:xtx_vs_xbeta} reports the ratio of the time of the two-product form to the time of the Gram form, so values above 1 favor the Gram form.

\paragraph{Many configurations.} With $\prod C = 1{,}000$, the operation count predicts the faster form well. Along $p \approx n$ the measured ratio is 2.0--2.7, against a predicted value of 2. The measured crossover lies 1.1 to 1.4 times above $p^\ast(n)$ for $n \geq 300$, so the Gram form wins slightly more often than the operation count predicts; the largest deviation is at $n = 100$ (1.8 times). Additional measurements at $\prod C = 100$ and $\prod C = 10{,}000$ for $n \in \{300, 1{,}000, 3{,}000\}$ place the crossover between 0.96 and 1.41 times $p^\ast(n)$, except at $\prod C = 100$, $n = 300$, where the Gram form still wins at $p = 1{,}000$ ($p^\ast = 598$). Forming $G$ costs at most 2.4\% of the total time of the Gram form at $\prod C = 1{,}000$ (at most 1.4~s, at $n = p = 30{,}000$). When both forms are limited by arithmetic, their throughput is similar (about 37 and 36 TFLOP/s at $n = p = 30{,}000$).

\paragraph{A single configuration.} With $\prod C = 1$, every iteration is a matrix--vector product, and the operation count does not predict the outcome. For $p \leq 1{,}000$ the ratio is at least about 2.2 for every $n$, because each iteration is dominated by kernel launch overhead and the two-product form launches two kernels instead of one. For $p \leq 3{,}000$, $G$ (at most 36 MB) fits in the 96 MB L2 cache of the L40S while $X$ may not, which makes the Gram form up to 12 times faster at large $n$ where the operation count predicts the opposite. For $p \geq 10{,}000$, $G$ (0.4 to 3.6 GB) must be read from memory at every iteration, and the Gram form is much slower when $n$ is small (the ratio is 0.003 at $n = 100$, $p = 30{,}000$). The measured crossover lies 2.3 to 15 times above $p^\ast(n)$, and at $n = 30{,}000$ the Gram form wins for every $p$ we tested.

\paragraph{Memory.} The measured additional memory matches Table~\ref{tab:xtx_vs_xbeta_runtime_storage} to within about 1 MB in every case. At $\prod C = 1$ and $p = 30{,}000$, the Gram form needs 3.6 GB and the two-product form 0.24 MB. At $\prod C = 10{,}000$, $n = 3{,}000$ and $p = 3{,}000$, the order is reversed (156 MB against 240 MB).

\paragraph{Recommendation.} For the batch sizes for which ParaDeAl is intended ($\prod C$ in the hundreds or more), we use the Gram form when $p < p^\ast(n) \approx 2n$ and $X^TX$ fits in GPU memory, and the two-product form otherwise. In particular, all experiments in this paper with $p \geq 10^4$ and $p \gg n$ use $X^T(X\beta)$. For very small batches the choice should be made by measurement, because launch overhead and caching dominate.

\section{Proofs} \label{app:proofs}

\subsection*{Proof of Proposition~\ref{prop:thresholding_form}}

\begin{proof}
For a set $\mathcal{S} \subseteq [p]$ with $|\mathcal{S}| \leq k$, consider the objective of~\eqref{eq:t_thresholding} over vectors $\beta$ with $\beta_i = 0$ for $i \notin \mathcal{S}$. It is separable across coordinates. A coordinate $i \notin \mathcal{S}$ contributes $\frac{1}{2}u_i^2$. For a coordinate $i \in \mathcal{S}$, the scalar problem $\min_{b \in \mathbb{R}} \frac{1}{2}(b - u_i)^2 + \lambda |b|$ is solved by $b = \sign(u_i) \max(|u_i| - \lambda, 0)$, with optimal value $\frac{1}{2}u_i^2 - \frac{1}{2}\big(\max(|u_i| - \lambda, 0)\big)^2$. Hence the optimal value of~\eqref{eq:t_thresholding} restricted to $\mathcal{S}$ is
\begin{equation*}
\frac{1}{2}\|u\|_2^2 - \frac{1}{2}\sum_{i \in \mathcal{S}} g(|u_i|), \qquad g(t) = \big(\max(t - \lambda, 0)\big)^2 .
\end{equation*}
Since every vector with $\|\beta\|_0 \leq k$ is supported on some $\mathcal{S}$ with $|\mathcal{S}| \leq k$, minimizing~\eqref{eq:t_thresholding} amounts to choosing $\mathcal{S}$ with $|\mathcal{S}| \leq k$ that maximizes $\sum_{i \in \mathcal{S}} g(|u_i|)$. The function $g$ is nonnegative and nondecreasing on $[0, \infty)$, so a maximizer is given by the indices of the $k$ largest values of $|u_i|$, that is, $\mathcal{S} = \mathcal{R}_k(u)$ (adding an index never decreases the sum). Setting the coordinates in $\mathcal{R}_k(u)$ to their scalar minimizers and the others to zero gives $\hat{\beta}$, which therefore belongs to $T(u, k, \lambda)$. When several indices tie in $|u_i|$ at the $k$th position, any choice among them gives the same value, which is why $T$ can be set-valued.
\end{proof}

\subsection*{Monotonicity of hard-soft thresholding}

\begin{proposition} \label{prop:descent}
Let $\lambda_1, \lambda_2 \geq 0$, let $L' \geq L_D = 2(\sigma_{\max}^2(X) + \lambda_2)$ and $\eta = 1/L'$, and let $\beta_t$ satisfy $\|\beta_t\|_0 \leq k$. If $\beta_{t+1} \in T(\beta_t - \eta \nabla F_D(\beta_t), k, \eta\lambda_1)$, then
\begin{equation*}
F(\beta_{t+1}) \leq F(\beta_t) - \frac{L' - L_D}{2} \|\beta_{t+1} - \beta_t\|_2^2 \leq F(\beta_t).
\end{equation*}
\end{proposition}

\begin{proof}
Write $g = \nabla F_D(\beta_t)$ and $\Delta = \beta_{t+1} - \beta_t$. By~\eqref{eq:deriving_final_opt_problem}--\eqref{eq:final_opt_problem}, $\beta_{t+1}$ minimizes
\begin{equation*}
Q(\beta) = \langle g, \beta - \beta_t \rangle + \frac{L'}{2}\|\beta - \beta_t\|_2^2 + \lambda_1 \|\beta\|_1
\end{equation*}
over $\{\beta : \|\beta\|_0 \leq k\}$, since $Q$ differs from $L'$ times the objective of~\eqref{eq:t_thresholding}, with $u = \beta_t - \eta g$ and $\lambda = \eta\lambda_1$, only by a constant. Because $\beta_t$ is feasible, $Q(\beta_{t+1}) \leq Q(\beta_t) = \lambda_1 \|\beta_t\|_1$. By the descent lemma~\eqref{eq:descent_lemma} with constant $L_D$,
\begin{align*}
F(\beta_{t+1}) &= F_D(\beta_{t+1}) + \lambda_1\|\beta_{t+1}\|_1 \leq F_D(\beta_t) + \langle g, \Delta \rangle + \frac{L_D}{2}\|\Delta\|_2^2 + \lambda_1\|\beta_{t+1}\|_1 \\
&= F_D(\beta_t) + Q(\beta_{t+1}) - \frac{L' - L_D}{2}\|\Delta\|_2^2 \leq F(\beta_t) - \frac{L' - L_D}{2}\|\Delta\|_2^2 . \qedhere
\end{align*}
\end{proof}

\subsection*{Proof of Proposition~\ref{prop:converge}}

Proposition~\ref{prop:converge_precise} below is the full version of Proposition~\ref{prop:converge}. It adapts Proposition~6 of \citet{bkm16}, who study hard-thresholded gradient descent, the same update without the $L_1$ penalty, for a convex loss with a Lipschitz gradient. Our update adds soft thresholding to every step, so the objective $F$ is no longer differentiable, and on a fixed support the iteration is a proximal gradient method rather than gradient descent. Parts (a), (b) and (d) below adapt Parts (a), (b) and (c) of \citet[Proposition~6]{bkm16}, and Part (c), on limit points, is added. Parts (d) and (e) of \citet[Proposition~6]{bkm16} treat the case in which the $k$th largest coefficient tends to zero; the gradient of the loss then vanishes along a subsequence, so the iterates approach unconstrained minimizers of the loss. \citet{bkm16} regard this case as of little statistical interest, and we do not adapt those parts.

Throughout, $L_D = 2(\sigma_{\max}^2(X) + \lambda_2)$ is the Lipschitz constant of $\nabla F_D$ in~\eqref{eq:lipschitz_inequality}, and $S(\beta)$ is the support of $\beta$ (Section~\ref{subsec:notation}). For $t \geq 1$, let $\alpha_t$ denote the $k$th largest entry of $|\beta_t|$, so that $\alpha_t > 0$ if and only if $\|\beta_t\|_0 = k$. For $\mathcal{S} \subseteq [p]$, let $V_\mathcal{S} = \{\beta \in \mathbb{R}^p : \beta_i = 0 \text{ for } i \notin \mathcal{S}\}$.

\begin{proposition} \label{prop:converge_precise}
Let $\lambda_1, \lambda_2 \geq 0$, $L' > L_D$ and $\eta = 1/L'$. Let $\beta_0 \in \mathbb{R}^p$ be arbitrary and $\beta_{t+1} \in T(\beta_t - \eta \nabla F_D(\beta_t), k, \eta\lambda_1)$ for $t \geq 0$. Then the following hold.
\begin{enumerate}
\item[(a)] For every $t \geq 1$, $F(\beta_t) - F(\beta_{t+1}) \geq \frac{L' - L_D}{2}\|\beta_{t+1} - \beta_t\|_2^2$. The sequence $F(\beta_t)$, $t \geq 1$, is nonincreasing and converges.
\item[(b)] $\sum_{t \geq 1} \|\beta_{t+1} - \beta_t\|_2^2 \leq 2F(\beta_1)/(L' - L_D)$. In particular, $\beta_{t+1} - \beta_t \to 0$.
\item[(c)] Every limit point $\bar\beta$ of $\{\beta_t\}$ is a fixed point of the update: $\bar\beta \in T(\bar\beta - \eta\nabla F_D(\bar\beta), k, \eta\lambda_1)$.
\item[(d)] If $\liminf_{t \to \infty} \alpha_t > 0$, then the support stabilizes after finitely many iterations: there exist a set $\mathcal{S}$ with $|\mathcal{S}| = k$ and an index $M$ such that $S(\beta_t) = \mathcal{S}$ for all $t \geq M$. Moreover, $\beta_t$ converges to a point $\beta^\ast$ with $S(\beta^\ast) = \mathcal{S}$ that minimizes $F$ over $V_\mathcal{S}$ and is a fixed point of the update, as in Part (c).
\end{enumerate}
\end{proposition}

\begin{proof}
Write $u_t = \beta_t - \eta\nabla F_D(\beta_t)$ and $\phi_t(\beta) = \frac{1}{2}\|\beta - u_t\|_2^2 + \eta\lambda_1\|\beta\|_1$. By definition of $T$ in~\eqref{eq:t_thresholding},
\begin{equation} \label{eq:T_optimality}
\phi_t(\beta_{t+1}) \leq \phi_t(\beta) \quad \text{for every } \beta \text{ with } \|\beta\|_0 \leq k .
\end{equation}

\emph{Part (a).} The starting point $\beta_0$ need not satisfy $\|\beta_0\|_0 \leq k$, but every $\beta_t$ with $t \geq 1$ is an element of $T(\cdot, k, \cdot)$ and therefore does. Proposition~\ref{prop:descent} applies to every $t \geq 1$ and gives the inequality. Since $F \geq 0$, the nonincreasing sequence $F(\beta_t)$ is bounded below and converges.

\emph{Part (b).} Summing the inequality of Part (a) over $t = 1, \dots, N$ gives $\frac{L' - L_D}{2}\sum_{t=1}^{N}\|\beta_{t+1} - \beta_t\|_2^2 \leq F(\beta_1) - F(\beta_{N+1}) \leq F(\beta_1)$. Since $L' > L_D$, dividing by $(L' - L_D)/2$ and letting $N \to \infty$ proves the bound. The terms of a convergent series tend to zero.

\emph{Part (c).} Let $\beta_{t_j} \to \bar\beta$. By Part (b), $\beta_{t_j + 1} \to \bar\beta$ as well, and since $\nabla F_D$ is continuous, $u_{t_j} \to \bar u := \bar\beta - \eta\nabla F_D(\bar\beta)$. Both sides of~\eqref{eq:T_optimality} are continuous in $(\beta_{t+1}, u_t)$, so taking $t = t_j$ and $j \to \infty$ gives $\frac{1}{2}\|\bar\beta - \bar u\|_2^2 + \eta\lambda_1\|\bar\beta\|_1 \leq \frac{1}{2}\|\beta - \bar u\|_2^2 + \eta\lambda_1\|\beta\|_1$ for every $\beta$ with $\|\beta\|_0 \leq k$. The set $\{\beta : \|\beta\|_0 \leq k\}$ is closed, so $\|\bar\beta\|_0 \leq k$, and therefore $\bar\beta \in T(\bar u, k, \eta\lambda_1)$.

\emph{Part (d).} Let $\underline{\alpha} = \liminf_{t} \alpha_t > 0$, and choose $M_0$ such that $\alpha_t \geq \underline{\alpha}/2$ for all $t \geq M_0$. For these $t$, $\|\beta_t\|_0 = k$. Suppose that $S(\beta_{t+1}) \neq S(\beta_t)$ for some $t \geq M_0$. Both supports have $k$ elements, so some index $i$ lies in $S(\beta_t)$ but not in $S(\beta_{t+1})$, and
\begin{equation*}
\|\beta_{t+1} - \beta_t\|_2 \geq |\beta_{t,i}| \geq \alpha_t \geq \underline{\alpha}/2 .
\end{equation*}
By Part (b), this can happen for only finitely many $t$. Hence there exist $M \geq M_0$ and $\mathcal{S}$ with $|\mathcal{S}| = k$ such that $S(\beta_t) = \mathcal{S}$ for all $t \geq M$. This step follows the proof of Part (c) of \citet[Proposition~6]{bkm16}.

For $t \geq M$, both $\beta_t$ and $\beta_{t+1}$ lie in $V_\mathcal{S}$. Every vector in $V_\mathcal{S}$ has at most $k$ nonzeros, so by~\eqref{eq:T_optimality} the point $\beta_{t+1}$ also minimizes $\phi_t$ over $V_\mathcal{S}$. On $V_\mathcal{S}$ the function $\phi_t$ is strongly convex and separable across coordinates, so its minimizer is unique and is given by $(\beta_{t+1})_i = \sign(u_{t,i})\max(|u_{t,i}| - \eta\lambda_1, 0)$ for $i \in \mathcal{S}$ and $(\beta_{t+1})_i = 0$ otherwise. This is one step of the proximal gradient method, with step size $\eta$, applied to the convex problem
\begin{equation} \label{eq:restricted_problem}
\min_{\beta \in V_\mathcal{S}} \ F_D(\beta) + \lambda_1\|\beta\|_1 ,
\end{equation}
whose smooth part $F_D$ has an $L_D$-Lipschitz gradient. Thus, from iteration $M$ onward, hard-soft thresholding coincides with the proximal gradient method (ISTA) for~\eqref{eq:restricted_problem}, started from $\beta_M$.

Problem~\eqref{eq:restricted_problem} has a minimizer. If $\lambda_1 > 0$ or $\lambda_2 > 0$, its objective is continuous and tends to infinity as $\|\beta\|_2 \to \infty$. If $\lambda_1 = \lambda_2 = 0$, the problem is least squares on the columns $X_\mathcal{S}$, whose normal equations $X_\mathcal{S}^TX_\mathcal{S}\beta_\mathcal{S} = X_\mathcal{S}^Ty$ always have a solution. For a convex problem that has a minimizer, the iterates of the proximal gradient method with a constant step size in $(0, 2/L_D)$ converge to a minimizer \citep{com05}; see also \citet[Chapter~10]{bec17}. Since $\eta = 1/L' < 1/L_D$, the sequence $\beta_t$ converges to a minimizer $\beta^\ast$ of~\eqref{eq:restricted_problem}, which is the same as a minimizer of $F$ over $V_\mathcal{S}$. For $i \in \mathcal{S}$ and $t \geq M$, $|\beta_{t,i}| \geq \alpha_t \geq \underline{\alpha}/2$, so $|\beta^\ast_i| \geq \underline{\alpha}/2 > 0$ and $S(\beta^\ast) = \mathcal{S}$. Finally, $\beta^\ast$ is a limit point of $\{\beta_t\}$, so Part (c) shows that $\beta^\ast$ is a fixed point of the update.
\end{proof}

\paragraph{Remarks.} (1) Part (d) also shows that, under its condition, a stopping rule based on the support, such as the support-unchanged rule of Section~\ref{subsubsec:support_selection}, is met after finitely many updates. The condition $\liminf_t \alpha_t > 0$ concerns the iterates and cannot be checked in advance.

(2) The convergence of the whole sequence $\{\beta_t\}$ can also be obtained without the condition of Part (d), at the cost of more advanced tools. The update~\eqref{eq:iterative_T} is a forward--backward (proximal gradient) step for the sum of $F_D$ and the nonconvex function $\lambda_1\|\beta\|_1 + \delta_k(\beta)$, where $\delta_k$ is $0$ on $\{\beta : \|\beta\|_0 \leq k\}$ and $+\infty$ elsewhere. This sum is semi-algebraic, so it has the Kurdyka--{\L}ojasiewicz property, and the results of \citet{att13} on forward--backward splitting for such functions, which also cover iterative thresholding for sparse recovery, show that a bounded sequence of this kind converges to a critical point, in the sense of the limiting subdifferential. This conclusion concerns the iterates only: when the limit has fewer than $k$ nonzeros, the support need not stabilize. The sequence is bounded when $\lambda_1 > 0$ or $\lambda_2 > 0$, because the level sets of $F$ are then bounded. When $\lambda_1 = \lambda_2 = 0$, the sequence is bounded whenever every set of $k$ columns of $X$ is linearly independent: there are finitely many supports of size at most $k$, the least squares objective restricted to each of them is coercive, and $F(\beta_t) \leq F(\beta_1)$ for $t \geq 1$.

(3) The same argument as in Proposition~\ref{prop:descent}, with the constraint set $\{\beta : \|\beta\|_0 \leq k\}$ replaced by the convex set $\{\beta : \beta_i = 0, \ i \notin \mathcal{K}\}$, shows that the value modification update~\eqref{eq:iterative_vm} is also monotone when its step size is at most $1/L_D$. In that convex case, $F(\beta_t)$ converges to the optimal value of~\eqref{eq:vm_equation} \citep{beck2009fast}, and, for any constant step size in $(0, 2/L_D)$, $\beta_t$ converges to a minimizer \citep{com05}.

(4) Proposition~\ref{prop:converge_precise} requires $\eta < 1/L_D$. Step sizes larger than $1/L_D$, such as the steps selected by the search of Section~\ref{subsec:paradeal_finds_eta}, are not covered by Propositions~\ref{prop:descent} and~\ref{prop:converge_precise}.

\section{Number of Unchanged Updates in the Stopping Rule} \label{app:stability_r}

Support selection stops once one update leaves the support unchanged (Section~\ref{subsubsec:support_selection}). Here we examine the more general rule that stops once $r$ consecutive updates have left the support unchanged, for $r \in \{1, 2, 3, 5, 10\}$. The rule remains support-only for every $r$. We judge each $r$ by the support it returns, after value modification on that support.

\paragraph{Setup.} We use $(n, p) \in \{(500, 5{,}000), (1{,}000, 20{,}000)\}$. The rows of $X$ are Gaussian with correlation $\rho^{|i-j|}$ between features $i$ and $j$, $\rho \in \{0, 0.5\}$, and every column is standardized to mean 0 and standard deviation 1. The true coefficient vector has 25 nonzeros at random positions with standard Gaussian values, and $y = X\beta^\star + \epsilon$ with signal-to-noise ratio 5. A test set of the same size is drawn from the same model. Each of 3 seeds uses a grid of 180 configurations: $k \in \{10, 25, 50, 100\}$; $\lambda_1 \in \{0, 10^{-3}, 10^{-2}\} \times 2\|X^Ty\|_\infty$; $\lambda_2 \in \{0, 10^{-3}, 10^{-2}\} \times \sigma_{\max}^2(X)$; and the five starting points of Section~\ref{subsec:benefits}. We use three step sizes: the conservative $\eta = 0.9/(\sigma_{\max}^2(X) + \lambda_2)$ of Section~\ref{subsec:resolve_cmtp}, and $\eta_{\text{pred}}/4$ and $\eta_{\text{pred}}/2$, where $\eta_{\text{pred}} = 0.9/(n(1+\sqrt{k/n})^2 + \lambda_2)$ is the restricted prediction of Section~\ref{subsec:active_set_results}. (At $\eta_{\text{pred}}$ itself almost every configuration diverged, as expected from the probe results of Section~\ref{subsec:active_set_results}, which select between $1/8$ and $1/2$ of $\eta_{\text{pred}}$.) At $\eta_{\text{pred}}/2$, between 49 and 323 of the 540 configurations per instance type diverged, mostly with $k = 100$; they are excluded. No configuration diverged at the other two step sizes. The cap is $M_{SS} = 1{,}000$. Value modification uses PCG when $\lambda_1 = 0$ and proximal gradient otherwise, both to a relative tolerance of $10^{-6}$. Computation is in fp32 on one NVIDIA L40S.

\begin{table}[htbp]
\centering
\small
\begin{tabular}{|c|c|c|c|c|c|c|}
\hline
\textbf{Step size} & \textbf{Instance} & $r$ & \textbf{iterations} & \textbf{same support} & \textbf{median $\Delta F$} & \textbf{true found} \\
 & & & (median / 90\%) & \textbf{as $r = 1$} & \textbf{vs $r = 1$} & (median, of 25) \\ \hline
\multirow{4}{*}{conservative} & \multirow{2}{*}{$p = 5{,}000$, $\rho = 0$} & 1 & 2 / 2 & 1.00 & 0 & 0 \\
 & & 10 & 41 / 76 & 0.42 & $-68\%$ & 13 \\ \cline{2-7}
 & \multirow{2}{*}{$p = 20{,}000$, $\rho = 0$} & 1 & 2 / 2 & 1.00 & 0 & 0 \\
 & & 10 & 11 / 11 & 0.94 & 0 & 0 \\ \hline
\multirow{8}{*}{$\eta_{\text{pred}}/4$} & \multirow{4}{*}{$p = 5{,}000$, $\rho = 0$} & 1 & 2 / 12 & 1.00 & 0 & 0 \\
 & & 2 & 13 / 26 & 0.41 & $-20\%$ & 13 \\
 & & 3 & 19 / 34 & 0.21 & $-66\%$ & 14 \\
 & & 10 & 34 / 76 & 0.03 & $-71\%$ & 16 \\ \cline{2-7}
 & \multirow{4}{*}{$p = 20{,}000$, $\rho = 0.5$} & 1 & 2 / 19 & 1.00 & 0 & 1 \\
 & & 2 & 12 / 31 & 0.51 & 0 & 17 \\
 & & 3 & 20 / 39 & 0.22 & $-23\%$ & 18 \\
 & & 10 & 34 / 67 & 0.01 & $-50\%$ & 18 \\ \hline
\multirow{4}{*}{$\eta_{\text{pred}}/2$} & \multirow{2}{*}{$p = 5{,}000$, $\rho = 0$} & 1 & 8 / 21 & 1.00 & 0 & 15 \\
 & & 10 & 22 / 35 & 0.48 & $-0.8\%$ & 16 \\ \cline{2-7}
 & \multirow{2}{*}{$p = 20{,}000$, $\rho = 0.5$} & 1 & 11 / 27 & 1.00 & 0 & 17 \\
 & & 10 & 22 / 36 & 0.81 & 0 & 17 \\ \hline
\end{tabular}
\caption{Effect of $r$, the number of consecutive updates that must leave the support unchanged. $\Delta F = (F_r - F_1)/F_1$ is the relative change in the training objective after value modification; negative values mean that $r$ found a better support than $r = 1$. Pooled over 3 seeds and all kept configurations. Selected rows; Figure~\ref{fig:stability_r} shows every instance and $r$.}
\label{tab:stability_r}
\end{table}

\begin{figure}[htbp]
    \centering
    \includegraphics[width=0.9\textwidth]{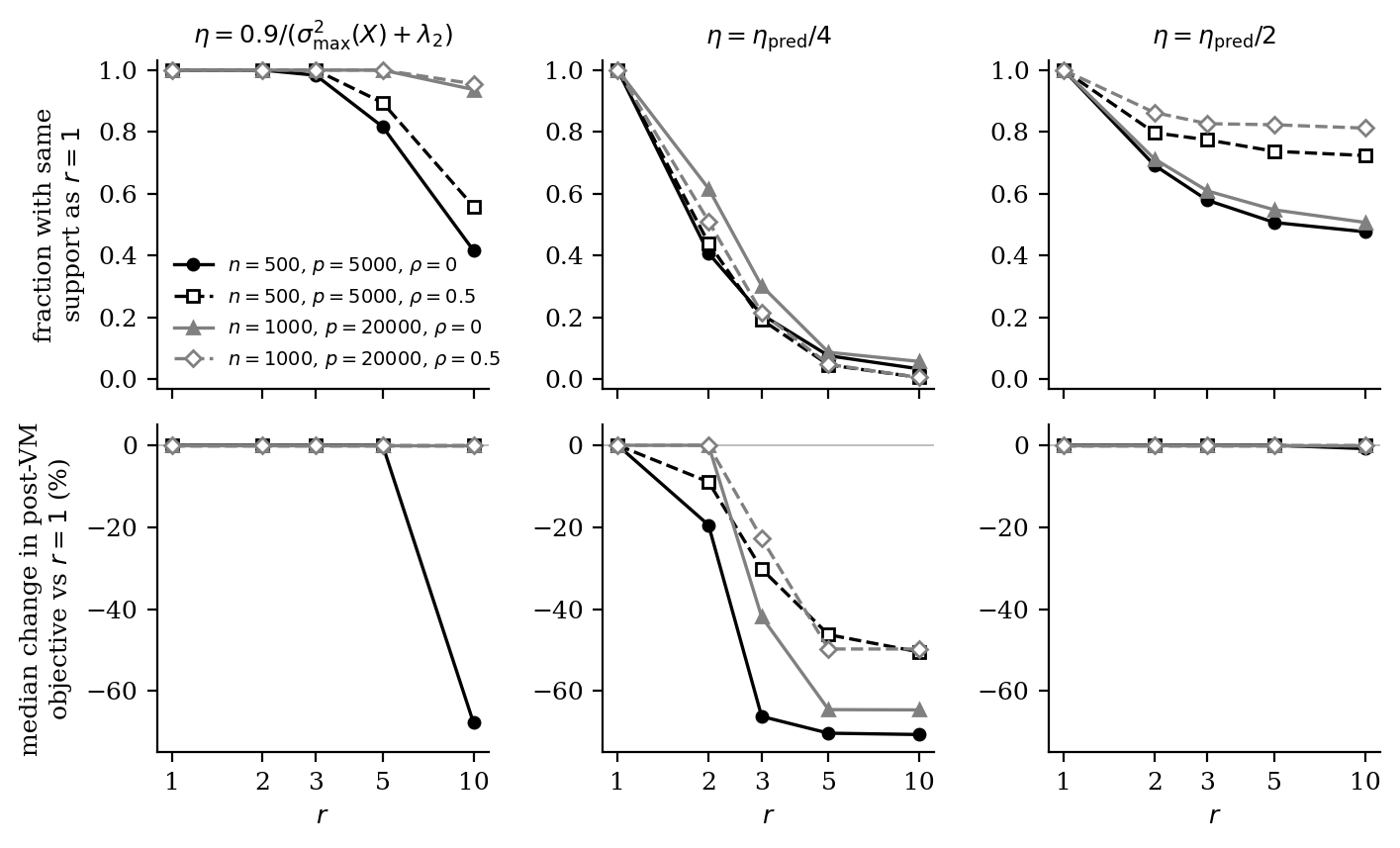}
    \caption{Top: fraction of configurations whose support under rule $r$ equals the support under $r = 1$. Bottom: median relative change in the training objective after value modification, relative to $r = 1$. Columns correspond to the three step sizes, with $\eta_{\text{pred}} = 0.9/(n(1+\sqrt{k/n})^2 + \lambda_2)$.}
    \label{fig:stability_r}
\end{figure}

\paragraph{Results.} Table~\ref{tab:stability_r} and Figure~\ref{fig:stability_r} summarize the results. The choice of $r$ matters, and how much it matters depends on the step size and the starting point.
\begin{itemize}
    \item \textbf{Nonzero starting points with small steps.} From a nonzero $\beta_o$, the first update keeps approximately the $k$ largest entries of $|\beta_o|$. When $\eta |\nabla F_D|$ is small relative to $|\beta_o|$, the second update leaves this support unchanged, and $r = 1$ stops after two updates on a support determined by the starting point rather than by the data. At $\eta_{\text{pred}}/4$, larger $r$ allows the iterates to move away from this support: the median objective after value modification is 50--71\% lower at $r = 10$ than at $r = 1$, and the median number of true features recovered rises from 0--6 to 15--18. From $\beta_o = 0$ the effect is smaller, with median improvements of 8--30\% at $r = 10$.
    \item \textbf{The conservative step.} With $\eta = 0.9/(\sigma_{\max}^2(X) + \lambda_2)$ and $p = 20{,}000$, even $r = 10$ keeps the $r = 1$ support in 94--96\% of configurations, and no $r \leq 10$ recovers true features from a nonzero start. Here the limiting factor is the step size, in line with Section~\ref{subsec:active_set_results}.
    \item \textbf{The largest stable step.} At $\eta_{\text{pred}}/2$, the median change is zero in every instance type, but individual configurations still improve by up to 9--36\%.
    \item \textbf{Direction and cost.} Increasing $r$ almost never worsened the training objective (6 of 8{,}640 comparisons, by at most 0.6\%); test error moved in the same direction in the median but not in every configuration. Supports were still changing between $r = 5$ and $r = 10$ in 26--53\% of configurations at $\eta_{\text{pred}}/4$, so these experiments do not identify a value of $r$ beyond which results stop changing. The cost of a larger $r$ is small: going from $r = 1$ to $r = 10$ adds roughly 10 to 50 updates per configuration, and support selection took between 1 and 121 ms per batch of 180 configurations for every $r$.
\end{itemize}
Configurations whose supports cycle are not stopped by any $r$ and run to $M_{SS}$; the check against the support two updates back used in Section~\ref{subsec:active_set_results} handles cycles of length two.

\section{Bitwise Differences in Dynamic Shrinking} \label{app:bitwise_differences}

Dynamic Shrinking does not give bitwise the same results as the \textit{Flag Approach}. Each resize changes the width of the batched matrix product, so the floating-point operations are summed in a different order, and the computed $\beta$s differ slightly. This appendix looks at how often that leads to a different support, and whether it matters.

In the first experiment of Section~\ref{subsec:resolve_cdpp} ($\prod C = 16{,}807$, $p = 10{,}000$), at most 35 configurations (0.208\%) end on a different support than under the \textit{Flag Approach}, depending on the shrink factor. In the second experiment ($\prod C = 3{,}125$, $p = 20{,}000$), the figure is 18 configurations (0.576\%).

A different support is not a defect in itself. When two coordinates have nearly the same magnitude at the hard-thresholding boundary, single precision cannot decide which of them is larger, and the problem does not distinguish between the two resulting supports, so either answer is equally valid. At $p = 10{,}000$, the evidence points to this explanation. On the configurations whose supports disagree, the two runs still share 95 to 98\% of their selected entries. The objective $F(\beta)$ in~\eqref{eq:our_estimator}, including both penalty terms, is also not systematically worse under \textit{Dynamic Shrinking}. Across $s \in \{0.5, 0.9, 0.99\}$, there are 60 disagreeing configurations. \textit{Dynamic Shrinking} reaches the lower objective on 33 of them and the higher objective on 27. The median relative difference is of order $10^{-4}$, and the worst case is 7.9\%. Compared with how much the objective varies across the configuration grid itself, the largest discrepancy is only $1.886 \times 10^{-3}$ of that range.

Two more observations support this reading. First, the number of disagreements does not grow with the number of gathers. Across the 9 shrink factors it stays between 4 and 35, while the number of gathers goes from 3 to 2436. If the disagreements came from rounding error building up with each resize, the count would rise with the number of gathers, and it does not. Second, the coefficients that differ are exactly the ones at the hard-thresholding boundary.

\section{Generating $y$ with a Given SNR} \label{app:generate_snr}

The experiment in Section~\ref{subsec:benefits} uses one synthetic instance $(X, y)$ with $n = 500$ observations, $p = 1{,}000$ features, a true coefficient vector $\beta^\star$ with $k^\star = 25$ nonzeros, and a signal-to-noise ratio (SNR) of 25. We define the SNR as
\begin{equation*}
\text{SNR} = \frac{\|X\beta^\star\|_2^2}{\|\epsilon\|_2^2},
\end{equation*}
where $y = X\beta^\star + \epsilon$. The same instance is used for all 3,125 configurations. Algorithm~\ref{alg:data_generation} gives the procedure.

\begin{algorithm}[h]
\caption{Synthetic Data Generation}
\label{alg:data_generation}
\begin{algorithmic}[1]
\REQUIRE number of observations $n$, number of features $p$, true sparsity $k^\star$, target signal-to-noise ratio $\text{SNR}$
\ENSURE $X \in \mathbb{R}^{n \times p}$, $y \in \mathbb{R}^n$, $\beta^\star \in \mathbb{R}^p$
\STATE Draw $Z \in \mathbb{R}^{n \times p}$ with entries $Z_{ij}$ independent standard Gaussian.
\FOR{$j = 1, \dots, p$}
\STATE $\bar{z}_j \gets \frac{1}{n}\sum_{i=1}^n Z_{ij}$ and $s_j \gets \big(\frac{1}{n}\sum_{i=1}^n (Z_{ij} - \bar{z}_j)^2\big)^{1/2}$
\STATE $X_{ij} \gets (Z_{ij} - \bar{z}_j)/s_j$ for $i = 1, \dots, n$ \COMMENT{column $j$ now has mean 0 and standard deviation 1}
\ENDFOR
\STATE Choose a set $S^\star$ of $k^\star$ indices uniformly at random from $\{1, \dots, p\}$, without replacement.
\STATE Set $\beta^\star_j$ to an independent standard Gaussian draw for $j \in S^\star$, and $\beta^\star_j \gets 0$ for $j \notin S^\star$.
\STATE $\mu \gets X\beta^\star$ \COMMENT{the noiseless signal}
\STATE Draw $u \in \mathbb{R}^n$ with entries independent standard Gaussian.
\STATE $\epsilon \gets \dfrac{\|\mu\|_2}{\sqrt{\text{SNR}}\,\|u\|_2}\, u$ \COMMENT{rescale the noise to the target ratio}
\STATE $y \gets \mu + \epsilon$
\RETURN $X$, $y$, $\beta^\star$
\end{algorithmic}
\end{algorithm}

For Section~\ref{subsec:benefits} we call Algorithm~\ref{alg:data_generation} with $n = 500$, $p = 1{,}000$, $k^\star = 25$ and $\text{SNR} = 25$. The step that forms $\epsilon$ rescales the noise by the norms of the drawn vectors $\mu$ and $u$, so the target is met exactly for this instance and not only in expectation:
\begin{equation*}
\frac{\|\mu\|_2^2}{\|\epsilon\|_2^2} = \frac{\|\mu\|_2^2}{\|\mu\|_2^2 \, \|u\|_2^2 / (\text{SNR}\, \|u\|_2^2)} = \text{SNR}.
\end{equation*}
The direction of $\epsilon$ is that of $u$, so the noise is still an isotropic Gaussian direction; only its length is fixed.

\section{Choosing the Active Set Size} \label{app:as_size}

Table~\ref{tab:as_main} reports one active set size $|\mathcal{A}|$ per $(p, k)$. This appendix gives the full sweep behind that choice, describes how much the results vary between runs, and lists the limitations of the procedure. The setup is that of Section~\ref{subsec:active_set_results}.

\begin{table}[htbp]
    \centering
    \small
    \setlength{\tabcolsep}{4pt}
    \begin{tabular}{|c|c|c|c|c|c|c|c|c|}
    \hline
        $p$ & $k$ & $|\mathcal{A}| = 2k$ & $3k$ & $4k$ & $5k$ & $8k$ & $10k$ & $15k$ \\ \hline
        \multirow{3}{*}{$10^5$} & 10 & 1.66 / $+0.2$ & 1.67 / $+0.2$ & 1.62 / $+0.2$ & 1.54 / $+0.2$ & 1.59 / $+0.2$ & 1.59 / $+0.2$ & 1.57 / $+0.5$ \\ \cline{2-9}
         & 50 & 1.81 / $-0.6$ & 1.92 / $-0.7$ & 1.66 / $-1.2$ & 2.09 / $-1.0$ & 1.65 / $-0.9$ & 1.51 / $-0.5$ & 1.52 / $-0.6$ \\ \cline{2-9}
         & 100 & 0.39 / $-69.5$ & 0.33 / $-70.1$ & 0.30 / $-71.5$ & 0.29 / $-71.8$ & 1.03 / $-2.7$ & 1.00 / $+1.5$ & 0.28 / $-6.5$ \\ \hline
        \multirow{5}{*}{$3 \times 10^5$} & 10 & 2.05 / $+0.4$ & 2.02 / $+0.3$ & 2.02 / $+0.1$ & 2.01 / $+0.0$ & 2.01 / $+0.1$ & 2.00 / $+0.1$ & 1.97 / $+0.1$ \\ \cline{2-9}
         & 50 & 2.71 / $+0.1$ & 2.36 / $+0.1$ & 2.70 / $+0.0$ & 3.03 / $-0.3$ & 3.06 / $+0.3$ & 2.94 / $+0.0$ & 2.76 / $+0.1$ \\ \cline{2-9}
         & 100 & 0.56 / $-58.6$ & 0.51 / $-59.2$ & 0.48 / $-58.7$ & 0.37 / $-14.7$ & 0.25 / $-28.3$ & 0.29 / $-27.1$ & 0.18 / $-42.0$ \\ \cline{2-9}
         & 250 & 2.13 / $+10.6$ & 1.81 / $+7.5$ & 1.20 / $+12.7$ & 1.95 / $+9.5$ & 2.02 / $+25.5$ & 0.80 / $+15.4$ & 1.12 / $+5.3$ \\ \cline{2-9}
         & 500 & 3.83 / $-16.5$ & 4.73 / $-17.1$ & 4.35 / $-32.1$ & 4.85 / $-30.7$ & 4.40 / $-28.3$ & 0.59 / $-7.2$ & 0.43 / $+61.5$ \\ \hline
        \multirow{5}{*}{$10^6$} & 10 & 2.14 / $-0.0$ & 2.16 / $-0.1$ & 2.14 / $-0.1$ & 2.14 / $-0.1$ & 2.14 / $+0.0$ & 2.15 / $+0.0$ & 2.14 / $-0.1$ \\ \cline{2-9}
         & 50 & 2.66 / $-0.5$ & 2.62 / $-0.2$ & 2.65 / $-0.6$ & 2.62 / $-0.4$ & 2.60 / $-0.4$ & 2.59 / $-0.2$ & 2.57 / $-0.2$ \\ \cline{2-9}
         & 100 & 6.98 / $-34.9$ & 6.76 / $-35.1$ & 6.63 / $-35.4$ & 6.49 / $-35.5$ & 6.24 / $-35.3$ & 10.86 / $+15.8$ & 10.39 / $+7.8$ \\ \cline{2-9}
         & 250 & 2.50 / $+33.7$ & 2.68 / $+1.9$ & 2.66 / $+1.7$ & 2.31 / $+6.0$ & 2.98 / $+19.6$ & 4.44 / $+21.7$ & 2.09 / $+23.8$ \\ \cline{2-9}
         & 500 & 15.72 / $-1.7$ & 21.55 / $+9.9$ & 21.81 / $+18.7$ & 5.89 / $+24.8$ & 7.39 / $+6.2$ & 16.12 / $+155.8$ & 19.06 / $+157.4$ \\ \hline
    \end{tabular}
    \caption{The full $|\mathcal{A}|$ sweep. Each cell gives the solve-only speedup of the active set over full-width support selection (learning rate selection excluded from both arms), followed by the objective difference in percent, as defined in Table~\ref{tab:as_main}. Negative differences mean the active set found the better support.}
    \label{tab:as_sweep}
\end{table}

\paragraph{Speed.} Along most rows the solve-only speedup changes little with $|\mathcal{A}|$. Even $15k$ features is small next to $p$, so the run time is set mainly by the number of full-width updates, not by the size of $X_\mathcal{A}$. The exceptions are rows where some sizes need many more rounds than others, such as $p = 3 \times 10^5$, $k = 500$, where $10k$ and $15k$ need 17 and 25 rounds against 2 for the smaller sizes.

\paragraph{Accuracy.} The objective depends on $|\mathcal{A}|$ much more than the run time does. For $k \leq 50$ every size gives an objective within about 1\% of full width. For $k \geq 100$ the differences are large, and the largest active sets are the most likely to be worse. With a large $\mathcal{A}$ the inner loop can settle on a support that the full-width check accepts after only one or two rounds, so the outer loop stops before the support has been improved. Smaller active sets force more rounds. The pattern is not monotone, and it varies between seeds, so we do not identify a single best size. In every row except $p = 3 \times 10^5$, $k = 250$, at least one size between $2k$ and $5k$ is within 2\% of full width or better, and these are the sizes used in Table~\ref{tab:as_main}. No single size is safe everywhere, however: $2k$ is 34\% worse at $p = 10^6$, $k = 250$, and $3k$ to $5k$ are 10 to 25\% worse at $p = 10^6$, $k = 500$. At $p = 3 \times 10^5$, $k = 250$, every size we tried gives a worse objective than full width.

\paragraph{Variation between runs.} The objective differences above are averages over five starting points, and the medians over three seeds. Individual configurations vary more. For example, at $p = 10^6$, $k = 500$ and $|\mathcal{A}| = 2k$, the average is 1.7\% in favor of the active set, but one of the 15 configurations is 22\% worse. Much of this variation does not come from the active set. Full-width support selection is itself sensitive to small changes, since the problem is nonconvex and hard thresholding is discontinuous (Appendix~\ref{app:as_supports}). To measure this, we reran full-width support selection with its learning rate multiplied by 0.9 and nothing else changed, over five seeds. Measured as in Table~\ref{tab:as_main}, this alone changes the objective by $-62\%$ to $+96\%$ for $k \geq 100$, and by up to 33\% at $k = 10$. A difference between the two arms of a few percent in either direction is therefore within the variation of full-width support selection itself.

\paragraph{Limitations.} First, the cycle check only detects cycles of length two, which are the only ones we observed. A longer cycle would run until the cap $M_{AS}$; no run here reached it, but this need not hold for other data. Second, the speedup depends on how many full-width updates support selection would need without the active set. This is not known in advance and ranges here from 3 to 80 for the same design. The cost of the active set itself is much more predictable.

\section{Why the Active Set Returns a Different Support} \label{app:as_supports}

The active set is intended to carry out the same procedure as full-width support selection at a lower cost. Its inner loop only rearranges features within $\mathcal{A}$, and the full-width check verifies that no feature outside $\mathcal{A}$ should enter, so the returned support satisfies the same stopping condition as that of the arm without an active set. Nevertheless, the two arms often disagree. Across all 1{,}365 configurations of Section~\ref{subsec:active_set_results}, the median fraction of the full-width support that the active set also selects is 0.60, and the supports agree completely in only 6.5\% of configurations. The overlap depends strongly on $k$. It is 0.90 to 0.96 at $k = 10$, and between 0.17 and 0.30 for $k \geq 100$. At $p = 10^6$ and $k = 100$, the four nonzero starting points give an overlap of 0.00 to 0.01, so the two supports of 100 features share at most one feature. The start $\beta_o = 0$ always gives the largest overlap of the five, as expected, since from zero both arms begin by ranking the features by $|X^Ty|$.

The disagreement is a consequence of the nonconvexity of Problem~\eqref{eq:our_estimator}. Because of the constraint $\|\beta\|_0 \leq k$, the feasible set is a union of $\binom{p}{k}$ subspaces, and the problem has many local minima. Our method, like any tractable method for this problem, returns one of them. The stopping rule is a stationarity condition that many supports satisfy, so two procedures with the same stopping rule can end at different supports, depending on the path each one takes. The paths of the two arms differ by construction, since the inner loop of the active set moves only within $\mathcal{A}$ and uses its own step size. Moreover, hard thresholding is discontinuous, so a small difference between two iterates can lead to different supports after a single update, after which the two runs evolve independently.

Rounding error does not explain these differences. Section~\ref{subsec:resolve_cmtp} shows its effect in this code: three methods that perform the same arithmetic in different orders return identical supports, with coefficients that agree to $10^{-6}$. Rounding can decide between two nearly tied features, but it cannot explain why one arm finds the signal while the other misses it. Table~\ref{tab:as_support} shows what happens for one instance. For each arm, it reports how many of the 50 true features are recovered and whether the returned support is a fixed point, that is, whether one more full-width update leaves it unchanged.

\begin{table}[htbp]
    \centering
    \small
    \begin{tabular}{|c|c|c|c|c|c|c|}
    \hline
        \multirow{2}{*}{$k$} & \multirow{2}{*}{$\beta_o$} & \multicolumn{2}{c|}{\textbf{No active set}} & \multicolumn{2}{c|}{\textbf{Active set}} & \multirow{2}{*}{\textbf{overlap}} \\ \cline{3-6}
         & & \textbf{true found} & \textbf{fixed pt.} & \textbf{true found} & \textbf{fixed pt.} & \\ \hline
        \multirow{5}{*}{100} & zeros & 19/50 & yes & 21/50 & yes & 0.85 \\ \cline{2-7}
         & ones & \cellcolor{red!15}0/50 & no & 22/50 & yes & 0.01 \\ \cline{2-7}
         & neg. ones & 29/50 & yes & 23/50 & yes & 0.44 \\ \cline{2-7}
         & Gaussian & \cellcolor{red!15}0/50 & yes & 22/50 & yes & 0.01 \\ \cline{2-7}
         & uniform & \cellcolor{red!15}0/50 & no & 23/50 & yes & 0.01 \\ \hline
        \multirow{5}{*}{500} & zeros & 23/50 & yes & 22/50 & yes & 0.40 \\ \cline{2-7}
         & ones & 27/50 & no & 17/50 & yes & 0.13 \\ \cline{2-7}
         & neg. ones & 28/50 & no & 27/50 & yes & 0.18 \\ \cline{2-7}
         & Gaussian & 27/50 & no & 15/50 & no & 0.16 \\ \cline{2-7}
         & uniform & 24/50 & no & 27/50 & no & 0.21 \\ \hline
    \end{tabular}
    \caption{What each arm returns at $p = 10^6$, $|\mathcal{A}| = 5k$, for one instance (seed 0). ``true found'' counts the features of the true signal in the returned support, ``fixed pt.'' records whether one more full-width update leaves the support unchanged, and ``overlap'' is the fraction of the full-width support that the active set also selects.}
    \label{tab:as_support}
\end{table}

The low overlap has two different causes, which appear at different values of $k$. At $k = 100$, from three of the five starting points, full-width support selection recovers \emph{none} of the 50 true features, whereas the active set started from the same point recovers 22 or 23. The overlap of 0.01 in these rows therefore does not reflect a choice between interchangeable features. This accounts for the objective difference of about $-35\%$ at $p = 10^6$, $k = 100$, which holds for every $|\mathcal{A}|$ from $2k$ to $8k$ (Tables~\ref{tab:as_main} and~\ref{tab:as_sweep}). Here the cheaper procedure is also the more reliable one. These failures are local optima rather than early stops. Every full-width run in the table stops through the ordinary rule and none through the cycle check, and the Gaussian row is a fixed point that contains no true feature, so additional updates would not help. With a dense starting point, $\beta_o$ dominates the first update, so the top $k$ is drawn from $\beta_o$ rather than from the data, and the full-width step size is too small to move away from this support. The active set avoids this because its inner loop works on $X_\mathcal{A}$, where it can take much larger steps and rebuild the support before the full-width check.

At $k = 500$ the cause is near-ties. Both arms recover part of the true signal (23 to 28 true features without the active set and 15 to 27 with it), but different parts. There are 113 features within 1\% of the 500th largest $|X^Ty|$, compared with 28 at $k = 100$ and only 2 at $k = 10$. With so many nearly equal candidates, the selection among them is close to arbitrary. This explains the low overlap at $k = 500$, why neither arm is consistently better there, and why the two arms agree almost perfectly at $k = 10$.

The fixed point column shows a further difference. The active set returns a fixed point in five of five configurations at $k = 100$ and three of five at $k = 500$, compared with three and one for the arm without it. The stopping rule compares supports only after a single update while the coefficients continue to change, so a support that is unchanged after one update may still change later. The active set satisfies the rule in this stronger sense more often, because it solves to stability within $\mathcal{A}$ before each full-width check.

In summary, both causes stem from the nonconvexity of the problem. At $k = 100$ one arm reaches a poor local optimum that the other avoids, and at $k = 500$ there are many local optima of nearly equal quality. In neither case is the difference due to rounding error, and in neither case does the active set depart from the method it implements.

\section{Amortizing Learning Rate Selection} \label{app:as_amortize}

Learning rate selection (Section~\ref{subsec:paradeal_finds_eta}) depends only on the design and on $(k, \lambda_2)$. One selection therefore serves the whole batch, however many values of $\beta_o$ and $\lambda_1$ it contains, and its share of the run time should fall as $\prod C$ grows. To test this, Table~\ref{tab:as_amort} grows $\prod C$ along the $\lambda_1$ axis at $p = 10^6$ and $k = 100$, with $|\mathcal{A}| = 5k$ and one seed.

\begin{table}[htbp]
    \centering
    \begin{tabular}{|r|c|c|c|c|c|c|}
    \hline
        \multirow{2}{*}{$\prod C$} & \multicolumn{2}{c|}{\textbf{No active set}} & \multicolumn{2}{c|}{\textbf{Active set}} & \multicolumn{2}{c|}{\textbf{Speedup}} \\ \cline{2-7}
         & \textbf{total (s)} & \textbf{LR share} & \textbf{total (s)} & \textbf{LR share} & \textbf{total} & \textbf{solve} \\ \hline
        5 & 0.55 & 30.3\% & 0.18 & 75.2\% & 3.03$\times$ & 8.52$\times$ \\ \hline
        25 & 0.74 & 22.6\% & 0.20 & 68.8\% & 3.68$\times$ & 9.11$\times$ \\ \hline
        125 & 1.94 & 8.5\% & 0.36 & 40.5\% & 5.45$\times$ & 8.39$\times$ \\ \hline
        625 & 7.79 & 2.1\% & 1.09 & 14.2\% & 7.15$\times$ & 8.16$\times$ \\ \hline
    \end{tabular}
    \caption{Growing $\prod C$ along the $\lambda_1$ axis at $p = 10^6$, $k = 100$, $|\mathcal{A}| = 5k$. ``LR share'' is the fraction of the total run time spent on learning rate selection, which takes 0.17 s in every row.}
    \label{tab:as_amort}
\end{table}

Learning rate selection takes 0.17 s in every row, while the solve time grows with $\prod C$. Its share of the run without an active set falls from 30\% to 2\%, and the end-to-end speedup rises from 3.0 to 7.2, approaching the solve-only speedup of about 8. At $\prod C = 5$, learning rate selection takes up most of the active set's run time, so the end-to-end figure understates the benefit of the active set by nearly a factor of three. Since ParaDeAl is designed for large batches, the solve-only speedup is the more relevant figure in practice.

\section{Parallelizing the Active Set Across $k$} \label{app:as_across_k}

In Section~\ref{subsec:active_set_results} each batch has a single value of $k$, and different values of $k$ are run one after another. The reason is that $k$ is the only hyperparameter that makes the active set ragged. If $|\mathcal{A}|$ depends on $k$, as with $|\mathcal{A}| = 5k$, then a batch with several values of $k$ produces blocks $X_\mathcal{A}$ of different widths. These blocks must be padded to the largest width and handled with an overlay, as in \textit{Resizing to Candidate Set} (Section~\ref{subsec:resolve_cmtp}). A configuration with $k = 10$ in a batch whose largest $k$ is 100 then does ten times more arithmetic than it needs. Unlike in value modification, this padded gather is repeated in every round.

The problem disappears if $|\mathcal{A}|$ is the same for every configuration in the batch. Every block then has size $n \times |\mathcal{A}|$, so no padding is needed, and the only remaining difference between configurations is the $k$ used in the hard threshold, which Algorithm~\ref{alg:topk_parallel} already handles. A configuration with a small $k$ then has a larger active set than it needs. This does no harm, since $\mathcal{A}$ only needs to contain the support that is eventually selected.

We test this at $p = 10^6$ with $k \in \{10, 50, 100\}$ and a shared $|\mathcal{A}| = 500$. We compare one batch of 15 configurations (five starting points for each $k$) with three batches of 5 run one after another. The single batch is 1.82 times faster on the solve (0.088 s against 0.159 s), because the three values of $k$ share their full-width updates, which reduces their number from 11 to 5. The end-to-end gain is only 1.14 times, because at this $\prod C$ both runs are dominated by learning rate selection, which is still run once per value of $k$. The two approaches return the same supports, and their objectives after value modification agree to within $10^{-7}$ in relative terms on every configuration and seed. The cost is memory: 0.567 GiB against 0.191 GiB, since all 15 configurations are held at once.

Batching across $k$ is therefore available when it is needed, at the cost of holding $|\mathcal{A}|$ fixed and of the memory to keep every configuration in the batch at once. For the experiments of Section~\ref{subsec:active_set_results} we keep $k$ sequential. This costs little, since a single case at $p = 10^6$ takes well under a second.

\section{The Competitors in Detail} \label{app:competitors}

This appendix describes how each competitor of Section~\ref{sec:compare_results} works and how we run it. All CPU methods run on a single core of the same node, and each method, dataset and grid runs in its own process.

\subsection*{L0Learn}

L0Learn \citep{haz20,haz23} fits $L_0$-regularized linear models. For squared error loss it minimizes
\[
\tfrac12\|y - X\beta\|_2^2 + \lambda_0\|\beta\|_0 + \gamma_1\|\beta\|_1 + \gamma_2\|\beta\|_2^2,
\]
with the $L_0$ penalty on its own or paired with one of the other two penalties (the \texttt{L0L1} and \texttt{L0L2} penalties). The main difference from ParaDeAl is that sparsity is penalized, not constrained: the user cannot ask for a support of size $k$, only for a value of $\lambda_0$, and L0Learn reports whatever support that value produces.

The core algorithm of L0Learn is cyclic coordinate descent. Each coordinate is set to the exact minimizer of the objective with all other coordinates fixed, which for the $L_0$ penalty is a thresholding rule. \citet{haz20} show that the points to which coordinate descent converges, the coordinate-wise minima, form a subset of the fixed points of iterative hard thresholding. Coordinate descent can therefore stop at better solutions than IHT, but never at worse ones, which is part of why L0Learn is a strong accuracy benchmark. An optional local search that swaps features in and out of the support (the \texttt{CDPSI} algorithm) improves the solutions further; we use the default, plain coordinate descent (\texttt{CD}).

L0Learn computes a whole path of solutions per call. For each value of $\gamma$ (a grid of \texttt{nGamma} values between \texttt{gammaMin} and \texttt{gammaMax}), L0Learn chooses a decreasing sequence of \texttt{nLambda} values of $\lambda_0$, each set so that the support changes, shrinking by the factor \texttt{scaleDownFactor} $= 0.8$. Every solution on the path is warm-started from the previous one. Two further devices keep each fit cheap: an active set (coordinate descent cycles only over the current support until the support stabilizes, then checks the remaining features) and correlation screening (at most \texttt{screenSize} $= 1000$ candidate features are considered at a time). The path stops early once the support exceeds \texttt{maxSuppSize}. The solver is written in C++ with the Armadillo library and called from R.

\paragraph{Settings.} Version 2.1.0 in R 4.4.0, one core, \texttt{algorithm = "CD"}, \texttt{nGamma} $=$ \texttt{nLambda} $= G$ for $G \in \{10, 30, 100\}$, \texttt{maxSuppSize} $= 500$, \texttt{intercept = FALSE}, and the remaining settings at their defaults (\texttt{maxIters} $= 200$, \texttt{rtol} $= 10^{-6}$, active set and screening on).

\subsection*{glmnet}

glmnet \citep{fri10} is the standard solver for the lasso and the elastic net. For the Gaussian family it minimizes
\[
\tfrac{1}{2n}\|y - X\beta\|_2^2 + \lambda\Big(\alpha\|\beta\|_1 + \tfrac{1-\alpha}{2}\|\beta\|_2^2\Big),
\]
where $\alpha \in [0,1]$ sets the mix of the two penalties and $\lambda$ their overall strength. glmnet has no sparsity constraint and no $L_0$ penalty, so the comparison with glmnet is on the elastic net alone.

glmnet is built around the regularization path. For a fixed $\alpha$, the path starts at the smallest $\lambda$ for which $\beta = 0$, decreases $\lambda$ geometrically over \texttt{nlambda} values down to \texttt{lambda.min.ratio} times that value (0.01 when $n < p$, as here), and warm-starts each fit from the previous solution. Neighboring solutions on the path are close, so each warm-started fit needs few passes. Within a fit, glmnet runs cyclic coordinate descent with two devices that make it fast on wide data. Sequential strong rules \citep{tib12} discard most features before the fit starts, based on the previous solution, and check the discarded features afterwards. Active-set cycling iterates over the current nonzeros until they converge before sweeping all features. With $p \ge 500$, glmnet uses the ``naive'' updates, which keep the residual up to date instead of forming inner products of features. The path can also end early, once the explained deviance stops improving (by less than $10^{-5}$) or reaches 0.999. The solver is written in Fortran and called from R.

This design is the main obstacle for ParaDeAl. A warm-started path with screening does very little work per additional $\lambda$, while ParaDeAl solves every configuration from zero at full width. ParaDeAl has to win through parallelism and through stopping each configuration as soon as that configuration has converged.

\paragraph{Settings.} Version 5.0 in R 4.4.0, one core, a \texttt{for} loop over $G_\alpha$ values of $\alpha$ evenly spaced in $[0.05, 0.95]$, each call computing a path of \texttt{nlambda} $= G_\lambda$ values, with \texttt{standardize = FALSE}, \texttt{intercept = FALSE} and the default convergence threshold \texttt{thresh} $= 10^{-7}$.

\subsection*{scikit-learn}

The \texttt{ElasticNet} class of scikit-learn \citep{ped18} minimizes
\[
\tfrac{1}{2n}\|y - Xw\|_2^2 + a\,r\,\|w\|_1 + \tfrac{a(1-r)}{2}\|w\|_2^2 ,
\]
the same objective as glmnet in a different parametrization. The solver is cyclic coordinate descent written in Cython and stops when the largest coordinate update falls below the tolerance and the duality gap confirms convergence. scikit-learn is the most widely used Python implementation and serves here as the plain reference: one configuration at a time, no path, no screening. (scikit-learn does offer a path function; the comparison fits each configuration separately, as a user sweeping a grid would.)

\paragraph{Settings.} Version 1.2.2, one core, \texttt{fit\_intercept = False}, \texttt{precompute = False}, \texttt{selection = "cyclic"}, \texttt{tol} $= 10^{-4}$ (the default) and \texttt{max\_iter} $= 10{,}000$ so that every fit converges (none reported a convergence warning). Fitting a $100\times100$ grid this way would take days on the larger datasets, so we fit a fixed random subset of at most 50 of glmnet's configurations and report the per-configuration average times the size of the grid.

\subsection*{skglm}

skglm \citep{ber22} is a recent Python library designed for speed on sparse generalized linear models. The library combines a working-set strategy (a sequence of small subproblems on the features most likely to be nonzero, with the set growing until the optimality conditions hold for all features) with coordinate descent accelerated by Anderson extrapolation, and compiles its inner loops with Numba. The \texttt{ElasticNet} class of skglm uses the parametrization and interface of scikit-learn, so skglm can replace scikit-learn directly.

\paragraph{Settings.} Version 0.3.1 (the latest release compatible with scikit-learn 1.2.2 on the cluster), one core, \texttt{fit\_intercept} $=$ \texttt{False}, \texttt{tol} $= 10^{-4}$, on the same subset of configurations as scikit-learn, with one untimed warm-up fit so that Numba compilation is not counted.

\subsection*{Converting the penalty parameters}

Every solution that a competitor returns defines one ParaDeAl configuration. Both conversions below were verified on the returned solutions, through stationarity on the support (median relative violation $10^{-9}$ for L0Learn and $10^{-3}$ for glmnet, the latter at glmnet's own tolerance).
\begin{itemize}
\item L0Learn rescales $y$ to unit norm before solving, so a returned $\gamma$ corresponds to $\lambda_1 = 2\gamma\|y\|_2$ (L0+L1) or $\lambda_2 = 2\gamma$ (L0+L2) in ParaDeAl's objective, and the support size of the solution becomes $k$.
\item glmnet fits $y/s_y$ with $s_y = \sqrt{\mathrm{mean}(y^2)}$ and reports $\lambda$ in the units of $y$, which is exact for the $L_1$ term but not for the ridge term. A returned pair corresponds to $\lambda_1 = 2n\lambda\alpha$ and $\lambda_2 = n\lambda(1-\alpha)/s_y$. Without the factor $1/s_y$, ParaDeAl and glmnet would solve different problems.
\item scikit-learn and skglm: $a\,r = \lambda_1/2n$ and $a(1-r) = \lambda_2/n$.
\end{itemize}

\subsection*{What each comparison measures}

The four competitors test different things. Against L0Learn, the question is how much speed ParaDeAl gains by parallelizing over the whole grid, and how much accuracy ParaDeAl gives up by using IHT with a sparsity constraint instead of coordinate descent with a penalty. Against glmnet, the question is whether solving every configuration at once, in parallel, can beat a sequential path that reuses each solution to start the next. Against scikit-learn and skglm, the question concerns the more common situation of a user who fits one configuration at a time.

\section{Construction of the Real Datasets} \label{app:data_construction}

\paragraph{Why the feature count is expanded.} In practice, the number of samples $n$ is far easier to scale than the number of features $p$: collecting additional data (people, measurements) is a matter of equipment and funding, while devising new characteristics to measure is much harder. Public regression datasets therefore have few features (Section~\ref{subsec:datasets}). Feature expansion should preserve the information in the original features. We considered two techniques and use both.
\begin{itemize}
    \item \textbf{Powers and interactions.} The $w$th power includes all terms of the form $\prod_{i=1}^{p} x_i^{r_i}$, with $x_i$ the $i$th feature and $\sum_{i=1}^p r_i = w$. The total number of features up to power $w$ is $\sum_{i=1}^w {i+p-1 \choose i} = {w+p \choose w} - 1 \approx p^w/ w!$; for $p = 3$ and $w = 50$ there are already $23{,}425$ features. This construction keeps the original features and lets a linear model capture nonlinear relationships and interactions. The experiments use $w = 2$ without squares: the original columns and every product $x_ix_j$, $i<j$.
    \item \textbf{Nearest-neighbor features.} When $n$ is large, each sample can be augmented with the features of its $h$ nearest neighbors, giving $(h+1)p$ features; for $p = 25$ and $n = 1{,}000$, $h = 50$ gives $1{,}275$ features. The neighbors may be scaled or transformed to dilute their influence on the original row. The same construction applies to temporal data, by flattening the previous $h$ rows into one row vector.
\end{itemize}
Categorical variables are encoded numerically and missing values are imputed or removed.

\paragraph{The three datasets.} Neighbors are found in the z-scored base features and exclude the row itself.
\begin{itemize}
\item \textbf{Diabetes} \citep{uci_diabetes}: 47 features (every column except the two IDs and the response, time in hospital), pairwise interactions (1{,}128 columns), then the interaction features of the 250 nearest neighbors found among all 101{,}766 rows ($p = 283{,}128$); 2{,}500 training and 1{,}000 test rows.
\item \textbf{Communities and Crime} \citep{uci_crime}: 122 features with mean imputation, response \texttt{Violent\-Crimes\-Per\-Pop}, pairwise interactions (7{,}503 columns), then 10 nearest neighbors ($p = 82{,}533$); 1{,}500 training and 250 test rows.
\item \textbf{Air Quality} \citep{uci_airquality}: 12 features predicting the temperature \texttt{T}, 15 nearest neighbors (192 columns), then pairwise interactions ($p = 18{,}528$); 1{,}000 training and 1{,}000 test rows.
\end{itemize}
Every column of the training $X$ has mean 0 and $L_2$ norm 1, the test set receives the training transformation, and $y$ is centered with the training mean.

\paragraph{The three polynomial datasets.} These datasets use no neighbor features. Each base feature is replaced by its normal score: its values are ranked (over all $n$ rows, without using $y$), rank $r$ is mapped to $\Phi^{-1}((r - 1/2)/n)$, where $\Phi$ is the standard normal distribution function, and the result is standardized. This keeps the order of the values but removes skewness and outliers, so that products of features are not dominated by a few extreme rows. The features are then the base features, the products of two distinct base features and the products of three distinct base features (no squares or cubes), each scaled as above.
\begin{itemize}
\item \textbf{Communities and Crime}: the same 122 base features, response and rows as above; all $122 + 7{,}381 + 295{,}240 = 302{,}743$ products.
\item \textbf{NHANES} (cycles 2015--16 and 2017--18): adults aged 20 or over with a systolic blood pressure reading. The response is the mean of up to three systolic readings, and the 52 base features are examination and laboratory values. A value is left out if more than 30\% of its entries are missing, below the detection limit, or equal to its most common value, or if it has correlation above 0.95 with a value already kept; rows with missing entries are removed, and values with skewness above 2 are log-transformed. All $23{,}478$ products; 2{,}500 training and 1{,}000 test rows.
\item \textbf{Weather}: the weather prediction dataset of \citet{huber_weather}, daily measurements from 18 European cities in 2000--2010. The response is the next day's mean temperature in Basel, as a deviation from its seasonal average; the 128 base features are the same-day measurements of the other cities, each also expressed as a deviation from its seasonal average, so that the shared yearly cycle does not make every column look alike. Measurements with correlation above 0.95 with one already kept are left out. Of the $349{,}632$ possible products, $100{,}000$ are sampled at random. The training rows are the first 2{,}500 days and the test rows are the remaining 1{,}123 days after a gap of 30 days, so the test set lies entirely in the future.
\end{itemize}

\paragraph{Clipping the Diabetes test features.} Many Diabetes features are rare categories (medications almost nobody takes, uncommon admission codes). Repeating these features through 251 neighbor blocks produces columns that are constant (35{,}975 of them, set to 0) or nearly constant on the 2{,}500 training rows. Scaling a nearly constant column to unit norm multiplies the column by a very large factor, so a test row that differs from the training rows in such a column takes a huge value: 332 Diabetes test columns reach magnitudes up to 4{,}700, against at most 1 in training. Any coefficient on such a column then produces a large test error, for every method. On Crime and Air Quality, the test values stay within a factor of 4 of the training range. We therefore also report the test RE with each test feature clipped to the training range of its column, which is what one would do at prediction time; the solutions are the same, and only the test features used for scoring differ.

\section{Experimental Setup in Detail} \label{app:setup_details}

\paragraph{Grids.} L0Learn: \texttt{nGamma} $=$ \texttt{nLambda} $= 10, 30, 100$ and \texttt{maxSuppSize} $= 500$; each solution gives one ParaDeAl configuration, with $k$ equal to L0Learn's support size (support size 0 and repeated tuples are dropped). glmnet: 10, 30 or 100 values of $\alpha$ evenly spaced in $[0.05, 0.95]$, each with a path of 10, 30 or 100 values of $\lambda$, plus the two additional grids $10\alpha\times100\lambda$ and $30\alpha\times100\lambda$. L0+L1+L2 uses glmnet's $(\lambda_1, \lambda_2)$ pairs with $k = 50$.

\paragraph{ParaDeAl.} fp32 throughout (TF32 off), with all configurations of a run solved together on one L40S, in chunks sized to a 30\,GB working budget. In the L0 modes, support selection is proximal gradient with hard thresholding, with at most 100 iterations per sweep and the support-unchanged stopping rule, inside the active set procedure. The step size of every configuration comes from the parallel search for the largest stable $\eta$: a geometric ladder of 17 multiples of $\eta_0 = 1/(2(\sigma_{\max}^2(X) + \lambda_2))$, 25 probe iterations, and one probe per distinct pair $(k, \lambda_2)$. Support selection returns the iterate with the lowest objective seen in any sweep or active-set round (Section~\ref{subsec:lowest_objective}). Value modification on the fixed support uses conjugate gradients when $\lambda_1 = 0$ and proximal gradient otherwise. Pair swaps (Section~\ref{subsec:pair_swap}) then run with $M = 16$ candidate additions, $D = 8$ candidate drops, a relative improvement threshold of $10^{-7}$ and at most 200 rounds; each accepted move is refit by value modification at the same tolerance, warm-started from the current coefficients.

\paragraph{Final configuration.} The results of Section~\ref{subsec:real_results} come from a second full run of ParaDeAl; the competitor runs are unchanged, since their results do not depend on ParaDeAl. Relative to the first run, the second run changes three things.
\begin{itemize}
\item \textbf{Dynamic Shrinking everywhere.} Every loop in which configurations iterate (IHT sweeps, active-set rounds, conjugate gradients, value modification and L1+L2) freezes a configuration once that configuration has stopped. Whenever another 5\% of the current batch has stopped, the loop writes the finished configurations out and compacts every per-configuration tensor to the configurations still running. Compaction stops once 50 or fewer configurations are left, where its indexing costs more than it saves. The first run used the \textit{Flag Approach}, which keeps computing finished configurations until the whole batch is done.
\item \textbf{A looser stopping tolerance.} Every value solver (L1+L2, value modification, conjugate gradients) stops a configuration when its relative change $\|\beta_{t+1}-\beta_t\|/\|\beta_{t+1}\|$ (for conjugate gradients, its relative residual) falls below the tolerance. The first run used $10^{-6}$, which fp32 cannot resolve, so nearly every configuration ran to the iteration cap. L1+L2 is now run at $10^{-3}$ and at $5\cdot10^{-4}$, and the L0 modes at $5\cdot10^{-4}$. Support selection keeps its support-unchanged rule.
\item \textbf{The active set in every L0 run}, including Air Quality, with $|A| = \max(250, 5k)$ instead of $\max(250, 2k)$.
\end{itemize}

\paragraph{Three solvers for L1+L2.} L1+L2 has no thresholding and is run in three variants: plain proximal gradient with at most 1{,}000 and with at most 5{,}000 iterations, and FISTA with at most 5{,}000 iterations. Section~\ref{subsec:real_results} reports FISTA.

\paragraph{FISTA.} Plain proximal gradient (ISTA) takes each step from the current iterate:
\[
\beta_{t+1} = \operatorname{soft}\big(\beta_t - \eta\nabla F(\beta_t),\ \eta\lambda_1\big),
\qquad \nabla F(\beta) = 2X^T(X\beta - y) + 2\lambda_2\beta .
\]
FISTA \citep{beck2009fast} takes the identical step, but from a point pushed further along the direction in which the iterates were already moving:
\[
z_t = \beta_t + \frac{s_{t-1} - 1}{s_t}\,(\beta_t - \beta_{t-1}),\qquad
\beta_{t+1} = \operatorname{soft}\big(z_t - \eta\nabla F(z_t),\ \eta\lambda_1\big),\qquad
s_t = \frac{1 + \sqrt{1 + 4s_{t-1}^2}}{2},\ s_0 = 1 .
\]
The gradient, the thresholding and the objective are unchanged, and each iteration still costs one product $X^T(Xz)$. The momentum weight $s_t$ is a single scalar shared by every configuration, so FISTA batches over configurations exactly like the plain method: the \textit{Flag Approach} or Dynamic Shrinking, chunking and the tensor layout carry over without change. The two methods differ in four ways.
\begin{itemize}
\item \textbf{Rate.} The objective error falls like $O(1/t)$ for ISTA and like $O(1/t^2)$ for FISTA. On strongly correlated designs this difference is the one that matters. On Air Quality ($\sigma^2_{\max}(X) = 7{,}089$ with unit-norm columns), a pilot on the $10\times10$ grid gave a glmnet-to-ParaDeAl train RE ratio of 0.915 for plain proximal gradient after 1{,}000 iterations and 0.990 after 10{,}000, while FISTA matched glmnet to four decimals after 5{,}000.
\item \textbf{Step size.} The guarantee of FISTA needs $\eta = 1/L$ with $L = 2(\sigma^2_{\max}(X) + \lambda_2)$; the plain method may go up to $2/L$ and is run at $0.9 \cdot 2/L$. The steps of FISTA are therefore about half as long, and the momentum makes up for the shorter steps many times over.
\item \textbf{Memory.} One extra $p$-vector per configuration, the previous iterate.
\item \textbf{Monotonicity.} The objective may rise briefly because of the overshoot, whereas plain proximal gradient decreases the objective at every step.
\end{itemize}
FISTA is used only in the L1+L2 mode, where there is no hard thresholding. Momentum combined with $H_k$ has no standard convergence guarantee, so the L0 modes use the method of Section~\ref{sec:architecture} unchanged.

\paragraph{Alternatives tested and not adopted.} All of the following were run with Dynamic Shrinking on L1+L2.
\begin{itemize}
\item \emph{An active set for L1+L2.} Each configuration starts from the (at most 250) features that violate optimality at $\beta = 0$, solves on its block, then checks optimality at full width and adds the violators, until none remain, which certifies the full solution. The procedure is exact but slower than full width: 60\,s against 32\,s on Crime $30\times30$ (FISTA, tolerance $10^{-4}$), and 858\,s against glmnet's 570\,s on Diabetes $30\times30$. The supports are too large for the procedure to pay off. On Diabetes the block grew to 25{,}250 columns and some configurations were still adding features after 100 rounds, and growing the block by five times the support size instead of by 250 per round ran out of GPU memory. With supports in the thousands there is little to screen out, and every round pays for a full-width optimality check and restarts the momentum of FISTA.
\item \emph{FISTA with adaptive restart} \citep{odo15}, which resets the momentum of a configuration when the momentum points against the step. This variant is faster (5.6$\times$ glmnet instead of 4.0$\times$ on Air Quality $100\times100$) but less accurate (train RE 0.9\% and test RE 3\% worse than glmnet): after a restart the step is short, so the relative-change rule stops the configuration too early.
\item \emph{One chunk instead of several} on Air Quality $100\times100$: no measurable difference.
\end{itemize}

\paragraph{Measurements.} Accuracy is computed by one Python routine for every method, from the returned coefficients. Runtime and storage are measured in each method's own language, since R and Python time differently. R uses \texttt{system.time} (elapsed) around the solver calls. Python uses \texttt{perf\_counter}, with a CUDA synchronization before each reading for ParaDeAl and after one untimed warm-up. ParaDeAl's time includes everything from the configurations to the final coefficients ($\sigma_{\max}$, the step-size search and both stages). Storage is peak memory above the resident data: GPU memory for ParaDeAl and process resident memory for the CPU methods, with the peak counter reset right before the solve. scikit-learn and skglm fit one configuration at a time, so on every grid a fixed random subset of at most 50 configurations is timed and scored, and the full-grid time is the per-configuration average times the grid size.

\section{Results on the Polynomial Datasets} \label{app:poly_results}

Table~\ref{tab:res_soft} repeats the comparison of Section~\ref{subsec:real_results} on the three polynomial datasets of Section~\ref{subsec:datasets}, which are built without neighbor features and whose correlations are weaker. It supports the conclusions drawn there.

\begin{table}[p]\centering\scriptsize\setlength{\tabcolsep}{3.5pt}
{\footnotesize (a) Runtime, and accuracy over all configurations}\\[3pt]
\begin{tabular}{llrrrcc}\toprule
 & & \multicolumn{3}{c}{Runtime (s) (speedup)} & \multicolumn{2}{c}{Competitor $\div$ ParaDeAl} \\ \cmidrule(lr){3-5}\cmidrule(lr){6-7}
Mode & Method & 10$\times$10 & 30$\times$30 & 100$\times$100 & Objective ($\uparrow$) & Test error ($\downarrow$) \\
\midrule
\multicolumn{5}{l}{\textbf{Communities and Crime} ($p = 302{,}743$)} & $n = 1{,}500$ & $n = 250$ \\
L0+L1 & ParaDeAl (pair swaps) & 0.365 & 2.30 & 429 & 1 & 1 \\
 & ParaDeAl (no swaps) & 0.344 & 1.82 & 87.0 & 0.638 & 1.33 \\
 & L0Learn (CPU) & 26.8 (73.3$\times$) & 263 (114$\times$) & 2{,}101 (4.89$\times$) & 0.947 & 1.07 \\
L0+L2 & ParaDeAl (pair swaps) & 1.42 & 25.1 & 348 & 1 & 1 \\
 & ParaDeAl (no swaps) & 1.35 & 24.6 & 242 & 0.655 & 1.15 \\
 & L0Learn (CPU) & 29.4 (20.7$\times$) & 272 (10.9$\times$) & 2{,}638 (7.59$\times$) & 0.927 & 1.04 \\
L1+L2 & ParaDeAl (FISTA) & 10.2 & 46.8 & 488 & 1 & 1 \\
 & glmnet (CPU) & 61.2 (6.00$\times$) & 223 (4.77$\times$) & 1{,}537 (3.15$\times$) & 0.999 & 1.00 \\
 & scikit-learn (CPU) & 11{,}396 (1{,}117$\times$) & 85{,}545 (1{,}829$\times$) & 683{,}923 (1{,}400$\times$) & 0.999 & 1.00 \\
 & skglm (CPU) & 235 (23.0$\times$) & 2{,}341 (50.0$\times$) & 25{,}717 (52.7$\times$) & 1.016 & 1.00 \\
L0+L1+L2 & ParaDeAl & -- & -- & 502 & -- & -- \\
\midrule
\multicolumn{5}{l}{\textbf{NHANES} ($p = 23{,}478$)} & $n = 2{,}500$ & $n = 1{,}000$ \\
L0+L1 & ParaDeAl (pair swaps) & 0.203 & 2.12 & 330 & 1 & 1 \\
 & ParaDeAl (no swaps) & 0.137 & 0.728 & 22.2 & 0.835 & 1.12 \\
 & L0Learn (CPU) & 2.96 (14.6$\times$) & 33.0 (15.6$\times$) & 326 (0.988$\times$) & 0.984 & 1.00 \\
L0+L2 & ParaDeAl (pair swaps) & 0.246 & 2.65 & 77.8 & 1 & 1 \\
 & ParaDeAl (no swaps) & 0.201 & 2.25 & 30.2 & 0.856 & 1.10 \\
 & L0Learn (CPU) & 2.47 (10.0$\times$) & 44.5 (16.8$\times$) & 402 (5.16$\times$) & 0.986 & 1.00 \\
L1+L2 & ParaDeAl (FISTA) & 0.428 & 1.45 & 15.0 & 1 & 1 \\
 & glmnet (CPU) & 14.5 (33.9$\times$) & 74.7 (51.7$\times$) & 474 (31.5$\times$) & 1.000 & 1.00 \\
 & scikit-learn (CPU) & 319 (744$\times$) & 2{,}064 (1{,}427$\times$) & 17{,}658 (1{,}174$\times$) & 1.000 & 1.00 \\
 & skglm (CPU) & 61.6 (144$\times$) & 462 (320$\times$) & 4{,}817 (320$\times$) & 1.000 & 1.00 \\
L0+L1+L2 & ParaDeAl & -- & -- & 52.5 & -- & -- \\
\midrule
\multicolumn{5}{l}{\textbf{Weather} ($p = 100{,}000$)} & $n = 2{,}500$ & $n = 1{,}123$ \\
L0+L1 & ParaDeAl (pair swaps) & 0.177 & 0.782 & 427 & 1 & 1 \\
 & ParaDeAl (no swaps) & 0.160 & 0.640 & 55.5 & 0.718 & 1.05 \\
 & L0Learn (CPU) & 47.5 (268$\times$) & 111 (141$\times$) & 1{,}514 (3.54$\times$) & 0.980 & 1.03 \\
L0+L2 & ParaDeAl (pair swaps) & 1.26 & 10.0 & 179 & 1 & 1 \\
 & ParaDeAl (no swaps) & 1.21 & 9.74 & 108 & 0.732 & 1.00 \\
 & L0Learn (CPU) & 16.7 (13.3$\times$) & 111 (11.0$\times$) & 1{,}369 (7.63$\times$) & 0.962 & 1.01 \\
L1+L2 & ParaDeAl (FISTA) & 2.58 & 11.0 & 115 & 1 & 1 \\
 & glmnet (CPU) & 33.7 (13.1$\times$) & 212 (19.3$\times$) & 1{,}159 (10.1$\times$) & 1.000 & 1.00 \\
 & scikit-learn (CPU) & 2{,}186 (846$\times$) & 16{,}944 (1{,}538$\times$) & 143{,}669 (1{,}247$\times$) & 1.000 & 1.00 \\
 & skglm (CPU) & 126 (48.9$\times$) & 1{,}586 (144$\times$) & 14{,}496 (126$\times$) & 1.000 & 1.00 \\
L0+L1+L2 & ParaDeAl & -- & -- & 210 & -- & -- \\
\bottomrule\end{tabular}

\vspace{4pt}
{\footnotesize (b) Accuracy by support size $k$ in the L0 modes, pooled over the three grids}\\[2pt]
\begin{tabular}{llcccccccccccc}\toprule
 & & \multicolumn{6}{c}{Objective (L0Learn $\div$ ParaDeAl)} & \multicolumn{6}{c}{Test error (L0Learn $\div$ ParaDeAl)} \\ \cmidrule(lr){3-8}\cmidrule(lr){9-14}
Mode & & 1--5 & 6--20 & 21--50 & 51--100 & $>$100 & All & 1--5 & 6--20 & 21--50 & 51--100 & $>$100 & All \\
\midrule
\multicolumn{14}{l}{\textbf{Communities and Crime, products of up to 3}} \\
L0+L1 & Configurations & 545 & 608 & 409 & 382 & 1489 & 3433 & & & & & & \\
 & ParaDeAl & 0.991 & 0.991 & 0.984 & 0.968 & 0.897 & 0.947 & 1.01 & 1.07 & 1.04 & 1.05 & 1.11 & 1.07 \\
 & No pair swaps & 0.938 & 0.860 & 0.777 & 0.662 & 0.394 & 0.638 & 0.94 & 0.97 & 1.02 & 1.18 & 1.75 & 1.33 \\
L0+L2 & Configurations & 560 & 685 & 528 & 587 & 2398 & 4758 & & & & & & \\
 & ParaDeAl & 0.993 & 0.988 & 0.978 & 0.963 & 0.873 & 0.927 & 1.00 & 1.04 & 1.01 & 1.02 & 1.05 & 1.04 \\
 & No pair swaps & 0.930 & 0.882 & 0.786 & 0.703 & 0.486 & 0.655 & 0.92 & 0.94 & 0.95 & 1.02 & 1.34 & 1.15 \\
\midrule
\multicolumn{14}{l}{\textbf{NHANES}} \\
L0+L1 & Configurations & 510 & 659 & 373 & 375 & 1414 & 3331 & & & & & & \\
 & ParaDeAl & 1.001 & 1.002 & 0.998 & 0.994 & 0.964 & 0.984 & 1.00 & 1.01 & 1.00 & 0.99 & 1.00 & 1.00 \\
 & No pair swaps & 0.972 & 0.911 & 0.922 & 0.860 & 0.720 & 0.835 & 0.97 & 0.95 & 1.03 & 1.07 & 1.29 & 1.12 \\
L0+L2 & Configurations & 544 & 707 & 484 & 577 & 2196 & 4508 & & & & & & \\
 & ParaDeAl & 1.001 & 1.001 & 0.999 & 0.995 & 0.972 & 0.986 & 1.00 & 1.00 & 1.00 & 0.99 & 1.00 & 1.00 \\
 & No pair swaps & 0.983 & 0.926 & 0.949 & 0.882 & 0.775 & 0.856 & 0.98 & 0.95 & 1.01 & 1.05 & 1.20 & 1.10 \\
\midrule
\multicolumn{14}{l}{\textbf{Weather}} \\
L0+L1 & Configurations & 525 & 412 & 368 & 343 & 1290 & 2938 & & & & & & \\
 & ParaDeAl & 1.009 & 1.011 & 1.005 & 0.997 & 0.946 & 0.980 & 1.00 & 1.00 & 1.02 & 1.03 & 1.04 & 1.03 \\
 & No pair swaps & 0.990 & 0.906 & 0.803 & 0.696 & 0.530 & 0.718 & 0.97 & 0.92 & 0.89 & 0.93 & 1.19 & 1.05 \\
L0+L2 & Configurations & 543 & 624 & 532 & 578 & 2359 & 4636 & & & & & & \\
 & ParaDeAl & 1.006 & 1.008 & 0.994 & 0.982 & 0.928 & 0.962 & 1.00 & 0.99 & 1.01 & 1.02 & 1.02 & 1.01 \\
 & No pair swaps & 0.990 & 0.924 & 0.795 & 0.737 & 0.607 & 0.732 & 0.97 & 0.93 & 0.88 & 0.88 & 1.08 & 1.00 \\
\bottomrule\end{tabular}

\caption{Results on the three polynomial datasets, built without neighbor features (Section~\ref{subsec:datasets}). Panel (a) and the settings are as in Table~\ref{tab:res_hard}. Panel (b) gives the ratios of Table~\ref{tab:res_hard}(b) as numbers: ``ParaDeAl'' is ParaDeAl with pair swaps, and ``Configurations'' gives the number of configurations of each size.}
\label{tab:res_soft}
\end{table}

\section{Trying Several Supports Before Value Modification} \label{app:multisupport}

Section~\ref{subsec:two_stage_desc} notes that the support chosen by support selection is not always the support with the smallest objective after value modification. This appendix measures how often that happens and how much is gained by trying several supports instead of one.

\paragraph{Setup.} For each configuration we collect about 15 candidate supports from three sources, refit every candidate with ParaDeAl's own value modification, and keep the one with the smallest objective. (i) The last five distinct supports visited during support selection, together with the support of the lowest-objective iterate. (ii) Four one-for-one exchanges at the edge of the hard threshold: the weakest feature in the support is replaced by one of the strongest features just outside it. (iii) Support selection rerun with the starting step size, a quarter and a sixteenth of the chosen step size, and from a random starting point. The chosen support is the one returned by Algorithm~\ref{alg:support_selection} with the last-iterate rule. We use $30\times30$ grids in both L0 modes, on the three real datasets of Section~\ref{sec:compare_results} and on synthetic data ($n = 1{,}000$, $p = 20{,}000$, 50 true features, SNR 5) with either independent features or columns whose correlation is $0.5^{|i-j|}$, as in Appendix~\ref{app:stability_r}. Pair swaps (Section~\ref{subsec:pair_swap}) are run from the same chosen support for comparison.

\begin{table}[htbp]
\centering\footnotesize\setlength{\tabcolsep}{4pt}
\begin{tabular}{l|ccccc|ccccc}
\toprule
 & \multicolumn{5}{c|}{L0+L1} & \multicolumn{5}{c}{L0+L2} \\
Data & Chosen & Bo15 & Pair swap & Gap closed & Not best & Chosen & Bo15 & Pair swap & Gap closed & Not best \\
\midrule
Air Quality & 0.233 & 0.710 & 0.991 & 63\% & 99\% & 0.603 & 0.842 & 0.976 & 64\% & 94\% \\
Communities and Crime & 0.933 & 0.964 & 1.003 & 44\% & 66\% & 0.954 & 0.974 & 0.999 & 44\% & 80\% \\
Diabetes & 0.872 & 0.908 & 0.983 & 32\% & 88\% & 0.910 & 0.935 & 0.992 & 30\% & 84\% \\
Synthetic, independent & 0.999 & 1.000 & 1.000 & $\approx$100\% & 23\% & 0.995 & 0.999 & 1.000 & 83\% & 55\% \\
Synthetic, correlated & 0.942 & 0.989 & 1.000 & 81\% & 32\% & 0.991 & 0.997 & 1.000 & 67\% & 56\% \\
\bottomrule
\end{tabular}
\caption{Trying about 15 supports per configuration. The first three columns of each mode give the optimization objective of L0Learn divided by that of ParaDeAl, averaged over configurations (`Bo15' means `Best of 15', and 1 means the same objective; above 1, ParaDeAl is better), for the chosen support, the best of the candidates, and pair swaps. ``Gap closed'' is the share of pair swap's improvement over the chosen support that the best candidate recovers. ``Not best'' is the share of configurations for which some candidate beats the chosen support.}
\label{tab:multisupport}
\end{table}

\paragraph{Results.} Table~\ref{tab:multisupport} answers both questions. The chosen support is often not the best one: on the real datasets, some candidate beats it for 66 to 99\% of configurations. The gain from picking the best candidate is real but partial. It recovers 30 to 64\% of the improvement that pair swaps achieve on the real datasets, and most of it comes from the exchanges at the edge of the hard threshold, at small $k$ (between 2 and 20). With independent features the chosen support is already within 0.5\% of L0Learn, and the gaps are tiny. On five further real datasets with weaker correlations, the picture is the same: a better candidate exists for 50 to 89\% of configurations, and the best candidate recovers 24 to 70\% of the pair swap improvement.

The candidates are also not cheaper than pair swaps. Refitting about 15 supports per configuration added 0.3 to 3.9\,s per grid, often more than the whole solve without swaps, whereas pair swaps added 0.07 to 2.2\,s and reached an objective at least as good in every case. Trying several supports can be parallelized in the same way as everything else in ParaDeAl, but pair swaps are the better use of the extra computation, which is why they, and not multiple supports, are part of ParaDeAl.

\section{Pair Swaps on Correlated Synthetic Data} \label{app:pair_swap_results}

This appendix gives the full results behind the paragraph ``Effect'' in Section~\ref{subsec:pair_swap}. The data are synthetic, with $n = 1{,}000$, $p = 20{,}000$, 50 true features and SNR 5, generated as in Appendix~\ref{app:stability_r}, except that features $i$ and $j$ have correlation $\rho^{|i-j|}$ with $\rho \in \{0, 0.25, 0.5, 0.9, 0.99\}$. At $\rho = 0.99$ each feature is nearly a copy of its neighbors, which is the hardest case for support selection. For each mode we solve the $30 \times 30$ grid of configurations chosen by L0Learn, and compare three versions of ParaDeAl: the last support of support selection, the lowest-objective support (Section~\ref{subsec:lowest_objective}), and the lowest-objective support followed by pair swaps. Times are the median of three runs on one L40S GPU; L0Learn runs on one CPU core.

\begin{table}[htbp]
\centering\footnotesize\setlength{\tabcolsep}{4pt}
\begin{tabular}{cl|ccc|c|cc|cc}
\toprule
 & & \multicolumn{3}{c|}{Objective ratio} & Test ratio & \multicolumn{2}{c|}{Seconds} & \multicolumn{2}{c}{True features found} \\
$\rho$ & Mode & Last & Lowest & + Pair swap & + Pair swap & Pair swap adds & L0Learn & L0Learn & + Pair swap \\
\midrule
0 & L0+L1 & 0.999 & 1.000 & 1.000 & 1.000 & 0.06 & 6.7 & 0.68 & 0.68 \\
 & L0+L2 & 0.995 & 0.997 & 1.000 & 1.000 & 0.10 & 7.4 & 0.57 & 0.55 \\
0.25 & L0+L1 & 0.999 & 1.000 & 1.000 & 1.001 & 0.07 & 6.8 & 0.64 & 0.64 \\
 & L0+L2 & 0.996 & 0.998 & 1.000 & 1.000 & 0.16 & 6.9 & 0.51 & 0.51 \\
0.5 & L0+L1 & 0.942 & 0.985 & 1.000 & 1.000 & 0.11 & 8.5 & 0.62 & 0.62 \\
 & L0+L2 & 0.991 & 0.992 & 1.000 & 1.000 & 0.08 & 8.7 & 0.46 & 0.46 \\
0.9 & L0+L1 & 0.714 & 0.726 & 1.002 & 1.017 & 0.32 & 5.7 & 0.56 & 0.52 \\
 & L0+L2 & 0.794 & 0.794 & 1.000 & 1.004 & 0.13 & 6.3 & 0.34 & 0.34 \\
0.99 & L0+L1 & 0.616 & 0.620 & 0.986 & 0.964 & 0.71 & 7.4 & 0.10 & 0.17 \\
 & L0+L2 & 0.689 & 0.692 & 0.999 & 0.968 & 0.25 & 6.2 & 0.16 & 0.16 \\
\bottomrule
\end{tabular}
\caption{Pair swaps on synthetic data with correlation $\rho^{|i-j|}$. The objective and test ratios are L0Learn's value divided by ParaDeAl's, averaged over configurations, so 1 is a tie and a value below 1 means L0Learn is better. ``Last'', ``Lowest'' and ``+ Pair swap'' are the three versions of ParaDeAl described in the text. ``Pair swap adds'' is the time of the pair swaps alone, per grid. ``True features found'' is the share of the 50 true features in the selected support, averaged over the configurations with $k = 50$ (or the nearest $k$ on L0Learn's path).}
\label{tab:pair_swap_synthetic}
\end{table}

\paragraph{Results.} Table~\ref{tab:pair_swap_synthetic} shows that the weakness of ParaDeAl without pair swaps is a weakness under correlation. With independent features or $\rho = 0.25$, ParaDeAl already matches L0Learn to within 0.4\%, and with $\rho = 0.5$ it falls short only slightly. At $\rho = 0.9$ and $0.99$ it reaches only 62 to 79\% of L0Learn's objective, and the shortfall is concentrated at moderate support sizes: at $\rho = 0.99$ in the L0+L1 mode, the mean ratio is 0.49 for $k$ between 6 and 20 and 0.33 for $k$ between 21 and 50. Choosing the lowest-objective support helps little here, since the problem is not which iterate is returned but which support is found. Pair swaps remove the shortfall: the ratio becomes 0.99 to 1.00 for every $\rho$, and 90\% of the configurations lie within 0.98 to 1.03 at $\rho = 0.9$ and within 0.82 to 1.06 at $\rho = 0.99$. The test error follows the objective, and pair swaps find about as many of the true features as L0Learn; over the three random datasets they find more in seven of the twelve comparisons, the same number in three and fewer in two. Few of the true features are found by any method at $\rho = 0.99$, because each true feature has near-copies that fit almost equally well.

Pair swaps are cheap. A configuration needs 4 to 7 rounds on average at $\rho = 0.9$ and $0.99$, and at most 31; no configuration in any run was still improving at the limit of 200 rounds. With pair swaps, ParaDeAl solves each grid in 0.2 to 1.7 seconds, against 5.7 to 8.7 seconds for L0Learn.

\paragraph{By support size.} Table~\ref{tab:pair_swap_synthetic_k} splits the objective ratios by support size. Without pair swaps, the shortfall under strong correlation grows with $k$ up to about 50; pair swaps bring every bin back to between 0.977 and 1.014. Supports larger than 100 occur in at most 8 configurations per grid and change the means by at most 0.003, so they are not shown.

\begin{table}[htbp]
\centering\footnotesize\setlength{\tabcolsep}{5pt}
\begin{tabular}{clcccc}
\toprule
$\rho$ & Mode & $k$: 1--5 & 6--20 & 21--50 & 51--100 \\
\midrule
0 & L0+L1 & 1.000 / 1.000 & 0.999 / 1.000 & 0.998 / 1.000 & -- \\
 & L0+L2 & 0.997 / 1.000 & 0.993 / 1.000 & 0.996 / 1.000 & 0.996 / 1.000 \\
0.25 & L0+L1 & 1.000 / 1.000 & 0.999 / 1.000 & 0.999 / 1.000 & -- \\
 & L0+L2 & 0.998 / 1.000 & 0.996 / 1.000 & 0.996 / 1.000 & 0.994 / 1.000 \\
0.5 & L0+L1 & 1.000 / 1.000 & 0.999 / 1.000 & 0.771 / 1.000 & -- \\
 & L0+L2 & 0.999 / 1.000 & 0.998 / 1.000 & 0.977 / 1.000 & 0.997 / 1.000 \\
0.9 & L0+L1 & 0.867 / 0.999 & 0.734 / 0.995 & 0.503 / 1.014 & -- \\
 & L0+L2 & 0.906 / 0.999 & 0.784 / 0.997 & 0.670 / 1.004 & 0.990 / 1.000 \\
0.99 & L0+L1 & 0.875 / 1.000 & 0.492 / 0.977 & 0.326 / 0.997 & -- \\
 & L0+L2 & 0.895 / 0.999 & 0.557 / 0.996 & 0.633 / 1.004 & 0.884 / 1.001 \\
\bottomrule
\end{tabular}
\caption{Objective ratio (L0Learn divided by ParaDeAl) on the synthetic data of Table~\ref{tab:pair_swap_synthetic}, by support size $k$. In each cell, the first number is for the last support without pair swaps, and the second for the lowest-objective support with pair swaps. A dash means that the grid has no configuration of that size.}
\label{tab:pair_swap_synthetic_k}
\end{table}

\paragraph{Other random datasets.} We repeated the runs with $\rho \in \{0.9, 0.99\}$ on two further random datasets. Without pair swaps the mean objective ratio is between 0.61 and 0.80 on all six datasets, and with pair swaps it is between 0.986 and 1.029, so the results above do not depend on one particular draw.

\section{Returning the Last Versus the Lowest-Objective Support} \label{app:lowest_objective}

This appendix compares the two ways of choosing which iterate support selection returns (Section~\ref{subsec:lowest_objective}): the last iterate, and the iterate with the lowest objective. No pair swaps are used, so that the effect of the rule is seen on its own. The grids and settings are those of Appendix~\ref{app:pair_swap_results}, on the six real datasets of Section~\ref{sec:compare_results} and on the synthetic data of Appendix~\ref{app:pair_swap_results}.

\begin{table}[htbp]
\centering\footnotesize\setlength{\tabcolsep}{3.5pt}
\begin{tabular}{l|cccc|cccc}
\toprule
 & \multicolumn{4}{c|}{L0+L1} & \multicolumn{4}{c}{L0+L2} \\
 & \multicolumn{2}{c}{All $k$} & \multicolumn{2}{c|}{$k = 1$} & \multicolumn{2}{c}{All $k$} & \multicolumn{2}{c}{$k = 1$} \\
Data & Last & Lowest & Last & Lowest & Last & Lowest & Last & Lowest \\
\midrule
Diabetes & 0.872 & 0.872 & 1.000 & 1.000 & 0.910 & 0.910 & 1.000 & 1.000 \\
Communities and Crime & 0.933 & 0.942 & 1.000 & 1.000 & 0.954 & 0.954 & 1.000 & 1.000 \\
Air Quality & 0.233 & 0.241 & 0.959 & 1.000 & 0.603 & 0.610 & 0.996 & 1.000 \\
\midrule
Crime, products of up to 3 & 0.869 & 0.898 & 1.000 & 1.000 & 0.915 & 0.915 & 1.000 & 0.999 \\
NHANES & 0.879 & 0.918 & 0.871 & 1.000 & 0.940 & 0.940 & 0.998 & 1.000 \\
Weather & 0.887 & 0.974 & 1.000 & 1.000 & 0.970 & 0.969 & 1.000 & 0.998 \\
\midrule
Synthetic, $\rho = 0$ & 0.999 & 1.000 & 1.000 & 1.000 & 0.995 & 0.997 & 0.996 & 1.000 \\
Synthetic, $\rho = 0.5$ & 0.942 & 0.985 & 1.000 & 1.000 & 0.991 & 0.992 & 0.999 & 1.000 \\
Synthetic, $\rho = 0.9$ & 0.714 & 0.726 & 1.000 & 1.000 & 0.794 & 0.794 & 0.995 & 1.000 \\
Synthetic, $\rho = 0.99$ & 0.616 & 0.620 & 0.980 & 1.000 & 0.689 & 0.692 & 0.994 & 1.000 \\
\bottomrule
\end{tabular}
\caption{The last support versus the lowest-objective support, without pair swaps. Each entry is the optimization objective of L0Learn divided by that of ParaDeAl, averaged over all configurations of the grid or over the configurations with $k = 1$ only; 1 is a tie and a value below 1 means L0Learn is better. Synthetic data use the first of the random datasets of Appendix~\ref{app:pair_swap_results}.}
\label{tab:lowest_objective}
\end{table}

\paragraph{Results.} The lowest-objective rule is never meaningfully worse than the last-iterate rule: the largest loss in Table~\ref{tab:lowest_objective} is 0.001. It helps most in the L0+L1 mode, on Weather (0.887 to 0.974), the synthetic data with $\rho = 0.5$ (0.942 to 0.985), NHANES (0.879 to 0.918) and Crime with products of up to 3 (0.869 to 0.898); elsewhere the gain is between 0 and 1\%. The gain is largest, and most systematic, at $k = 1$, the case described in Section~\ref{subsec:lowest_objective}. On NHANES in the L0+L1 mode, 24 of the 29 configurations with $k = 1$ end more than 1\% behind L0Learn under the last-iterate rule, and on the synthetic data with $\rho = 0.99$, 27 of 29 (and 22 of 29 on a second random dataset). The lowest-objective rule fixes every one of these configurations, and under it the mean ratio at $k = 1$ is between 0.997 and 1.000 on every dataset and mode.

The two rules return a different support for between 0 and 40\% of the configurations (most often on Air Quality). The lowest-objective support is better after value modification in almost all of these cases, and worse in at most 4\% of all configurations. Since value modification can reorder two supports (Section~\ref{subsec:two_stage_desc} and Appendix~\ref{app:multisupport}), a lower objective during support selection does not guarantee a lower objective at the end, but it does so in most cases. The rule adds no measurable runtime: the differences between the two rules are within the variation between repeated runs.

\section{Agreement with L0Learn's Supports} \label{app:support_agreement}

Section~\ref{subsec:real_results} compares ParaDeAl and L0Learn through their objectives and test errors. This appendix compares the supports themselves, on the configurations of Section~\ref{sec:compare_results}, pooled over the three grids. For each configuration we record whether the two supports are identical, and their overlap $|A \cap B| / |A \cup B|$, where $A$ and $B$ are the supports returned by ParaDeAl and L0Learn. Table~\ref{tab:support_agreement} gives both by support size, for ParaDeAl with the lowest-objective support (Section~\ref{subsec:lowest_objective}) and with pair swaps added (Section~\ref{subsec:pair_swap}).

\begin{table}[htbp]
\centering\footnotesize\setlength{\tabcolsep}{4pt}
{(a) Share of configurations with identical supports (\%), without / with pair swaps}\\[2pt]
\begin{tabular}{llcccccc}
\toprule
Data & Mode & 1--5 & 6--20 & 21--50 & 51--100 & $>$100 & All \\ \midrule
Diabetes & L0+L1 & 0 / 0 & 0 / 0 & 0 / 0 & 0 / 0 & 0 / 0 & 0 / 0 \\
 & L0+L2 & 3 / 3 & 2 / 3 & 0 / 0 & 0 / 0 & 0 / 0 & 1 / 1 \\
Crime & L0+L1 & 25 / 32 & 0 / 10 & 0 / 4 & 0 / 0 & 0 / 0 & 5 / 8 \\
 & L0+L2 & 30 / 43 & 0 / 3 & 0 / 3 & 0 / 1 & 0 / 0 & 4 / 7 \\
Air Quality & L0+L1 & 4 / 35 & 0 / 0 & 0 / 0 & 0 / 0 & 0 / 0 & 1 / 7 \\
 & L0+L2 & 11 / 43 & 0 / 7 & 2 / 9 & 0 / 1 & 0 / 0 & 2 / 7 \\
Crime, products & L0+L1 & 25 / 40 & 0 / 3 & 0 / 1 & 0 / 0 & 0 / 0 & 4 / 7 \\
 & L0+L2 & 27 / 48 & 1 / 9 & 0 / 1 & 0 / 0 & 0 / 0 & 3 / 7 \\
NHANES & L0+L1 & 35 / 76 & 0 / 18 & 0 / 5 & 0 / 2 & 0 / 0 & 5 / 16 \\
 & L0+L2 & 46 / 83 & 1 / 20 & 1 / 24 & 0 / 9 & 0 / 2 & 6 / 18 \\
Weather & L0+L1 & 36 / 62 & 0 / 21 & 0 / 11 & 0 / 4 & 0 / 0 & 6 / 16 \\
 & L0+L2 & 44 / 71 & 11 / 26 & 0 / 3 & 0 / 1 & 0 / 0 & 7 / 12 \\
\bottomrule
\end{tabular}

\vspace{6pt}
{(b) Mean overlap $|A \cap B| / |A \cup B|$, without / with pair swaps}\\[2pt]
\begin{tabular}{llcccccc}
\toprule
Data & Mode & 1--5 & 6--20 & 21--50 & 51--100 & $>$100 & All \\ \midrule
Diabetes & L0+L1 & 0.02 / 0.05 & 0.05 / 0.15 & 0.02 / 0.29 & 0.01 / 0.26 & 0.01 / 0.20 & 0.02 / 0.20 \\
 & L0+L2 & 0.10 / 0.12 & 0.27 / 0.34 & 0.18 / 0.39 & 0.12 / 0.34 & 0.04 / 0.27 & 0.10 / 0.28 \\
Crime & L0+L1 & 0.48 / 0.56 & 0.16 / 0.41 & 0.06 / 0.35 & 0.05 / 0.35 & 0.02 / 0.22 & 0.14 / 0.34 \\
 & L0+L2 & 0.50 / 0.66 & 0.26 / 0.42 & 0.23 / 0.45 & 0.22 / 0.47 & 0.14 / 0.45 & 0.23 / 0.48 \\
Air Quality & L0+L1 & 0.26 / 0.67 & 0.08 / 0.36 & 0.06 / 0.33 & 0.03 / 0.18 & 0.02 / 0.09 & 0.09 / 0.33 \\
 & L0+L2 & 0.35 / 0.69 & 0.21 / 0.43 & 0.40 / 0.51 & 0.36 / 0.50 & 0.12 / 0.38 & 0.23 / 0.45 \\
Crime, products & L0+L1 & 0.44 / 0.56 & 0.11 / 0.22 & 0.04 / 0.26 & 0.03 / 0.25 & 0.01 / 0.17 & 0.10 / 0.26 \\
 & L0+L2 & 0.45 / 0.65 & 0.25 / 0.35 & 0.14 / 0.38 & 0.09 / 0.42 & 0.04 / 0.38 & 0.13 / 0.41 \\
NHANES & L0+L1 & 0.50 / 0.89 & 0.15 / 0.68 & 0.23 / 0.47 & 0.18 / 0.43 & 0.15 / 0.33 & 0.22 / 0.51 \\
 & L0+L2 & 0.64 / 0.92 & 0.26 / 0.69 & 0.40 / 0.67 & 0.26 / 0.56 & 0.22 / 0.50 & 0.30 / 0.61 \\
Weather & L0+L1 & 0.61 / 0.79 & 0.27 / 0.59 & 0.11 / 0.48 & 0.05 / 0.38 & 0.04 / 0.26 & 0.19 / 0.44 \\
 & L0+L2 & 0.65 / 0.86 & 0.43 / 0.69 & 0.18 / 0.51 & 0.11 / 0.57 & 0.12 / 0.46 & 0.23 / 0.56 \\
\bottomrule
\end{tabular}
\caption{Agreement between the supports of ParaDeAl and L0Learn, by support size $k$. In each cell, the first number is for ParaDeAl with the lowest-objective support and the second for ParaDeAl with pair swaps added. ``Crime, products'' is Communities and Crime with products of up to 3 features.}
\label{tab:support_agreement}
\end{table}

\paragraph{Results.} Identical supports are common only for the smallest supports. With pair swaps, ParaDeAl returns L0Learn's support for 62 to 83\% of the configurations with $k$ from 1 to 5 on NHANES and Weather, for 32 to 48\% on the two Crime datasets and Air Quality, and almost never on Diabetes. For $k$ above 20, at most a quarter of the supports are identical, and on most datasets almost none. The overlap tells the same story more gradually: pair swaps raise it on every dataset and in every bin, but even with pair swaps, two supports of the same size typically share only between a fifth and two thirds of their features. Since the objectives are nearly equal (Tables~\ref{tab:res_hard} and~\ref{tab:res_soft}), this means that the two methods often find different supports of almost the same quality, as expected when many features are close copies of one another. The agreement is highest on NHANES and Weather, whose features are the least correlated.

\end{appendices}

\end{document}